\documentclass[a4paper,11pt]{amsart}

\usepackage{geometry} 
\usepackage{amssymb, amsmath, mathrsfs} 
\usepackage{comment} 
\usepackage{hyperref} 
\usepackage[utf8]{inputenc} 
\usepackage[T1]{fontenc} 
\usepackage{enumerate} 
\usepackage{tcolorbox} 
\usepackage{color} 
\usepackage{srcltx} 
\usepackage[mathcal]{euscript} 
\usepackage{amsthm, amssymb} 
\usepackage{setspace} 
\usepackage{xfrac} 
\usepackage{bbold} 
\usepackage{tikz-cd} 
\usepackage{ifpdf} 
\usepackage{IEEEtrantools} 
\usepackage{todonotes} 
\usepackage{url} 
\usetikzlibrary{calc} 
\usepackage{bbm} 

\newcommand{\N}{\mathbb{N}} 
 
\newcommand{\RR}{\mathrm{R}}
\newcommand{\CC}{\mathrm{C}} 

\newcommand{\dist}{\mathrm{dist}}
\newcommand{\Hau}{\mathcal{H}} 
\newcommand{\Leb}{\mathcal{L}} 
\renewcommand{\d}{\,\mathrm{d}} 
\newcommand{\dd}{\d} 
\newcommand{\D}{\mathrm{D}}
\newcommand{\WW}{\mathrm{W}}
\newcommand{\LL}{\mathrm{L}}
\newcommand{\R}{\mathbb R}
\newcommand{\Tr}{\operatorname{tr}}
\newcommand{\tr}{\operatorname{tr}}
\newcommand{\BV}{\mathrm{BV}}

\newcommand{\B}{\mathbb B}

\newcommand{\weaklystar}{\text{weakly${}^\ast$ }}
\newcommand{\weakstar}{\text{weak${}^\ast$ }}
\newcommand{\A}{\mathbb A}
\newcommand{\Trace}{\mathrm{tr}}

\newcommand{\Esgnu}{\mathrm{E}_{\mathrm{sg}}^{1,2}}
\newcommand{\Eren}{\mathrm{E}_{\mathrm{ren}}}
\newcommand{\Erennu}{\mathrm{E}_{\mathrm{ren}}^{1,2}}
\newcommand{\sysN}{\mathrm{sys}(\mathcal N)}
\newcommand{\ds}{\displaystyle}
\renewcommand{\emptyset}{\varnothing}

\renewcommand{\limsup}{\varlimsup}
\renewcommand{\liminf}{\varliminf}

\providecommand{\eps}{\varepsilon}
\DeclareMathOperator{\VMO}{\mathrm{VMO}}
\DeclareMathOperator{\arcsinh}{\mathrm{arcsinh}}

\newcommand{\tildef}{G}

\newtheorem{theorem}{Theorem}[section]

\newtheorem{prop}[theorem]{Proposition}

\newtheorem{lemma}[theorem]{Lemma}
\newtheorem{proposition}[theorem]{Proposition}
\newtheorem{definition}[theorem]{Definition}

\makeatletter 

\@addtoreset{equation}{section}
\makeatother

\newcommand{\invlabel}[1]{\raisebox{\ht\strutbox}{\hypertarget{#1}{}}}
\newcommand{\invref}[2]{{\rm(\hyperlink{#1}{#2})}}

\theoremstyle{definition}
\newtheorem{ex}[theorem]{Example}
\makeatletter
\newtheoremstyle{italicRemark}
{}
{}
{\normalfont}
{}
{\itshape}
{.}
{ }
{}
\makeatother

\theoremstyle{italicRemark}
\newtheorem{Rmk}[theorem]{Remark}
\newenvironment{Proof}[3][\unskip]{\begin{proof}[Proof of #2 \ref{#3}{#1}]
  }{\end{proof}}
\providecommand{\eps}{\varepsilon}

\newcommand{\weaksto}{\overset{*}{\rightharpoonup}} 
\newcommand{\ren}{\mathrm{ren}} 
\newcommand{\sg}{\mathrm{sg}} 
\newcommand{\sys}{\mathrm{sys}(\mathcal N)}

\newcommand{\rms}{\mathrm{s}} 
\newcommand{\rma}{\mathrm{a}} 
\newcommand{\mres}{\mathbin{\vrule height 1.6ex depth 0pt width 0.13ex\vrule height 0.13ex depth 0pt width 1.3ex}}

\hypersetup{bookmarksdepth=2}

\title[Universality of renormalisable mappings in 2d]{Universality of renormalisable mappings in two dimensions: the case of polar convex integrands}
\date{\today}

\begin{document}

\author{Christopher Irving}
\address{Faculty of Mathematics\\ Technical University Dortmund\\ Vogelpothsweg 87\\ 44227\\ Dortmund\\ Germany}
\email{christopher.irving@tu-dortmund.de}

\author{Benoît Van Vaerenbergh}
\address{Université catholique de Louvain\\ Institut de Recherche en Mathématique et Physique\\ Chemin du cyclotron 2 bte L7.01.02\\ 1348 Louvain-la-Neuve\\ Belgium}
\email{benoit.vanvaerenbergh@uclouvain.be}


\begin{abstract}
  We establish universality of the renormalised energy for mappings from a planar domain to a compact manifold, by approximating subquadratic polar convex functionals of the form $\int_\Omega f(|\mathrm{D} u|)\,\mathrm{d} x$. 
  The analysis relies on the condition that the vortex map ${x}/{\lvert x\rvert}$ has finite energy and that $t\mapsto f (\sqrt{t})$ is concave. We derive the leading order asymptotics and provide a detailed description of the convergence of $\mathrm{W}^{1,1}$--almost minimisers, leading to a characterization of second-order asymptotics.
	At the core of the method, we prove a ball merging construction (following Jerrard and Sandier's approach) for a general class of convex integrands. 
  We therefore generalize the approximation by $p$--harmonic mappings when $p\nearrow 2$ and can also cover linearly growing functionals, including those of area-type. 
\end{abstract}

\maketitle

\begingroup\hypersetup{hidelinks}
\tableofcontents
\endgroup

\section{Introduction}

Given a bounded Lipschitz planar domain $\Omega \subset  \R^2$, a  compact Riemannian manifold $\mathcal  N \subset \R^\nu$ ($\nu \in \mathbb{N}\setminus\{0\}$), and boundary datum $g \in \WW^{\sfrac{1}{2},2}(\partial \Omega;\mathcal  N)$,  is a mapping $u \in \WW^{1,2}(\Omega;\mathcal  N) \doteq \WW^{1,2}(\Omega;\R^\nu)\cap \{\text{for a.e. } x \in \Omega, u(x)\in \mathcal  N\}$ solving the minimisation problem 
\begin{equation}
	\label{eq:dirichletproblem}
	\inf \left\{\int_\Omega \frac{|\D u|^2}{2} \d x : \begin{matrix}
		u \in  \WW^{1,2}(\Omega;\mathcal  N) \\ \Tr_{\partial\Omega}u = g
	\end{matrix}\right\}.
\end{equation}
If $\WW^{1,2}_g(\Omega;\mathcal  N) \doteq \WW^{1,2}(\Omega;\mathcal  N) \cap \{\Tr_{\partial\Omega}u = g\}$ contains at least one map, the existence of harmonic mappings in this Dirichlet class can be shown by the Direct Method of the Calculus of Variations. However, due to topological obstructions, there may exist $g \in \WW^{\sfrac{1}{2},2}(\partial \Omega;\mathcal  N)$ such that $\WW^{1,2}_g(\Omega;\mathcal  N) = \emptyset$. 
In particular, such a $g$ exists when $\mathcal N$ is not simply connected, or equivalently the fundamental group $\pi_1(\mathcal N)$ is non-trivial; see for instance \cite{bethuel1995extensions, van2024extension}.

While $\WW^{1,2}_g(\Omega;\mathcal N)$ may be empty, if we relax our admissible class we can find singular mappings $u \colon \Omega \to \mathcal N$ attaining the boundary condition.
Thus we may seek to understand \eqref{eq:dirichletproblem} in the obstructed case through a suitable sequence of approximate problems which we can solve, by analysing the behaviour of these approximate solutions in the singular limit.
Moreover, we seek to obtain a \emph{singular renormalisable harmonic map} in the limit; a notion introduced in \cite{MonteilEtAl2021,MonteilEtAl2022}.

More precisely, we consider a sequence $(f_n)_n$ of Young functions $f_n : \R_+ \to \R_+$ (referred to as integrands in the sequel) for which the \emph{vortex energy} 
\begin{equation}\label{eq:defOfV_n}
  \mathcal V(f_n) \doteq 2\int_0^1 f_n\Big( \frac1r \Big) r \d r = \frac{1}{\pi}\int_{\B}f_n(\lvert\D u_{\mathrm{V}}(x)\rvert) \d x  \quad\mbox{is finite,}
\end{equation}
where $u_{\mathrm{V}}(x) = \frac{x}{\lvert x\rvert}$ is the vortex map, defined on the unit ball $\B \subset \mathbb R^2$.
Condition \eqref{eq:defOfV_n} ensures the set
\begin{equation}\label{eq:admissible_class}
	\mathrm A_g(f_n) \doteq \Big\{ u \in \WW^{1,1}(\Omega;\mathcal N) :\int_{\Omega} f_n(\lvert \D u\rvert)\d x \text{ is finite}, \ \  \Tr_{\partial \Omega}u =g \Big\}
\end{equation} 
is non-empty (see Lemma \ref{lemma:existenceofmin}) and therefore
\begin{equation}\label{eq:theproblem}
  \inf_{v \in \mathrm{A}_g(f_n)} \int_{\Omega} f_n(\lvert\D v\rvert)\d x \mbox{ is finite,}
\end{equation}
so we obtain a sequence of well-posed minimisation problems.
To recover the Dirichlet case \eqref{eq:dirichletproblem}, we consider general families $(f_n)_n$ satisfying Definition \ref{defn:family_fn}, which will include
\begin{alignat*}{4}
  \invlabel{item:pcase}\mathrm{(a)} &\quad f_n(z) = \frac{\lvert z\rvert^{p_n}}{p_n}  \quad &&\mbox{where $p_n \nearrow 2$} \ \ &&\mbox{($p$--energy)}, \\
  \invlabel{item:area}\mathrm{(b)} &\quad  f_n(z) = \frac{\sqrt{1 + \delta^2_n \lvert z\rvert^2 } - 1}{\delta^2_n} \quad &&\mbox{where $\delta_n \searrow 0$} \ \ &&\mbox{(rescaled area functional)}, \\
  \invlabel{item:truncated}\mathrm{(c)} &\quad f_n(z) = \begin{cases} \frac12 \lvert z\rvert^2 &\mbox{if $\lvert z\rvert  < \kappa_n$,} \\ \kappa_n (\lvert z\rvert - \kappa_n/2) &\mbox{if $\lvert z\rvert \geq \kappa_n$}\end{cases} \quad &&\mbox{where $\kappa_n \rightarrow \infty$} \ \ &&\mbox{(truncated Dirichlet energy)},\\
  \invlabel{item:sublog}\mathrm{(d)} &\quad f_n(z) = \frac{\lvert z\rvert^2}{2(1+\eta_n\lvert \log \lvert z\rvert\rvert)^{2}} \quad &&\mbox{where $\eta_n \searrow 0$,} \ \ &&\mbox{($\mathrm{L} \log^{-2}\mathrm{L}$--energy)}.
\end{alignat*}

Observe that the aforementioned families convergence pointwise to the $\frac12 \lvert z\rvert^2$ as $n \to \infty$.
For such integrands we will show the following prototype theorem:

\begin{theorem}\label{thm:firstorder}
  Let $(f_n)_n$ be as in Definition \ref{defn:family_fn}, which includes any of \invref{item:pcase}{a}--\invref{item:sublog}{d}.
  Then for any $g \in \WW^{\sfrac12,2}(\partial\Omega;\mathcal N)$ we have
  \begin{equation}
    \lim_{n \to \infty} \frac{1}{\mathcal V(f_n)} \inf\left\{ \int_{\Omega} f_n(\lvert\D v\rvert)\d x : \begin{matrix} v \in \WW^{1,1}(\Omega;\mathcal N)\\ \Tr_{\partial\Omega}v=g\end{matrix}\right\}= \Esgnu([g]).
  \end{equation} 
  Here $\Esgnu([g])$ is the singular energy; \textit{cf.}\,Definition \ref{def:singen}.
\end{theorem}
Since $\lim_{n \to \infty} \mathcal V(f_n) = \infty$ by Fatou's lemma, Theorem \ref{thm:firstorder} quantifies the rate of blow-up of \eqref{eq:theproblem} up to first order by the singular energy, which is a homotopy invariant quantity that does not depend on the choice of the sequence $(f_n)_n$.
When $\Omega = \B\subset \R^2$ is the unit ball and $\mathcal N = \mathbb S^1$, we have $\Esgnu([g]) = \pi |\deg(g)|$, where $\deg(g)$ is the topological degree of $g$ viewed as a map $\partial \mathbb B = \mathbb S^1 \to \mathbb S^1$.

Moreover we prove a stronger result; namely we are able to completely characterise the next leading order of the asymptotic expansion in $n$, and establish convergence results for asymptotically minimising sequences.

\begin{theorem}\label{thm:2ndorder}
  Under the assumptions of Theorem \ref{thm:firstorder}, suppose $(u_n)_n \subset \WW^{1,1}(\Omega;\mathcal N)$ is an asymptotically minimising sequence in the sense that $\Tr_{\partial\Omega}u_n = g$ for all $n$ and
\begin{equation}\label{eq:uscren_intro}
	 \int_{\Omega} f_n(\lvert \D u_n\rvert) \d x - \inf\left\{ \int_{\Omega} f_n(\lvert\D v\rvert)\d x : \begin{matrix} v \in \WW^{1,1}(\Omega;\mathcal N)\\ \Tr_{\partial\Omega}v=g\end{matrix}\right\} \xrightarrow{n \to \infty} 0.
\end{equation} 
Then, passing to a subsequence, $u_n$ converges strongly in $\WW^{1,1}$ to a limit map $u_{\ast} \in \WW^{1,1}(\Omega;\mathcal N)$ with $\Tr_{\partial \Omega}u_{\ast} = g$ satisfying 
\begin{equation}u_{\ast} \in \WW^{1,2}\Big(\Omega\setminus \bigcup_{i = 1}^\kappa \B(a_i,\rho);\mathcal N\Big) \quad\mbox{for all $\rho>0$},\end{equation} 
where $\kappa \in \mathbb N$ and $a_1,\dotsc,a_{\kappa} \in \Omega$. 
Moreover, for all $\rho>0$ sufficiently small, 
\begin{enumerate}[\rm(i)]
  \item\label{item:minimisingaway} $u_*$ is a minimising harmonic map on $\Omega\setminus \bigcup_{i = 1}^\kappa \B(a_i,\rho)$ with respect to its own boundary conditions, 
	\item 	\label{item:strongintthmintro}\( \ds
	\int_{\Omega \setminus \bigcup_{i=1}^{\kappa} \B(a_i,\rho)} f_n(\lvert \D u_{\ast} - \D u_n\rvert) \d x \xrightarrow{n\to \infty} 0,
	\)
\item $u_{\ast}$ is renormalisable in that the renormalised energy
	\begin{equation}\label{eq:defereninintro}
    \mathrm{E}_{\ren}^{1,2}(u_{\ast})= \lim_{\rho \to 0}\left[\int_{\Omega \setminus \bigcup_{i = 1}^k \B(a_i,\rho)}\frac{\lvert \D u_*\rvert^2}{2} \d x - \Esgnu([g])\log \frac{1}{\rho}\right] \quad\mbox{exists and is finite,}
	\end{equation}
  and we have
  \begin{equation} \label{eq:secondorder_intro}
  \lim_{n \to \infty} \left[\int_{\Omega} f_n(\lvert\D u_n\rvert) \d x - \mathcal{V}(f_n) \Esgnu(g)\right] = \mathrm{E}_{\ren}^{1,2}(u_{\ast}) + \mathrm H([u_{\ast},a_i])_{i=1}^{\kappa},\end{equation}
	where $\mathrm H([u_{\ast},a_i])_{i=1}^{\kappa}$ is defined in \eqref{eq:entropy_remainder}.
\end{enumerate}
\end{theorem}
The goal of this paper is to identify general conditions for sequences of integrands $(f_n)_n$, covering the four examples \invref{item:pcase}{a}--\invref{item:sublog}{d}, such that Theorems \ref{thm:firstorder} and \ref{thm:2ndorder} hold, 
working in the spirit of $\Gamma$-convergence (also known as variational convergence). These results show that \emph{the renormalised energy is universal}: regardless of the sequence of approximating integrands $(f_n)_n$, the limiting mappings have finite renormalised energy \eqref{eq:defereninintro}. 
This universality suggests that the manifold constrained harmonic extension of a map $g : \partial \Omega \to \mathcal N$ should be defined through the minimisation of the renormalised energy $\Erennu(u)$; this quantity was introduced in \cite{MonteilEtAl2022}. Under certain assumptions on the domain, the manifold and the boundary data, we can explicitly compute the manifold constrained harmonic extension:

\begin{theorem}\label{thm:univrotex}
  Suppose $\Omega = \mathbb{B}^2$ is the unit ball centred at the origin, $\mathcal{N} = \mathbb{S}^1$  and $g : \partial \mathbb{B} = \mathbb{S}^1 \to \mathbb{S}^1$ is defined by $g(x) \doteq x$ for all $x \in \partial \mathbb B$,
	and let $(f_n)_n$ satisfy Definition \ref{defn:family_fn} (which includes  \invref{item:pcase}{a}--\invref{item:sublog}{d}).
  Then for any asymptotically minimising sequence $(u_n)_n \subset \WW^{1,1}(\mathbb{B}, \mathbb{S}^1)$ satisfying $\operatorname{tr}_{\partial \mathbb{B}} u_n = g$ for all $n$ and \eqref{eq:uscren_intro}, we have
  \begin{equation}
    u_n \to u_{\mathrm V}(x) = \frac{x}{\lvert x\rvert} \quad\mbox{strongly in $\WW^{1,1}(\B,\mathbb S^1)$}.
  \end{equation}
\end{theorem}
Theorem \ref{thm:univrotex} asserts \emph{universality of the vortex map $u_{\mathrm{V}}$} and can be extended to other manifolds and boundary data, \textit{cf.}\,Theorem \ref{thm:universalityofthevortexmap}. We underline that in Theorem \ref{thm:univrotex} the full sequence $u_n$ converges instead of some subsequence.

\medbreak

The assumptions we use are detailed in Section \ref{sec:growth_structures}, but we will consider
convex polar integrands of the form $f_n(\lvert z\rvert)$, with $f_n : t \in \R_+ \to \R_+$; here, and throughout the paper, we will use $\lvert \cdot \rvert$ to denote the Euclidean / Frobenius norm.
 We impose that the sequence $f_n(t)$ converges pointwise to $t^2/2$, that each $f_n$ has finite vortex energy $\mathcal V(f_n)$ (see \eqref{eq:defOfV_n}), and that $t \in \R_+ \mapsto f_n \circ \sqrt{t}$ is concave. 
These assumptions covers the case of $p$-harmonic mappings, which was studied by \textsc{Hardt} and \textsc{Lin} in \cite{HardtLin1987} for $\mathcal N = \mathbb S^1$, and for general $\mathcal N$ in \cite{vanschaftingen2023asymptotic}.
Our assumptions allows for a much more general class of integrands however, including linear growth functionals such as \invref{item:area}{b}, and integrands which are not strictly convex such as \invref{item:truncated}{c}, which is $|z|^2/2$ truncated at $z = \kappa$ (see Definition \ref{defn:kappa_truncation}). 

 The family \invref{item:area}{b} is motivated by the observation that the Dirichlet energy can be obtained from the expansion of the area functional, namely
\begin{equation}\label{eq:area_expansion}
  \sqrt{1 + \lvert z\rvert^2} = 1 + \frac{\lvert z\rvert^2}2 + O(\lvert z\rvert^4),\quad  z \to 0.
\end{equation}
Note that this integrand arises by considering the area of the graph of $\lvert u\rvert$; that is $\mathcal H^2(\{ (x,\lvert u(x)\rvert) : x \in \Omega\})$.
To capture the second term in the expansion \eqref{eq:area_expansion}, we consider the change of variables $u_{\delta}(x) = u(\delta x)$ for $u \in \WW^{1,1}(\B;\mathcal N)$ which satisfies
\begin{equation}\label{eq:rescaledareafunctional}
  \int_{\B(0,1/\delta)} \sqrt{1 + \lvert \D u_{\delta}\rvert^2} -1 \d x = \int_{\B} \frac{\sqrt{1+\delta^2\lvert \D u\rvert^2} - 1}{\delta^2} \d x,
\end{equation} 
which leads us to consider $f_n$  defined as in \invref{item:area}{b}, which converges pointwise to $\frac12 \lvert z\rvert^2$ as $\delta_n \searrow 0$. 
This approximation has recently appeared in \cite{palmer2024lifting} as a relaxation to the problem of quadrilateral and hexagonal meshing, formulated as a minimisation problem over suitable fiber bundles. 
This approximation was also used in 
\cite[\S A.3(2) p.185]{polya1951isoperimetric} to prove the P\'olya-Szegő inequality for $p=2$.

One could sacrifice the convexity of the integrands and the intrinsic character of the problem by using a Ginzburg-Landau approach, where for every $\varepsilon > 0$ one considers
\begin{equation}
	\inf \left \{ \int_\Omega \frac{|\D u|^2}{2} + \frac{F(u)}{\varepsilon^2} \, \d x : u \in \WW_g^{1,2}(\Omega; \R^\nu)\right\}.
\end{equation}
Here $F(u)$ behaves as $\dist(u; \mathcal N)^2$ near $\mathcal N$ (see \cite[(1.4), \S2]{MonteilEtAl2021} for detailed assumptions on $F$), and $\mathcal N$ is a connected compact submanifold of $\R^\nu$. 
Results similar to Theorems \ref{thm:firstorder} and \ref{thm:2ndorder} can be obtained in that setting (see \cite[Thm.~7.1, 7.3]{MonteilEtAl2021}), and this relaxation has also been applied to the theory of meshing; see for instance \cite{beaufort2017computing}.
However, this approach is extrinsic since $F(u)$ depends on an isometric embedding $\mathcal N \subset \R^\nu$, which always exists \cite{nash1956imbedding}, but is highly non-unique. Moreover, $F(u)$ is typically non-convex, as the classical case attests \cite{bethuel2017ginzburg} for $\mathcal N = \mathbb S^1$, where $F(u) = (|u|^2 - 1)^2 / 4$. The extrinsic behaviour is recorded in the obtained renormalised limit; see the $Q$-term in \cite[Thm.~7.3(vi)]{MonteilEtAl2021}.

\textit{Innovations}. 
While many prior works focused on specific approximations such as the $p$-energy and the Ginzburg Landau relaxation, we initiate a systematic study of functionals $\mathcal F_n$ approximating the Dirichlet energy along which this asymptotic analysis can be carried out.
Our assumptions on the integrands $(f_n)_n$ is sufficiently general to allow for linearly growing integrands such as \invref{item:area}{b} and \invref{item:truncated}{c}, thus leading us to work in the framework of $\LL^1$ convergence.

The key results of our analysis are contained in the compactness theorem (Theorem \ref{thm:convsub}), valid for sequences whose energy is bounded in a renormalised sense.
The result crucially relies on a lower bound to control the singularities (see Proposition \ref{prop:genmergingballlemma}) based on the classical merging ball procedure of \textsc{Jerrard} \cite{Jerrard1999} and \textsc{Sandier} \cite{Sandier1998}.
While the a-priori bounds allow us to extract a subsequence of $(\D u_n)_n$ which converges \weaklystar (in the sense of measures), by separate asymptotic estimates of the sequence near and away from the singularities, we can improve this convergence to hold weakly in $\LL^1$ (see Proposition \ref{prop:compacitnoprblem} and Theorem \ref{thm:convsub}\eqref{item:L1bound}).
Moreover, in the case of asymptotically minimising sequences, in Theorem \ref{thm:strongconv} we can improve the convergence of gradients to hold strongly $\LL^1$, by exploiting the uniform convexity of the limiting $\LL^2$ scale (see Proposition \ref{prop:strong_conv}).

\textit{Future directions}.
Our claim of the universality of the renormalized energy naturally raises the question of whether more general classes of integrands can be considered (\textit{e.g.}\,to allow for non-polar integrands with an $(x,u)$-dependence), 
 or even discrete versions (see \cite{canevari2028defects} for $\mathcal N = \mathbb S^1$).
Further approximations based on variants of \textsc{Bourgain-Brézis-Mironescu}-type formulas \cite{bourgain2001anther} can also be considered, as discussed in \cite{solci2024nonlocal} for $\mathcal N = \mathbb S^1$. 
Additionally, our analysis is confined to mappings in $\WW^{1,1}$; however for linearly growing integrands we are naturally lead to consider minimisers in $\BV$ (see Section \ref{sec:existence_minima} for a further discussion).

\textit{Structure of the paper.} Assumptions on the sequence of approximating integrands $(f_n)_n$ are described in Section \ref{sec:growth_structures}. The singular energy and the renormalised energies are detailed in Section \ref{sec:energies-and-topological-quatities}. In Section \ref{sec:asymptotic_bounds}, we present key ingredients needed for the formal $\Gamma$-convergence analysis; namely an upper bound corresponding to the $\Gamma$-limsup, and a merging ball construction which is used to obtain asymptotic lower bounds. The main compactness theorem is proven in Section \ref{sec:compactness-results}, which relies on compactness results for unconstrained mappings and for traces, and corresponds to the $\Gamma$-liminf analysis. In Section \ref{sec:convergence-of-almost-minimisers}, we prove an extended version of Theorem \ref{thm:2ndorder} (see Theorem \ref{thm:strongconv}). Finally, in Section \ref{sec:universality-of-vortex-map}, we demonstrate that in the specific case where $\Omega = \B$ is the unit ball, $\mathcal N = \mathbb S^1$ and $g(z) = z$, the entire sequence $(u_n)_n$ in Theorem \ref{thm:2ndorder} converges to the vortex map $x/|x|$ asserting the \textit{universality of the vortex map}.

\section{Growth structures for the Calculus of Variations}\label{sec:growth_structures}

In this section we describe the set of assumptions we impose on our integrands $f_n$ (namely in Definition \ref{defn:family_fn}), and record some examples and useful properties thereof. 
We will write $\mathbb{R}_+ = [0,\infty)$ in what follows.

\subsection{Assumptions}

Let $f \colon \mathbb{R}_+ \to \mathbb R_+$ be a \emph{Young function}, that is a convex non-decreasing function such that $f(0)=0$.

We will impose that $f$ is subquadratic; either in the sense of the \emph{natural sense} that
\begin{equation}\label{hyp:ndec2} \tag*{$(\mathrm{Dec})^{\mathrm{n}}_2$}
  t \mapsto \frac1{t^2} f(t) \quad\mbox{is non-increasing on $[t_0,\infty)$},
\end{equation} 
or in the \emph{controlled sense} that
\begin{equation}\label{hyp:cdec2} \tag*{$(\mathrm{Dec})^{\mathrm{c}}_2$}
  t \mapsto \frac1{t} f'_+(t) \quad\mbox{is non-increasing on $[t_0,\infty)$},
\end{equation} 
where $t_0 \geq 0$; here $f_+'$ denotes the right-derivative as in Lemma~\ref{lemma:rightder}. 
The \ref{hyp:ndec2} condition is well--known in the Orlicz literature, see for instance \cite[\S 2.1]{HarjulehtoHasto2019}.
One may readily verify  by elementary calculus that \ref{hyp:cdec2} is equivalent to
\begin{equation}
  t \mapsto f(\sqrt{t}) \quad\mbox{is concave on $[t_0,\infty)$}.
\end{equation} 
Furthermore \ref{hyp:cdec2} implies \ref{hyp:ndec2} by \eqref{eq:convex_linearisation}, but is in general is a stronger condition, see Example \ref{ex:not_cdec2} below.

We will also impose the integrability condition
\begin{equation}\label{hyp:int}
  	\mathcal V(f) = 2\int_0^1 f\left( \frac1r \right) r \d r = 2\int_1^{\infty} \frac{f(s)}{s^2}\, \frac{\dd s}s < \infty, \tag*{$(\mathrm{vInt})$}
\end{equation} 
where the equality follows from the change of variable $s = \frac1r$. 
The factor $2$ is present to match the singular energy in the asymptotic expansion; see Theorem \ref{thm:firstorder}.
We call $\mathcal V(f)$ the \emph{vortex energy} of $f$; see Proposition \ref{prop:weakl2_integrable} below for consequences of this condition.

\begin{definition}\label{defn:family_fn}
  A \emph{sequence of approximating integrands} is a family
 $(f_n)_n$ of Young functions $f_n : \mathbb R_+ \to \R$ satisfying the following:
\begin{enumerate}[\rm (a)]
  \item There is $t_0 \geq 0$ such that each $f_n$ satisfies \ref{hyp:cdec2} with $t_0$ uniform in $n$,
  \item each $f_n$ satisfies \ref{hyp:int},
  \item the sequence $(f_n)_n$ converges pointwise to $\frac12 t^2$, that is
    \begin{equation}\label{eq:conv_pointwise}
  \lim_{n \to \infty} f_n(t) = \frac12 t^2 \quad\text{for all } t \geq 0.
\end{equation} 

\end{enumerate}
\end{definition}

Note that combining \ref{hyp:ndec2} and \ref{hyp:int} gives
\begin{equation}\label{eq:subquadratic}
	\lim_{t \to +\infty}\frac{f(t)}{t^2} = 0,
\end{equation}
so we infer that $f$ satisfies a subquadratic growth condition.
However under \ref{hyp:int}, the condition \ref{hyp:ndec2} is a stronger condition than \eqref{eq:subquadratic}, as illustrated by Example \ref{ex:not_ndec2} below. 

The \ref{hyp:int} condition is closely connected to $\LL^{2,\infty}$, whose definition we briefly recall.
\begin{definition}
  Let $\Omega \subset \mathbb R^2$ be open, the \emph{weak-$\LL^2$ space} $\LL^{2,\infty}(\Omega)$ is the space of measurable functions $u \colon \Omega \to \R$ such that
  \begin{equation}
    \lVert f \rVert_{\LL^{2,\infty}(\Omega)} \doteq \Big[\sup_{t>0} t^2 \,\Leb^2\left( \{ x \in \Omega : \lvert f(x)\rvert > t\} \right)\Big]^\frac{1}{2} <\infty,
  \end{equation} 
  where $\Leb^2$ denotes the Lebesgue measure on $\mathbb R^2$.
\end{definition}

Other names for weak-$\LL^2$ space are \emph{Marcinkiewicz space} \cite{marcinkiewicz1939interpolation} or \emph{Lorentz space} \cite{lorentz1950somenew}. Here, $\Leb^2$ denotes the $2$-dimensional Lebesgue measure. We point out that, despite the notation, $\lVert \cdot \rVert_{\LL^{2,\infty}(\Omega)}$ is not a norm as the triangle inequality does not hold. However, \(
	\lVert m_1 + m_2\rVert_{\LL^{2,\infty}(\Omega)} \leq  2\big(\lVert m_1\rVert_{\LL^{2,\infty}(\Omega)} +  \lVert m_2 \rVert_{\LL^{2,\infty}(\Omega)}\big).
\)

\begin{proposition}\label{prop:weakl2_integrable} Let $\Omega \subset \R^2$ be a domain and $\B(0,1) \subset \R^2$ be the unit ball. Then for a Young function $f : \R_+ \to \R$, the following conditions are equivalent:
	\begin{enumerate}[\rm (i)]
		\item \label{item:simpleintegral} the integrand $f$ satisfies \ref{hyp:int}; that is, $\mathcal V(f)$ is finite 
		\item \label{item:vortex} For the vortex map $u_{\mathrm{V}}(x) = x/|x|$ we have $\ds \int_{\B(0,1)} f(|\D u_{\mathrm{V}}|)\d x$ is finite, 
		\item \label{item:weakL2} for every $m \in \LL^{2,\infty}(\Omega)$, we have $\ds \int_{\Omega} f(\lvert m(x)\rvert)\d x$ is finite.
	\end{enumerate}
	Moreover if any of the above holds, for each non-zero measurable $m : \Omega \to \R_+$ we have the estimate
	\begin{equation}\label{eq:inequweakL2}
		\int_\Omega f\circ m \d x \leq 2  \|m\|_{\LL^{2,\infty}(\Omega)}^2 \int_{{\|m\|_{\LL^{2,\infty}(\Omega)}}/{\sqrt{\Leb^2(\Omega)}}}^\infty \frac{f(s)}{s^3}\d s.
	\end{equation}
\end{proposition}

\begin{ex}\label{eg:vortex_map}
  The prototypical example of a mapping $u$ for which $|\D u| \in \LL^{2,\infty}(\B_1)$ but $|\D u| \notin \LL^{2}(\B_1)$
  is the \emph{vortex map} $u_{\mathrm{V}} \colon \B_1 \to \mathbb S^1$, which is defined as
  \begin{equation}
    u_{\mathrm{V}}(x) = \frac{x}{\lvert x\rvert}.
  \end{equation} 
  Indeed since $\lvert \D u_{\mathrm{V}}\rvert(x) = \frac1{\lvert x\rvert}$, we have
  \(
    \lVert \D u_{\mathrm{V}} \rVert_{\LL^{2,\infty}(\B_1)} = \sqrt{\pi} < \infty, \) but \( \lVert \D u_{\mathrm{V}} \rVert_{\LL^2(\B_1)}^2 =  2\pi\int_0^1 {\d r}/{r}   = \infty,
  \) 
  where we used polar coordinates to compute the $\LL^2$-norm.
  For a general manifold $\mathcal N$, we can produce similar examples by taking an immersed curve $\gamma \colon \mathbb S^1 \to \mathcal N$ and considering $\gamma(x/\lvert x\rvert)$.
\end{ex}

\begin{Proof}{Proposition}{prop:weakl2_integrable}
	The equivalence between \eqref{item:vortex} and \eqref{item:simpleintegral} follows from integrating in polar coordinates and since $|\D u_{\mathrm{V}}|(x) = 1/|x|$; indeed
	\begin{equation}\label{eq:computingthevortexenergy}
    \int_{\B(0,1)} f(|\D u_{\mathrm{V}}|) = 2\pi \int_0^1 f\left(\frac{1}{r}\right) r\d r = \pi \mathcal V(f).
	\end{equation}
	Since $1/|x| \in \LL^{2,\infty}(\Omega)$ by Example \ref{eg:vortex_map}, we have \eqref{item:weakL2} implies \eqref{item:vortex}.
  To show that \eqref{item:simpleintegral} implies \eqref{item:weakL2}, it suffices to prove the estimate \eqref{eq:inequweakL2}.
  Let $m \in \LL^{2,\infty}(\Omega)$ be non-negative, and set $M \doteq{\|m\|_{\LL^{2,\infty}}}/\sqrt{\Leb^2(\Omega)}>0$. Since $s \mapsto f(s)/s^3$ is continuous on $(0,\infty)$ and \ref{hyp:int} holds, we have
  \(
    \int_M^{\infty} {f(s)}/{s^3} \d s < \infty.
  \) 
  Then by the layer cake representation (see \cite[Cor.~2.2.34]{willem2013functional}) and by a change of variables we have, 
  \begin{equation*}
  \int_\Omega f\circ m\d x = \int_{0}^\infty \Leb^2(\{x \in \Omega : f\circ m(x) > t\})\d t 
                          = \int_{0}^\infty f'(t)\, \Leb^2(\{x \in \Omega : m(x) > t\})\d t.
\end{equation*}
Note that since $f$ is convex, we have $f$ is locally Lipschitz, so this change of variables is justified (see for instance \cite[Thm.~3.9, 6.7]{EvansGariepy2015}).
We can split this integral as
	\begin{align*}
		\int_\Omega f\circ m &\leq  \Leb^2(\Omega)\int_{0}^M f'(t) \d t + \|m\|_{\LL^{2,\infty}(\Omega)}^2\int_{M}^\infty \frac{f'(t)}{t^2}\d t \\
		&= \Leb^2(\Omega)f(M) + \|m\|_{\LL^{2,\infty}(\Omega)}^2 \left[2\int_M^\infty \frac{f(s)}{s^3}\d s -\frac{f(M)}{M^2} \right].
	\end{align*}
 By our choice of $M$,  we obtain \eqref{eq:inequweakL2} and hence \eqref{item:weakL2}. This concludes the proof.
\end{Proof}

\subsection{Some examples}\label{sec:some-examples-and-their-the-vortex-energy}

\begin{ex}
We detail here explicitly the vortex energy for the examples from the introduction. For the $p$-case \invref{item:pcase}{a}, we compute, for $f_p(t) = |t|^p/p$,
\[
  \mathcal V(f_p) = \frac{2}{(2 - p)p} =  \frac1{2-p} + O(1) \quad\text{as } p\nearrow 2
\]
so that  Theorem \ref{thm:firstorder} reads 
  \begin{equation}
	\lim_{p \nearrow 2} (2 - p) \inf_{v \in \WW^{1,p}_g(\Omega;\mathcal N)} \int_{\Omega} \frac{\lvert\D v\rvert^p}{p}\d x = \Esgnu([g]).
\end{equation} 
For the rescaled area functional \invref{item:area}{b}, taking $f_\delta(t) = (\sqrt{1+\delta^2 t^2} - 1)/\delta^2$ we can compute
\begin{equation*}
	\mathcal V(f_\delta) =  2\int_0^1 \frac{\sqrt{1+ \delta^2r^{-2}} - 1}{\delta^2} \,r \d r 
  = f_{\delta}(1) + \arcsinh(\frac{1}{\delta})
  = \log\frac1{\delta} + O(1)\quad \text{as } \delta\searrow 0,
\end{equation*}
where $\arcsinh(t) = \log(t + \sqrt{1+t^2})$. In that case Theorem \ref{thm:firstorder} reads as 
  \begin{equation}
	\lim_{\delta \searrow 0} \frac{1}{\log{\frac{1}{\delta}}} \inf_{\WW^{1,1}_g(\Omega;\mathcal N)} \int_{\Omega} \frac{\sqrt{1 + \delta^2 |\D v|^2} - 1}{\delta^2}\d x  = \Esgnu([g]).
\end{equation} 
For the truncated functional \invref{item:truncated}{c}, namely
\(
	f_\kappa(z) = \frac{|z|^2}{2}\chi_{\{|z| < \kappa\}} + \kappa(|z| - \frac{\kappa}{2}) \chi_{\{|z| > \kappa\}}
\),
for $\kappa >1$,
\[
	\mathcal V(f_\kappa) = \int_{\frac{1}{\kappa}}^1 \frac{\d t}{t} + 2 \kappa\int_{0}^{\frac{1}{\kappa}}\left(\frac{1}{t} - \frac{\kappa}{2}\right)t \d t= \log \kappa + \frac{7}{4}.
\]
For \invref{item:sublog}{d} we can more generally consider $f_\eta(t) = \frac{t^2}{2(1 + \eta|\log t|)^\alpha}$ with $\eta>0$, for $\alpha >1$ fixed. 
By computing $f'_{\eta}(t)$ we see that each $f_{\eta}$ is a Young function satisfying \ref{hyp:cdec2} with $t_0=0$, and by a change of variables $s = - \eta \log r$ we can also compute
\begin{equation*}
	\mathcal V(f_\eta) = \int_0^1 \frac{\d r}{r(1 + \eta|\log r|)^\alpha}
  = \frac{1}{\eta}\int_{0}^{\infty} \frac{\d s}{(1 + s)^\alpha} = \frac1{(\alpha-1)\eta}. 
\end{equation*}
\end{ex}

We will also describe two counterexamples relating to the \ref{hyp:int}, \ref{hyp:ndec2} and \ref{hyp:cdec2} conditions.

\begin{ex}\label{ex:not_ndec2}
  The following is a variant of {\cite[\S 11, Eg.~6]{Maligranda1989}}.
  For $p \in (1,2)$ and $\beta>0$ to be determined, consider
  \begin{equation}
    f(t) = \frac1p t^p \left( 1 + \frac1{\beta} \sin(p \log t) \right) \quad t>0.
  \end{equation} 
  Then for the choice $\beta \geq 5$, we have $f$ is a Young function such that $\frac1p t^p \leq f(t) \leq \frac1p (1+\frac1{\beta})t^p$.
  Now we can compute 
  \begin{equation}
    f'(t) = t^{p-1} \left( 1 + \frac{\sqrt{2}}p \sin(p \log t + \pi/4) \right), 
  \end{equation} 
  and so
  \begin{equation}
    \frac{t_kf'(t_k)}{f(t_k)} = p \left( 1+ \frac{1}{\sqrt{\beta^2-1}} \right) \quad \text{ where } t_k = \exp\left( \frac1p ((2k+1)\pi-\arcsin(1/\beta) ) \right)
  \end{equation} 
  for all $k\in \mathbb N$ with $k \geq 1$.
  Since \ref{hyp:ndec2} is equivalent to $\frac{tf'(t)}{f(t)} \leq 2$ holding for sufficiency large $t>0$, we see this condition is violated for $\beta = 5$ and $p > \frac{4}{23}(12 - \sqrt{6}) \simeq 1.66\ldots$, thereby giving examples of Young functions satisfying \ref{hyp:int} and \eqref{eq:subquadratic}, but not \ref{hyp:ndec2}.

  To obtain a sequence of such functions $f_n$ converging to $\frac12 t^2$, we can choose $p_n = \frac{2n+1}{n+1}$ and $\beta_n = \sqrt{n^2-1}$ for $n \geq 6$.
\end{ex}

\begin{ex}\label{ex:not_cdec2}
  To construct a family of integrands satisfying \ref{hyp:ndec2} but not \ref{hyp:cdec2}, we will first define a function $g \colon \mathbb R_+ \to \mathbb R_+$ via
  \begin{equation}
    g(t) = g_p(t) \doteq t^{p-2} + \sum_{k=3}^{\infty} a_k\,h\left( \frac{t-k}{r_k} \right),
  \end{equation} 
  where $p \in (1,2)$ is sufficiently close to $2$, the ``spike function'' $h \colon \mathbb R \to \mathbb R$ is defined via
  \begin{equation}
    h(t) = \begin{cases} t &\text{if } 0 \leq t < \frac12, \\ 1-t &\text{if } \frac12 \leq t \leq 1, \\ 0 &\text{otherwise,} \end{cases}
  \end{equation} 
  and we set $a_k = \frac{2(2-p)}{k^2}$, $r_k = k^{1-p}$ for all $k \in \mathbb N$, $k \geq 3$.
  We then claim the family
  \begin{equation}
    f(t) = f_p(t) \doteq \int_0^t g_p(s) s \d s
  \end{equation} 
  is a sequence of Young functions satisfying \ref{hyp:ndec2}, \ref{hyp:int}, but does not satisfy \ref{hyp:cdec2} on $[t_0,\infty)$ for any $t_0\geq 0$.
  Note that setting $a_k = 0$ for all $k$ would recover the $p$-energy, so as $p \nearrow 2$ we have $f_p$ converges pointwise to $\frac12 t^2$.

  To establish this claim, note that since $g$ is non-negative we have $f$ is increasing, and we can compute
  \begin{equation}
    f_+''(t) = g(t) + tg_+'(t) \geq (p-1) t^{p-2} +\sum_{k=3}^{\infty} \frac{a_kt}{r_k} \left(\mathbbm{1}_{[k,k+r_k/2)}(t) - \mathbbm{1}_{[k+r_k/2,k+r_k)}(t)\right).
  \end{equation} 
  Here $f_{+}''$ denotes the right derivative of $f'$; see Lemma \ref{lemma:rightder}.
  To show this is non-negative,
  it suffices to check the endpoint cases $t = k+r_k$, which rearranges to give $\frac{p-1}{2-p} \geq 2(1+k^{-p})^{3-p}$. 
  By checking the $k=3$ case, this holds for all $p \in (p_0,2)$ with $p_0 \simeq 1.706\ldots$
  Hence, for this range of $p$, we have $f$ is a Young function.
  To see that $f$ does not satisfy \ref{hyp:cdec2}, we note that $\frac{f'(t)}t = g(t)$ by construction, and since $g'_+(k) = \frac{2-p}{k^{3-p}} > 0$ for each $k \geq 3$, it follows that $g$ is not non-increasing.

  To verify \ref{hyp:ndec2}, observe this is equivalent to showing that $g(t) = \frac{f'(t)}t \leq \frac{2f(t)}{t^2}$. Moreover, by integrating by parts we can compute
  \begin{equation}
    f(t) = \int_0^t g(s) s \d s = \frac12 t^2 g(t) - \frac12 \int_0^t g_+'(s) s^2 \d s,
  \end{equation} 
  so rearranging we see that \ref{hyp:ndec2} is equivalent to showing that $\int_0^t g_+'(s) s^2 \d s \leq 0$ for all $t \geq 0$.
  To compute this integral, we first calculate the contribution of each of the ``spike'' terms.
  For this we have
  \begin{equation}
    \begin{split}
    \int_{k}^{k+r_k} \frac{\dd}{\dd s} \left(  a_k\,h\left( \frac{s-k}{r_k} \right)  \right)  s^2 \d s 
    &=  - a_k \left(  \frac12 k r_k + \frac14 r_k^2 \right) \leq 0,
  \end{split}
  \end{equation} 
  so the contribution over a full spike is negative.
  Hence for any $t \in (k+r_k,k+1)$ we see that
  $  \int_0^t g'_+(s) s^2 \d s \leq -\frac{(2-p)}p k^p  < 0$,
  and if $t \in (k,k+r_k/2)$, then
  \begin{equation}
    \begin{split}
    \int_0^t g'_+(s) s^2 \d s 
    &\leq  - \frac{2-p}p\, k^{p} + \frac{a_k}{r_k}\int_{k}^{k+r_k/2} s^2  \d s 
    \leq -(2-p) \left( \frac1p k^p - 2\right),
    \end{split}
  \end{equation} 
  noting that $r_k \leq 1$, where the last line is non-negative since $k \geq 3$. The case $t \in (k+r_k/2,k+r_k)$ follows since the relevant terms are decreasing, thereby verifying \ref{hyp:ndec2}.

  Finally to show \ref{hyp:int} holds, we can compute
    $\int_{k}^{k+r_k} a_k\, h\left( \frac{t-k}{r_k} \right) \d t \leq \frac{2-p}2$ for each $k$,
  and hence we can estimate $f(t) \leq t^p + \frac12 (t+1)$ for all $t \geq 0$. 
  Thus it follows that $\mathcal V(f)$ is finite,
  verifying all of the claimed properties of $f$ for all $p_0 < p < 2$.
\end{ex}

\subsection{Convexity properties}
\begin{lemma}\label{lemma:rightder}
  Let $f$ be a convex function on $\mathbb{R}_+$.
  Then the right derivative
  \begin{equation}
    f_+'(t) \doteq \lim_{h \searrow 0} \frac{f(t+h)-f(t)}h
  \end{equation} 
  exists for all $t \geq 0$ and defines a non-decreasing and right-continuous function on $\mathbb{R}_+$.
\end{lemma}

We refer to \cite[Thm.~24.1]{Rockafellar1997} for a proof. The proof of this assertion relies on the \emph{slope inequality} for convex functions
\begin{equation}\label{eq:slope_inequality}
  \frac{f(s)-f(r)}{s-r} \leq \frac{f(t) - f(r)}{t-r} \leq \frac{f(t)-f(s)}{t-s}
\end{equation} 
valid for all $0 \leq r < s < t$.
From this we also infer that
\begin{align}\label{eq:convex_linearisation}
  f_+'(s)(t-s) \leq f(t) - f(s) \leq f_+'(t)(t-s) &\quad\mbox{for all $s \leq t$}.
\end{align}

  Analogously we can also define left derivatives $f_-'$, which will be non-decreasing and left-continuous on $(0,\infty)$. However, we have opted for right-derivatives to ensure $f'_+(t)$ is well-defined at $t=0$.

\begin{definition}\label{defn:kappa_truncation}
Let $f$ be a Young function, then for $\kappa>0$ define the \emph{$\kappa$-truncation} $T_{\kappa}f$ of $f$ as 
\begin{equation}\label{eq:kappa_truncation}
  T_{\kappa}f(t) = 
  \begin{cases}
    f(t) &\text{if } 0 \leq t \leq \kappa, \\
    f_+'(\kappa)(t - \kappa) + f(\kappa) &\text{if } t > \kappa.
  \end{cases}
\end{equation} 
which exists for all $\kappa \geq 0$.
\end{definition}

Observe that $T_{\kappa}f$ is itself a Young function. If $f$ satisfies \ref{hyp:ndec2} or \ref{hyp:cdec2}, then $T_{\kappa}f$ also does with the same $t_0$.
We also have $T_{\kappa}f(t) \leq f(t)$ for all $t \geq 0$ and $T_{\kappa}f$ converges pointwise monotonically to $f$ as $\kappa \to \infty$. We also note that $T_\kappa f$ always satisfies \ref{hyp:int} and $\mathcal V(T_\kappa f) \leq \mathcal V(f)$.

\begin{lemma}\label{lem:fn_conv}
  Let $(f_n)_n$ be a sequence of Young functions converging pointwise to $f(t) = t^2/2$ on $\R_+$.
  Then, as $n \to \infty$,
  \begin{enumerate}[\rm(a)]
    \item\label{item:local_uniform} $f_n$ converges locally uniformly to $f$ on $\R_+${\rm ;} \emph{i.e.}\,uniformly on $[0,T]$ for each $T>0$,
    \item\label{item:derivative_pointwise} $f_{n,+}'(t)$ converges pointwise to $f'(t) = t$ on $\R_+$.
  \end{enumerate}
\end{lemma}

\begin{Proof}{Lemma}{lem:fn_conv}
Assertion \eqref{item:local_uniform} is proven in \cite[Thm.~10.8]{Rockafellar1997}. For \eqref{item:derivative_pointwise} we can use \eqref{eq:convex_linearisation} to show that
\begin{equation}
  \frac{f_n(s)-f_n(r)}{s-r} \leq f_{n,+}'(s) \leq \frac{f_n(t)-f_n(s)}{t-s}
\end{equation} 
for all $0 \leq r < s < t$.
Sending $n \to \infty$ followed by $r,t \to s$ establishes the result.
\end{Proof}

\begin{Rmk}
  An alternative way to truncate $f$ is to define 
\begin{equation}\label{eq:truncation_envelope}
  \tilde T_Lf(t) =  \inf\{ f(s) + L\lvert s -t\rvert \colon s \in \R_+\}
\end{equation} 
for some $L>0$.
One can show that the choice $L = f_+'(\kappa)$ gives $\tilde T_Lf = T_{\kappa}f$; indeed using \eqref{eq:convex_linearisation} one shows that $\tilde T_Lf(t) = f(t)$ for all $t \in [0,\kappa]$, and if $t > \kappa$ one has the infimum in \eqref{eq:truncation_envelope} is attained with the choice $s = \kappa$.
\end{Rmk}

\begin{lemma}\label{lemma:quadgrowth}
	Let $f$ be a Young function satisfying \ref{hyp:ndec2}.
	Then
  \begin{equation}
    0\leq  f(t) \leq f(\max\{1,t_0\}) ( 1+t^2) \quad\mbox{for all $t \geq 0$}.
  \end{equation} 
  In particular, for a family $(f_n)_n$ as in Definition \ref{defn:family_fn}, we can choose $C>0$ independently of $n$ such that
  \begin{equation}\label{eq:quadupper_uniform}
    0\leq f_n(t)\leq C(1+t^2) \quad\mbox{for all $t \geq 0$}.
  \end{equation} 
\end{lemma}

\begin{Proof}{Lemma}{lemma:quadgrowth}
  We can assume that the \ref{hyp:ndec2}--condition holds with $t_0 \geq 1$.
  Then for $t \geq t_0$ we have
  \begin{equation}
    0 \leq f(t) \leq \frac{f(t_0)}{t_0^2} t^2 \leq f(t_0) t^2.
  \end{equation} 
  When $0 \leq t \leq t_0$, since $f$ is non-decreasing we have $0 \leq f(t) \leq f(t_0)$, from which the result follows.
  For a family $(f_n)$, since $\lim_{n \to \infty} f_n(t_0) = \frac12 t_0^2$ we can choose $C>0$ such that \eqref{eq:quadupper_uniform} holds uniformly in $n$.
\end{Proof}

The following is used in the proof of the upper bound (Proposition \ref{prop:upperBound}), whose proof crucially relies the stronger \ref{hyp:cdec2} condition.

\begin{lemma}\label{lemma:quad_comparison}
  Suppose $f$ is a Young function satisfying \ref{hyp:cdec2} for some $t_0 \geq 0$.
  Then for all $\sigma,\rho >0$ we have
  \begin{equation}\label{eq:dct_upperbound}
    f(\sigma) - f(\rho) \leq \frac{f_+'(\max\{\rho,t_0\})}{2\max\{\rho,t_0\}} (\sigma^2 - \rho^2) + f(t_0) + f_+'(t_0)t_0.
  \end{equation}
\end{lemma}

\begin{Proof}{Lemma}{lemma:quad_comparison}
  Define $g(t) = f(\sqrt{t})$, which is concave on $[t_0,\infty)$ by the \ref{hyp:cdec2} condition.
  By a variant of \eqref{eq:convex_linearisation} for concave functions, for $\sigma,\rho> t_0$ we have
  \begin{equation}
    g(\sigma^2) - g(\rho^2) \leq g_+'(\rho^2) (\sigma^2 - \rho^2). 
  \end{equation} 
  That is,
  \begin{equation}\label{eq:cdec_comparison}
    f(\sigma) - f(\rho) \leq \frac{f_+'(\rho)}{2\rho} (\sigma^2 - \rho^2)
  \end{equation} 
  holds.
  For the remaining cases, note that if $\rho \leq t_0 \leq \sigma$, then by \eqref{eq:cdec_comparison} with $t_0$ in place of $\sigma$ we have
  \begin{equation}
    f(\sigma) - f(\rho) \leq \frac{f_+'(t_0)}{2t_0} (\sigma^2 - t_0^2) + f(t_0)-f(\rho) \leq \frac{f_+'(t_0)}{2t_0} (\sigma^2 - \rho^2) + f(t_0),
  \end{equation} 
  whereas if $\sigma \leq t_0 \leq \rho$ applying \eqref{eq:cdec_comparison} similarly gives
  \begin{equation}
    f(\sigma) - f(\rho) \leq f(\sigma) - f(t_0) + \frac{f_+'(\rho)}{2\rho} (t_0^2 - \rho^2) \leq \frac{f_+'(\rho)}{2\rho}(\sigma^2-\rho^2) + \frac12 f_+'(t_0) t_0.
  \end{equation} 
  Finally if both $\rho,\sigma \leq t_0$ we simply bound
  \begin{align}
    f(\sigma) - f(\rho) \leq f(t_0) \leq \frac{f_+'(t_0)}{2t_0} (\sigma^2-\rho^2) + f(t_0) &\quad\mbox{if $\sigma \geq \rho$},\\
    f(\sigma) - f(\rho) \leq 0 \leq \frac{f_+'(t_0)}{2t_0} (\sigma^2-\rho^2) + \frac12 f_+'(t_0)t_0 &\quad\mbox{if $\sigma < \rho$},
  \end{align} 
  thereby verifying \eqref{eq:dct_upperbound} in all cases.
\end{Proof}

\begin{ex}\label{eq:weak2strong_dec2}
  Suppose $(f_n)_n$ is a sequence of approximating integrands as in Definition \ref{defn:family_fn}, except that each $f_n$ only satisfies \ref{hyp:ndec2} in place of \ref{hyp:cdec2}. Then we can define a modified family $(\bar f_n)_n$ by
  \begin{equation}
    \bar f_n(t) = \int_0^t \frac{2f_n(s)}{s} \d s \quad\mbox{for all $t \geq 0$},
  \end{equation} 
  which satisfies Definition \ref{defn:family_fn}, including \ref{hyp:cdec2}, and 
  that $f_n \leq \bar f_n$ pointwise using \ref{hyp:ndec2}.
  Moreover, integrating by parts gives
  \begin{equation}\label{eq:vint_diff}
    \mathcal V(\bar f_n) = \int_1^{\infty} \bar f_n(s) \,\frac{\dd s}{s^3} = \mathcal V(f_n) + \frac12 \bar f_n(1).
  \end{equation} 
  This construction can be used to prove Theorem \ref{thm:firstorder} assuming only \ref{hyp:ndec2}, which we will sketch in Remark \ref{rmk:firstorder_proof}.
\end{ex}

\subsection{Existence of minima}\label{sec:existence_minima}

We will briefly discuss the existence of minima for our approximate problems, using the Direct Method.
Here $\Omega \subset \mathbb R^2$ will be a bounded Lipschitz domain, and $\mathcal N$ will be a compact Riemannian manifold embedded in $\mathbb R^{\nu}$; we will make this precise in Section \ref{sec:energies-and-topological-quatities}. Compactness of $\mathcal N$ implies that there exists a constant $C>0$ depending only on $\mathcal N$ such that any measurable map $u : \Omega \to \mathcal N$ satisfies $\|u\|_{\LL^\infty(\Omega)} \leq C$; we will therefore never insist on the $\LL^p$--convergence of mappings and will express our results in terms of $\LL^1$--convergence.

\subsubsection{Existence of competitors}

\begin{lemma}\label{lemma:existenceofmin}
  Let $g \in \WW^{\sfrac12,2}(\partial\Omega; \mathcal N)$. If the vortex energy
  \(\mathcal V(f)\) is finite,
  the admissible class $\mathrm A_g(f)$ in \eqref{eq:admissible_class}
  is non-empty.
\end{lemma}

Two proofs of Lemma \ref{lemma:existenceofmin} are available. One is  established in \cite[Prop.~2.19]{vanschaftingen2023asymptotic} by showing that there exists a renormalisable map (see Proposition \ref{prop:renormalisedenergy}) realizing the constraints. 
Since renormalisable mappings have gradient in $\LL^{2,\infty}$ by \cite[Thm.~5.1]{MonteilEtAl2022}, combining this with Lemma \ref{prop:weakl2_integrable} we infer the result.
 Alternatively, we can use the fact that any map $g \in \WW^{\sfrac12,2}(\partial\Omega;\mathcal N)$ can be extended into a map whose gradient is  $\LL^{2,\infty}(\Omega)$, with a non-linear estimate relating the $\LL^{2,\infty}$-norm and the $\WW^{\sfrac12,2}(\partial\Omega;\mathcal N)$-norm; see \cite{bulanyi2023singular,petrache2017controlled}.

\subsubsection{The Direct Method}\label{sec:the-direct-method}
Since $f$ is assumed to be convex, we have $t \mapsto \frac{f(t)}t$ is non-decreasing on $[0,\infty)$. Hence we either have  $f$ is \emph{asymptotically linear} in that 
\begin{equation}\label{eq:asymptotic_linear}
	f^{\infty}(1) \doteq \lim_{t \to \infty} \frac{f(t)}t \quad\text{exists and is finite},
\end{equation} 
or $f$ is \emph{asymptotically superlinear} in that
\begin{equation}\label{eq:asymptotic_superlinear}
	\lim_{t \to \infty} \frac{f(t)}t = \infty.
\end{equation}
In the asymptotically linear case, the Direct Method naturally leads us to consider generalised minimisers in the class $\BV$ of functions of bounded variation;
this requires us to give a meaning to
\begin{equation}
	\mathcal F(u,\Omega) \doteq \int_\Omega f(\lvert \D u\rvert) \,\d x
\end{equation}  
for mappings $u \in \BV(\Omega;\mathcal N) \doteq \BV(\Omega;\mathbb R^{\nu}) \cap \{ u(x) \in \mathcal N\text{ a.e.\,in } \Omega\}$. 
This will be done through the \emph{$\LL^1(\Omega;\mathbb R^{\nu})$-relaxation}
\begin{equation}\label{eq:defn_relaxation}
	\overline{\mathcal F}(u,\Omega) = \inf\left\{ \liminf_{n \to \infty} \mathcal F(u_n,\Omega) : u \in \WW_g^{1,1}(\Omega;\mathbb R^{\nu}),\ u_n \to u \text{ strongly in } \LL^1(\Omega) \right\}.
\end{equation} 
Here $\WW^{1,1}_g(\Omega;\mathcal N)$ denotes the space of functions $u \in \WW^{1,1}(\Omega;\mathcal N)$ such that $\Tr_{\partial\Omega}u=g$, which corresponds to the admissible class \eqref{eq:admissible_class} with the choice $f(z) =|z|$.

\begin{proposition}\label{prop:existence}
	Let $f$ be a Young function and $g \in \WW^{\sfrac12,2}(\Omega;\mathcal N)$.
	\begin{enumerate}[\rm (a)]
		\item \label{item:bvcase} If \eqref{eq:asymptotic_linear} holds, then $\overline{\mathcal F}$ admits a minimiser $u \in \BV(\Omega;\mathcal N)$.
		\item \label{item:sobolevcase} If \eqref{eq:asymptotic_superlinear} holds, then $\mathcal F$ admits a minimiser $u \in \WW^{1,1}_g(\Omega;\mathcal N)$.
	\end{enumerate}
\end{proposition}

  
To treat both cases at once, we rely on the following representation formula (see \cite{AmbrosioDalMaso1992,GoffmanSerrin1964}):
\begin{lemma}\label{lemma:bv_relaxation}
	Let $f$ be a Young function which is asymptotically linear in that \eqref{eq:asymptotic_linear} holds.
	Then the $L^1(\Omega;\mathbb R^{\nu})$-relaxation $\overline{\mathcal F}(\cdot,\Omega)$ admits the measure representation
	\begin{equation}\label{eq:linear_relaxation}
		\overline{\mathcal F}(u,\Omega) = \int_{\Omega} f(\lvert \D^{\mathrm{a}}u\rvert) \d x + f^{\infty}(1) \lvert \D^{\mathrm{s}}u\rvert(\Omega) + f^{\infty}(1)\int_{\partial\Omega} \lvert \Tr_{\partial\Omega}u-g\rvert \d\mathcal{H}^1
	\end{equation} 
	for $u \in \BV(\Omega;\mathcal N)$, where $\D u = \D^{\mathrm{a}}u\,\Leb^2 + \D^{\mathrm{s}}u$ is the Lebesgue decomposition of the gradient.
	Moreover, this relaxation is \weaklystar  lower semicontinuous in $\BV(\Omega;\mathcal N)$.
\end{lemma}
\begin{Proof}{Lemma}{lemma:bv_relaxation}
  Denote the right-hand side of \eqref{eq:linear_relaxation} by $\widetilde{\mathcal F}(\cdot,\Omega)$, which by \cite[Thm.~3]{GoffmanSerrin1964} is weakly${}^{\ast}$ sequentially lower-semicontinuous in $\BV(\Omega;\mathbb R^{\nu})$, noting the recession function is given by $f^{\infty}(\lvert z\rvert) = f^{\infty}(1)\lvert z\rvert$.
  Since $\widetilde{\mathcal F}(v,\Omega) = \mathcal F(v,\Omega)$ for all $v \in \WW^{1,1}_g(\Omega;\mathbb R^{\nu})$, from the lower-semicontinuity it follows that $\widetilde{\mathcal F}(u,\Omega)\leq\overline{\mathcal F}(u,\Omega)$ for all $u \in \BV(\Omega;\mathbb R^{\nu})$.
  Moreover, if $u \in \BV(\Omega;\mathcal N)$, letting $\widetilde\Omega \Supset \Omega$ be a Lipschitz domain and $G \in \WW^{1,2}(\widetilde\Omega \setminus\overline\Omega;\mathbb R^{\nu})$ be an extension of $g$, by \cite[Lem.~B.2]{Bildhauer2003} there exists a sequence $U_k \in \WW^{1,2}(\widetilde\Omega;\mathbb R^{\nu})$ such that $U_k = G$ on $\widetilde\Omega \setminus \overline\Omega$ and $U_k$ converges $\BV$-strictly to
  \begin{equation}
    U(x) = \begin{cases} u(x) &\text{if } x \in \Omega, \\ G(x) &\text{if } x \in\widetilde\Omega\setminus\overline\Omega;\end{cases}
  \end{equation} 
  that is $U_k \to U$ strongly in $\LL^1$ and $\lvert\D U_k \rvert(\widetilde\Omega) \to \lvert \D U\rvert(\widetilde\Omega)$.
  Then by Reshetnyak's theorem (see \emph{e.g.}\,\cite[Thm.~10.3]{rindler2018calculus}) we infer the restrictions $u_k = U_k \rvert_{\Omega} \in \WW^{1,2}_g(\Omega;\mathbb R^{\nu})$ satisfies
  \begin{equation}
    \lim_{k \to \infty} \int_{\Omega} f(\lvert \D u_n\rvert) \d x = \widetilde{\mathcal F}(u,\Omega),
  \end{equation} 
  so the converse inequality follows, thereby establishing that $\widetilde{\mathcal F} = \overline{\mathcal F}$.
\end{Proof}

\begin{Proof}{Proposition}{prop:existence}
  If \eqref{eq:asymptotic_linear} holds, then by Lemma \ref{lemma:bv_relaxation} we have $\overline{\mathcal F}$ is \weaklystar se\-quen\-tial\-ly lower semicontinuous on $\BV(\Omega;\mathcal N)$, so by the Direct Method we can infer the existence of a minimiser.
  Indeed, since \eqref{eq:admissible_class} is non-empty by Lemma \ref{lemma:existenceofmin}, we can take a minimising sequence $(u_n)_n$, and assume $f$ is not identically zero (else every admissible map is a minimiser).
  Then, there is some $t_0 \geq 0$ such that $f'(t_0) >0$, and so $f(t) \geq f'(t_0) t$ for all $t \geq t_0$, so it follows that $(|\D u_n|)_n$ is bounded in $\LL^1(\Omega)$.
  Therefore, we can pass to a \weaklystar convergent subsequence $u_{n_k} \weaksto u$ in $\BV(\Omega;\mathcal N)$, and by lower-semicontinuity of $\overline{\mathcal F}$ it follows that $u$ is a minimiser.

  In the asymptotically superlinear case where \eqref{eq:asymptotic_superlinear} holds, for $\kappa>0$ we consider the truncations $T_{\kappa}f$ from Definition \ref{defn:kappa_truncation}, noting that $T_{\kappa}f^{\infty}(1) = f'_+(\kappa)$.
  Then if we take a minimising sequence $(u_n)_n$ for $\mathcal F$ in $\mathrm A_g(f)$, then analogously to the asymptotically linear case above, we obtain a weak${}^{\ast}$-limit $u \in \BV(\Omega;\mathcal N)$.
  By lower semicontinuity applied to the relaxation $\overline{\mathcal F_{\kappa}}$ associated to $T_{\kappa}f$, we have
  \begin{equation}
    \begin{split}
      \overline{\mathcal F_{\kappa}}(u,\Omega) 
      &= \int_{\Omega} T_{\kappa}f(\lvert \D^{\rma}u\rvert) \d x + f_+'(\kappa ) \lvert \D^{\rms}u\rvert(\Omega) + f_+'(\kappa) \int_{\partial\Omega} \lvert \tr_{\partial\Omega}u-g\rvert\d\Hau^1\\
      &\leq \liminf_{n \to \infty} \int_{\Omega} T_{\kappa}f(\lvert \D u_n\rvert) \d x \leq \liminf_{n \to \infty} \mathcal F(u_n,\Omega) < \infty.
    \end{split}
  \end{equation} 
  By \eqref{eq:asymptotic_superlinear} and \eqref{eq:convex_linearisation} we have
  $f_+'(\kappa) \to \infty$ as $\kappa \to \infty$,
   so it follows that $\lvert \D^{\rms}u\rvert(\Omega) = 0$ and that $\tr_{\partial\Omega}u =g$ on $\partial\Omega$.
  Therefore $u \in \WW^{1,1}(\Omega;\mathcal N)$, and by the monotone convergence theorem we have
  \begin{equation}
    \mathcal F(u,\Omega) = \lim_{\kappa\to\infty} \int_{\Omega} T_{\kappa}f(\lvert \D u\rvert)\d x \leq \liminf_{n \to \infty} \mathcal F(u_n,\Omega).
  \end{equation} 
  Since $(u_n)_n$ was a minimising sequence, it follows that $u$ is a minimiser of $\mathcal F$.
\end{Proof}

\begin{Rmk}
  Under \eqref{eq:asymptotic_linear}, the minima \( u \) in \( \BV(\Omega; \mathcal{N}) \) obtained in Proposition \ref{prop:existence}\eqref{item:bvcase} do not necessarily satisfy \( \Tr_{\partial \Omega} u = g \), where we consider the interior trace in the sense of \emph{e.g.}, \cite[Thm.~5.6]{EvansGariepy2015}.
  This may arise due to potential jumps across the boundary, which is penalised by the term
  $  f^{\infty}(1) \int_{\partial\Omega} \lvert \Tr_{\partial\Omega}u -g\rvert\d\mathcal H^1$
  in the relaxation $\overline{\mathcal F}(u)$. For this reason, in this setting our analysis will be confined to almost minimising sequences in  $\WW^{1,1}(\Omega;\mathcal N)$ (see Theorem \ref{thm:2ndorder}). In the superlinear case \eqref{eq:asymptotic_superlinear} the trace attainment is not an issue, as the proof of Proposition \ref{prop:existence}\eqref{item:sobolevcase} shows.
\end{Rmk}

\section{Energies and topological quatities}\label{sec:energies-and-topological-quatities}

In this section, we introduce scalar quantities that measure the non trivial-ness of the free homotopy class of a map $\partial \Omega \to \mathcal N$ as well as the renormalised energy, which the limiting energy maps will minimise.

Throughout the whole text, the domain $\Omega$ is assumed to be an open connected bounded set of $\R^2$ with Lipschitz boundary. 
Concerning the manifold $\mathcal N$, our setting covers the case where $\mathcal N$ is a smooth compact connected Riemannian manifold which is assumed to be isometrically embedded in the Euclidean space $\R^\nu$ (with canonical metric) for some $\nu\in \mathbb N$. Isometric embedding is not a restriction by Nash embedding theorem \cite{nash1956imbedding}.

We will also denote by $\B(a,\rho)$ the closed ball in $\mathbb R^2$ centred at $a \in \mathbb R^2$ with radius $\rho>0$, and we will sometimes denote the boundary by $\mathbb S^1(a,\rho) \doteq \partial \B(a,\rho)$.
The unit circle will be denoted by $\mathbb S^1 = \mathbb S^1(0,1)$.

\subsection{Topological quantities}
\subsubsection{Minimal length of the free homotopy class}

A curve of vanishing mean oscillation $\gamma \in \mathrm{VMO}(\mathbb S^1; \mathcal N)$ can be given a free homotopy class that extends the notion in the continuous category (see \cite{BrezisNirenberg1995}). In the sequel all the homotopy are assumed to be free. Every $\mathrm{VMO}(\mathbb S^1; \mathcal N)$-homotopy class contains a smooth length minimising geodesic in its homotopy class; indeed any $\mathrm{VMO}(\mathbb S^1; \mathcal N)$ curve is $\mathrm{VMO}$-homotopic to a smooth curve by \cite[Cor.~4]{BrezisNirenberg1995}, so the claim follows from the classical result for the smooth category (see \emph{e.g.}\,\cite[Prop.~6.28]{lee2018riem}). In particular, the following infimum is reached by such a geodesic.
\begin{definition}
Given $\gamma \in \mathrm{VMO}(\mathbb S^1; \mathcal N)$, denote by $[\gamma]$ its $\mathrm{VMO}(\mathbb S^1; \mathcal N)$-homotopy class. The \emph{minimal length of the free homotopy class $[\gamma]$} is 
	\begin{equation*}
		\lambda([\gamma]) = \inf \left\{\int_{\mathbb S^1}|\tilde \gamma'|\d \mathcal H^1 : \tilde \gamma \in \dot \WW^{1,1}(\mathbb{S}^1, \mathcal  N) \text{ and $\tilde\gamma$ is freely homotopic to } \gamma \right\}.
	\end{equation*}
By a \emph{minimising geodesic}, we mean a non-constant curve $\gamma \in \VMO(\mathbb S^1;\mathcal N)$ which attains this infimum.
By means of a reparametrisation, we will often assume that $\lvert \gamma'\rvert$ is constant.
\end{definition}
The \emph{systole} of the manifold $\mathcal{N}$ is the length of the shortest closed non-trivial geodesic on $\mathcal{N}$, namely,
\begin{equation}\label{eq:systole}
	\text{sys}(\mathcal{N})\doteq\inf\{\lambda([\gamma]): \gamma \in\VMO(\mathbb{S}^1, \mathcal  N)  \text{ is not freely homotopic to a constant}\}.
\end{equation}
The \emph{length spectrum} of the manifold $\mathcal{N}$ is the set  $\{\lambda([\gamma]) : \gamma \in \VMO(\mathbb{S}^1, \mathcal  N)\}\subset \{0\} \cup [\text{sys}(\mathcal{N}), +\infty)$.
\begin{lemma}\label{lemma:discreteness}
The length spectrum is a discrete set of the real line.
\end{lemma}
In particular, we have $\sys > 0$.
We refer to \cite[Prop.~3.2]{MonteilEtAl2021} for a proof; note that compactness of $\mathcal N$ is crucial for Lemma \ref{lemma:discreteness} to hold.

If $u \in \WW^{1,2}(\B(a,\sigma)\setminus \B(a,\rho);\mathcal N)$ for some $a\in \R^2$, $\sigma >0$ and all $\rho \in (0,\sigma)$, the trace $\Tr_{\mathbb S^1(a,\rho)}u \in \WW^{\sfrac12,2}(\mathbb S^1(a,\rho); \mathcal N) \subset \VMO(\mathbb S^1(a,\rho); \mathcal N)$ has a well defined homotopy class which is independent of $\rho$. Hence, minimal length of the free homotopy class of $u$ at $a$ is defined by
\begin{equation}\label{eq:lambdaua}
	\lambda([u,a]) \doteq \lambda(\Tr_{\mathbb S^1(a,\rho)}u).
\end{equation}

\subsubsection{Singular energy}
We will use the following variant of the notion of the topological resolution \cite[Def.~2.2]{MonteilEtAl2022}.
\begin{definition}\label{def:topres} A finite family of closed curves $\gamma_i \in \VMO(\mathbb S^1; \mathcal N)$, $i = 1, \dots,k$, with \(k \in \N\) is a \emph{topological resolution} of the map $g$
	whenever there exist $k$ non-intersecting closed balls $\B(a_i;\rho_i) \subset \Omega$ with \(\rho_i > 0\) and a map 
	\(u \in \WW^{1,2}(\Omega \setminus  \bigcup_{i = 1}^k\B(a_i;\rho_i); \mathcal N)\)
	such that $\tr_{\partial \Omega}u$ and $g$ are homotopic in $\VMO(\partial \Omega; \mathcal N)$ and $\Trace_{\mathbb{S}^1}u(a_i + \rho_i \,\cdot )$ and $\gamma_i$ are homotopic in $\VMO(\mathbb S^1; \mathcal N)$ for $i = 1,\dots,k$.
\end{definition}

\begin{definition}\label{def:singen}
	Let $g \in \VMO(\partial \Omega;\mathcal N)$, the \emph{singular energy} of $g$ is 
	\begin{equation}\label{eq:esg}
		\Esgnu([g]) = \inf\left\{\sum_{i = 1}^k \frac{\lambda([\gamma_i])^2}{4\pi } : \begin{matrix} k \in \mathbb N, \gamma_1, \dots, \gamma_k \in \VMO(\mathbb S^1; \mathcal N)\\ \text{ is a topological resolution of $g$}\end{matrix} \right\}.
	\end{equation}
\end{definition}
Observe that (see \cite[Prop.~2.6]{MonteilEtAl2022}) the singular energy is \emph{subadditive} in the sense that, if $(\gamma_i)_{i=1}^k$ is a topological resolution of $g$, we have 
\begin{equation}\label{eq:esg_subadditive}
  \Esgnu([g])  \leq \sum_{i=1}^k \Esgnu([\gamma_i]).
\end{equation} 

A \emph{minimal} topological resolution of $g$ is topological resolution that reach the infimum \eqref{eq:esg}. A curve $\gamma \in \VMO(\mathbb S^1; \mathcal N)$ is said \textit{atomic} if $\Esgnu([\gamma]) = {\lambda([\gamma])^2}/{4\pi}$. Often we will write $\Esgnu(g)$ instead of $\Esgnu([g])$.

A similar construction to \eqref{eq:lambdaua} can be done to define $\Esgnu([u,a])$ if $u \in \bigcup_{\rho>0}\WW^{1,2}(\B(a,\sigma)\setminus \B(a,\rho);\mathcal N)$ for some $\sigma>0$.

\subsection{Renormalisable mappings} 

In our analysis we will consider maps that are $\WW^{1,2}$-regular away from at most finitely many points, for which we introduce the following space.

\begin{definition}\label{defn:Rclass}
Given $\delta>0$, let $\RR^{1,2}_{\delta}(\Omega;\mathcal N)$ denote the space of mappings $u \in \WW^{1,1}(\Omega;\mathcal N)$ for which there is $k \in \mathbb N$ and points $a_1, \dotsc, a_k \in \Omega$ satisfying $\dist(\{a_i\}_{i=1}^k,\partial\Omega) \geq \delta$ such that
  \begin{equation}
    u \in \WW^{1,2}\left(\Omega \setminus \bigcup_{i=1}^k \B(a_i,\rho);\mathcal N\right) \text{ for all } \rho>0.
  \end{equation} 
  The points $a_1,\dotsc,a_{k}$ will be referred to as the \emph{singularities} of $u$.
  We also define $\RR^{1,2}(\Omega;\mathcal N) = \bigcup_{\delta>0}\RR^{1,2}_{\delta}(\Omega;\mathcal N)$.
\end{definition}

We will later show in Lemma \ref{lemma:density_of_the_R_class} that this class in dense in $\WW^{1,1}(\Omega;\mathcal N)$.
Also given finitely many points $a_1,\dots,a_k \in \Omega$, we will introduce the notation
\begin{equation}\label{eq:nonintersection}
  \rho_{\Omega}(a_1,\dots,a_k) \doteq \min\{ \dist(a_i,\partial\Omega), \lvert a_i - a_j\rvert : 1 \leq i<j \leq k\} >0.
\end{equation}

We also record the class of renormalisable mappings introduced in \cite[Def.~7.1]{MonteilEtAl2022}, which further specifies the behaviour of the mappings at the singular points.
The following is proven in \cite[Prop.~2.16]{vanschaftingen2023asymptotic} (see also \cite[Prop.~7.2]{MonteilEtAl2022}).

\begin{proposition}\label{prop:renormalisedenergy}
	Let $u \in  \WW^{1,1}(\Omega;\mathcal N)$ and let $a_1,\dots, a_k \in \Omega$ be $k \in \mathbb N$ distinct points. 
  Then the following are equivalent:
	\begin{enumerate}[\rm (i)]
		\item The {renormalised energy}
		\begin{equation}\label{eq:renormalisedenergy}
      \Erennu(u) = \lim_{\rho \to 0}\left[\int_{\Omega \setminus \bigcup_{i = 1}^k \B(a_i,\rho)}\frac{\lvert \D u\rvert^2}{2} \d x - \sum_{i = 1}^k \frac{\lambda([u,a_i])^2}{4\pi }\log \frac{1}{\rho}\right]
		\end{equation}
	exists as a limit and is finite.
		
		\item \label{item:actualuseofren} For each $0<\rho <\rho_\Omega(a_1,\dots,a_k)$, $|\D u| \in \LL^2(\Omega \setminus \bigcup_{i = 1}^k \B(a_i,\rho))$ and the map \begin{equation}\label{eq:renonthesing}
			 r  \mapsto \int_{\mathbb S^1(a_i,r)} \frac{\lvert \D u\rvert^2}{2} \d\Hau^1 - \frac{\lambda([u,a_i])^2}{4\pi r}
		\end{equation} is non-negative and integrable on $(0,\rho_\Omega(a_1,\dots,a_k))$.
	
	\end{enumerate}
	In any case 
  for all $0<\rho <\rho_\Omega(a_1,\dots,a_k)$,
	\begin{multline}\label{eq:erenwithoutlimit}
		\Erennu(u) + \sum_{i = 1}^k \frac{\lambda([u,a_i])^2}{4\pi}\log \frac{1}{\rho} \\= \int_{\Omega \setminus  \bigcup_{i = 1}^k\B(a_i,\rho)}\frac{\lvert \D u\rvert^2}{2} \d x + \sum_{i =1}^k \int_0^\rho\left [\int_{\mathbb S^1(a_i,r)} \frac{\lvert \D u\rvert^2}{2} \d\Hau^1 - \frac{\lambda([u,a_i])^2}{4\pi r} \right] \d r.
	\end{multline}
\end{proposition}

\begin{definition}\label{defn:renormalisable}
  A map $u \in \WW^{1,1}(\Omega;\mathcal N)$ is said to be \emph{renormalisable} if it satisfies any of the equivalent conditions in Proposition \ref{prop:renormalisedenergy}; that is $u \in \RR^{1,2}(\Omega;\mathcal N)$ and \eqref{eq:renormalisedenergy} is finite.
  We denote the space of all such mappings by $\WW^{1,2}_{\ren}(\Omega;\mathcal N)$.
\end{definition}

Prototypical examples of such mappings are \emph{vortex-type maps} taking the form $u_{\gamma}(x) = \gamma(x/\lvert x\rvert)$, where $\gamma \colon \mathbb S^1 \to \mathcal N$ is a minimising geodesic in its free homotopy class (\emph{i.e.}\,$\lambda([\gamma]) = \int_{\mathbb S^1} \lvert \dot\gamma\rvert \d\mathcal H^1$).
In this case, we have $u_{\gamma} \in \WW^{1,2}_{\ren}(\B(0,1);\mathcal N)$ with $\Eren^{1,2}(u_{\gamma}) = 0$.

Later in Section \ref{sec:compactness-results}, we will also use the notation $\mathrm{E}_{\ren}^{1,2}(u;\Omega)$ if we wish to specify the domain.
In particular, note that
\begin{equation}\label{eq:localised_renormalised}
  \mathrm{E}^{1,2}_{\ren}(u;\B(a_i,\rho)) + \frac{\lambda([u,a_i])^2}{4\pi} \log \frac1{\rho} = \int_{0}^{\rho}\left[\int_{\partial\B(a_i,r)} \frac12 \lvert \D u\rvert^2 \d\mathcal H^1 - \frac{\lambda([u,a_i])^2}{4\pi r}\right] \d r
\end{equation} 
for each $i = 1, \dotsc, k$, and $\rho < \rho_{\Omega}(a_1,\dotsc,a_{k})$.

In \cite[Prop.~7.2(iii)]{MonteilEtAl2022} it is shown that a property of renormalised maps is that for each $i \in \{1,\dots,k\}$, there exists a sequence $(\rho_\ell)_{\ell\in \mathbb N}$ converging to zero such that the sequence $(\tr_{\mathbb S^1}u(a_i + \rho_\ell\cdot))_{\ell \in \mathbb N}$ converges strongly to a geodesic $\gamma_i \in \mathrm C^1(\mathbb S^1; \mathcal N)$ in $\WW^{1,2}(\mathbb S^1; \mathcal N)$. We show in the next Example \ref{eg:geodesic_nonunique} that the selection of a sequence is in general necessary. 
However, the geodesics obtained along different subsequences lie in the same homotopy class (and moreover the same synharmony class; \emph{cf}.\,\cite[Def.~3.2, Prop.~7.2(iv)]{MonteilEtAl2022}) (see also \eqref{eq:synh_dist} in Section \ref{sec:universality-of-vortex-map}).

\begin{ex}[Non-uniqueness of the limiting geodesic]\label{eg:geodesic_nonunique}
  Let $\B = \B(0,1) \subset \mathbb R^2 \cong \mathbb C$ be the unit disc, and consider polar coordinates $(r,\theta) \mapsto z = r e^{i\theta}$. Then define $u \colon \B \to \mathbb S^1$ by
  \begin{equation}
    u(re^{i\theta}) = e^{i\theta + i \eta(r)}, \quad\mbox{where }\ \eta(r) = \log\log \left(1 + \frac1r\right).
  \end{equation} 
  Then since
  $  \lvert \D u\rvert^2 =
   \lvert \eta'(r)\rvert^2 + \frac1{r^2},$
  for $\rho \in (0,1)$ we have
  \begin{equation}
    \int_{\B \setminus \B(0,\rho)} \frac12 \lvert \D u\rvert^2 \,\d z - \pi \log\frac1\rho = \pi \int_{\rho}^1 r\lvert \eta'(r)\rvert^2 + \frac1{r} \,\d r - \pi \log \frac1{\rho},
  \end{equation} 
  which remains bounded as $\rho \to 0$ since $r \lvert \eta'(r)\rvert^2 \in \LL^1((0,1))$,
  thereby defining a map in $\WW^{1,2}_{\ren}(\B,\mathbb S^1)$.
  However if choose subsequences $\rho_k, \tilde\rho_k \searrow 0$ satisfying
  \begin{equation}
    \log\log \left(1+\frac1{\rho_k}\right) = 2k\pi, \quad \log \log \left(1+\frac1{\tilde\rho_k} \right) = 2k\pi - \frac{\pi}2
  \end{equation} 
  for all $k$, we see that
  \begin{equation}
    \Tr_{\mathbb S^1} u(\rho_n \cdot) \rightarrow \mathrm{id}_{\mathbb S^1}, \quad \Tr_{\mathbb S^1} u(\tilde\rho_n \cdot) \rightarrow e^{-i}\,\mathrm{id}_{\mathbb S^1},
  \end{equation} 
  exhibiting non-uniqueness of the limiting geodesic.

  We can also construct a similar example in the case of the flat torus $\mathcal N = \mathbb S^1 \times \mathbb S^1 \subset \mathbb C^2$, by defining \(
    v(re^{i\theta}) = (e^{i\theta}, e^{i\eta(r)})
  \) with $\eta$ as above.
  An analogous argument shows that $v \in \WW^{1,2}_{\ren}(\B; \mathcal N)$, and along the same subsequences $\rho_k, \tilde\rho_k$, we get two geodesics with disjoint images, namely $(\mathrm{id}_{\mathbb S^1},1)$ and $(\mathrm{id}_{\mathbb S^1},-1)$. 
\end{ex}

\section{Asymptotic bounds}\label{sec:asymptotic_bounds}

In this section, we derive bounds for our formal $\Gamma$-convergence analysis, starting with an upper bound in Proposition \ref{prop:upperBound}. Then, through a merging ball process, a lower bound for each fixed $n$ is established in Proposition \ref{prop:genmergingballlemma}, which allows us to establish asymptotic lower bounds detailed in Proposition \ref{prop:firstorderconve}.

\subsection{Upper bound}\label{sec:upper_bound}
We will show the renormalised limit along a fixed map is controlled by the renormalised energy, which corresponds to the $\Gamma$-limsup in the $\Gamma$-convergence analysis.
\begin{proposition}\label{prop:upperBound}
  Let $u \in  \WW^{1,2}_{\ren}(\Omega;\mathcal N)$ have $k \in \mathbb N$ distinct singularities $a_1,\dots, a_k \in \Omega$. Then if $(f_n)_n$ is a family of approximating integrands (as in Definition \ref{defn:family_fn}), we have
	\begin{equation}\label{eq:thelimit}
    \lim_{n \to \infty}\left[\int_{\Omega}f_n(|\D u|)\d x - \mathcal V(f_n)\sum_{i = 1}^k\frac{ \lambda([u,a_i])^2}{4\pi }\right] = \Erennu(u) + \mathrm{H}([u,a_i])_{i=1}^{\kappa},
	\end{equation}
  where we define
  \begin{equation}\label{eq:entropy_remainder}
    \mathrm{H}([u,a_i])_{i=1}^{\kappa} \doteq \sum_{i = 1}^k \frac{\lambda([u,a_i])^2}{4\pi}\log \frac{2\pi }{\lambda([u,a_i])}.
  \end{equation} 
\end{proposition}
We will record the following elementary calculation as we will need it several times later. Note that it also shows how the $\mathrm{H}([u,a_i])_{i=1}^{\kappa}$ term arises.

\begin{lemma}\label{lemma:entropy_error}
  Let $(f_n)_n$ be as in Definition \ref{defn:family_fn}. Then for all $\lambda, \rho>0$, we have
  \begin{equation}
    \lim_{n \to \infty} \left[\int_0^{\rho} 2 \pi r f_n \left( \frac{\lambda}{2\pi r} \right) \d r - \mathcal V(f_n) \frac{\lambda^2}{4 \pi}\right] = \frac{\lambda^2}{4 \pi}\log \frac{2 \pi \rho }{\lambda}.
  \end{equation} 
\end{lemma}
\begin{Proof}{Lemma}{lemma:entropy_error}
  For each $n$, change of variables $r \mapsto \frac{2\pi r}{\lambda}$ in the first integral gives
  \begin{equation}
      \int_0^{\rho} 2 \pi r f_n \left( \frac{\lambda}{2\pi r} \right) \d r - \mathcal V(f_n) \frac{\lambda^2}{4 \pi}
     =\frac{\lambda^2}{2\pi}\int_1^{\frac{2 \pi \rho}{\lambda}} f_n\left( \frac1r \right) r \d r,
  \end{equation} 
  and sending $n \to \infty$ using the dominated convergence theorem we have
  \begin{equation}
    \lim_{n \to \infty} \int_1^{\frac{2 \pi \rho}{\lambda}} f_n\left( \frac1r \right) r \d r = \frac12 \int_1^{\frac{2\pi \rho}{\lambda}} \frac{\dd r}{r} = \frac12 \log\frac{2 \pi \rho}{\lambda},
  \end{equation} 
  thereby establishing the result.
\end{Proof}

\begin{Rmk}
  In the case of the $p$-energy $f_p(t) = \frac{t^p}p$ (case {\rm(\hyperlink{item:pcase}{a})}), it was shown in
  \cite[Prop.~3.1]{vanschaftingen2023asymptotic} that
  \begin{equation*}
  \varlimsup_{p \nearrow 2}\left[\int_{\Omega} \frac{\lvert \D u\rvert^p}p  - \frac1{2-p} \sum_{i=1}^{\kappa} \frac{\lambda([u,a_i])^2}{4\pi}\right] \leq \mathrm{E}^{1,2}_{\ren}(u) + \sum_{i = 1}^k \frac{\lambda([u,a_i])^2}{4\pi}\Big[\frac12+\log \Big(\frac{2\pi }{\lambda([u,a_i])}\Big)\Big].\end{equation*} 
  This is consistent with Proposition \ref{prop:upperBound} above, and the discrepancy in the remainder term arises because $\mathcal V(f_p) = \frac2p\cdot \frac1{2-p} = \frac1{2-p} + \frac12 + o(1)$ as $p \nearrow 2$.
\end{Rmk}

\begin{Proof}{Proposition}{prop:upperBound}
  By Lemmas \ref{lemma:quadgrowth} and \ref{lemma:quad_comparison}, there exists $C_1,C_2>0$ such that
  \begin{align}
    \lvert f_n(t)\rvert &\leq C_1 (1 + t^2), \label{eq:quad_upper}\\
    f_n(t) - f_n(s) &\leq \frac{f'_{n,+}(\max\{s,t_0\})}{2\max\{s,t_0\}} (t^2 - s^2) + C_2,\label{eq:quad_upper2}
  \end{align}
  for all $n$ and $s,t \geq 0$, 
  where $t_0 \geq 1$ is as in the \ref{hyp:cdec2} condition, chosen uniformly in $n$.
	Fixing $\rho\in (0,\rho_\Omega(a_1,\dots,a_k)$, we will first show that
  \begin{equation} \label{eq:awayfromsing}
    \lim_{n \to \infty} \int_{\Omega\setminus \bigcup_{i = 1}^k \B(a_i,\rho)}f_n(|\D u|) = \int_{\Omega\setminus \bigcup_{i = 1}^k \B(a_i,\rho)}\frac{|\D u|^2}{2}.
  \end{equation} 
  This follows by the dominated convergence theorem; since $|\D u| \in \LL^2(\Omega \setminus \bigcup_{i=1}^k \B(a_i,\rho))$ by Proposition \ref{prop:renormalisedenergy}, using \eqref{eq:quad_upper} we obtain a uniform majorant.

  We now consider the contribution in $\B(a_i,r)$ for each $i = 1, \dotsc, k$.
  Using polar coordinates, we can write 
  \begin{equation}\label{eq:upperbound_mainestimate}
    \begin{split}
      &\int_{\B(a_i,r)} f_n(\lvert \D u\rvert) \d x - \mathcal V(f_n) \frac{\lambda([u,a_i])^2}{4\pi} \\
      &\quad= \int_0^{\rho} \int_{\mathbb S^1(a_i,r)} f_n(\lvert \D u\rvert) - f_n\left( \frac{\lambda([u,a_i])}{2\pi r} \right) \d\mathcal H^1\d r\\
           &\qquad+ \int_0^{\rho} 2 \pi r f_n\left( \frac{\lambda([u,a_i])}{2\pi r} \right) \d r- \frac{\lambda([u,a_i])^2}{2 \pi} \int_0^1 f_n(1/r) r \d r\\
           &\quad\doteq \mathrm{I}(n) + \mathrm{II}(n).
    \end{split}
  \end{equation} 
  For the first term, we will show that 
  \begin{equation}\label{eq:upperbound_I1_limit}
    \lim_{n \to \infty} \mathrm{I}(n) = \int_{0}^\rho \int_{\mathbb S^1(a_i,r)} \frac{|\D u|^2}{2} \d \mathcal H^1- \frac{\lambda([u,a_i])^2}{4\pi r}\d r
  \end{equation} 
  To see this, observe that for almost all $r \in (0,1)$ such that $\lvert\D u\rvert \in \LL^{2}(\partial\B(a_i,r))$, we have
  \begin{equation}
    \lim_{n \to \infty} \int_{\mathbb S^1(a_i,r)} f_n(\lvert\D u\rvert) - f_n\left( \frac{\lambda([u,a_i])}{2\pi r} \right) \d \mathcal H^1 = \int_{\mathbb S^1(a_i,r)} \frac12 \lvert \D u\rvert^2 \d\mathcal H^1 - \frac{\lambda([u,a_i])^2}{4 \pi r},
  \end{equation} 
  and using \eqref{eq:quad_upper2} we have 
  \begin{equation}\label{eq:whereweactuallyusecdec2}
    \begin{split}
    0 \leq &\int_{\mathbb S^1(a_i,r)} f_n(\lvert\D u\rvert) - f_n\left( \frac{\lambda([u,a_i])}{2\pi r} \right) \d \mathcal H^1 \\
           &\quad \leq \frac{f_{n,+}'(\max\{\lambda([u,a_i])/(2\pi r),t_0\})}{\max\{\lambda([u,a_i])/(2\pi r),t_0\}}\int_{\mathbb S^1(a_i,r)} \lvert \D u\rvert^2 - \frac{\lambda([u,a_i])^2}{4\pi^2r^2} \d\mathcal H^1 + 2 \pi C_2 r\\
           &\quad \leq  \left(\sup_{n \in \N} \frac{f_{n,+}'(t_0)}{t_0}\right)\left[ \int_{\mathbb S^1(a_i,r)} \frac12 \lvert \D u\rvert^2 \d\mathcal H^1 - \frac{\lambda([u,a_i])}{4 \pi r} \right]  + 2\pi C_2 r,
    \end{split}
  \end{equation} 
  using \ref{hyp:cdec2} and non-negativity of the integral in the last line, noting the supremum is finite by Lemma \ref{lem:fn_conv}\eqref{item:derivative_pointwise}.
  Thus, we have a uniform upper bound which is integrable on $(0,\rho)$ by Proposition \ref{prop:renormalisedenergy}\eqref{item:actualuseofren},
  from which \eqref{eq:upperbound_I1_limit} follows by the dominated convergence theorem.
  For the second term, we use Lemma \ref{lemma:entropy_error}, which gives
  \begin{equation}\label{eq:upperbound_I2_limit}
    \begin{split}
    \lim_{n \to \infty} \mathrm{II}(n) 
    &= \frac{\lambda([u,a_i])^2}{2\pi} \lim_{n \to \infty}\int_1^{\frac{2\pi \rho}{\lambda([u,a_i])}} f_n(1/r) r \d r \\
    &= \frac{\lambda([u,a_i])^2}{2\pi} \int_1^{\frac{2\pi \rho}{\lambda([u,a_i])}} \frac1{2r} \d r = \frac{1}{2}\log \frac{\lambda([u,a_i])}{2\pi\rho}.
  \end{split}
 \end{equation} 
  Combining \eqref{eq:upperbound_I1_limit} and \eqref{eq:upperbound_I2_limit} in \eqref{eq:upperbound_mainestimate}, the result follows.
\end{Proof}

Note that, instead of using \ref{hyp:cdec2} to establish \eqref{eq:whereweactuallyusecdec2}, one could have assumed that there exists $A,B>0$ such that for all $s,t \in \mathbb R_+$ and all $n \in \N$,
\(
	f_n(t) - f_n(s) \leq A\left(t^2 - s^2\right) + B.
\)

\subsection{Ball merging construction}\label{sec:ballmerging}

In this section, we will consider a Young function $f$ satisfying \ref{hyp:ndec2} and \ref{hyp:int}, and establish lower bounds for the functional $\int_{\Omega} f(\lvert \D u\rvert) \d x$ which will capture the singularities of $u$.
This will serve as the key step in our asymptotic analysis, for which we will follow the approach of \textsc{Jerrard} in \cite{Jerrard1999}, by identifying a suitable quantity $\Lambda_f(\cdot)$ to estimate the energy carried by the singularities.
In particular, this will lead to a lower bound estimate in the sequel (Proposition~\ref{prop:firstorderconve}\eqref{item:2ndorder}) which asymptotically matches the upper bound established in Proposition~\ref{prop:upperBound}.

For this we start with a fundamental lower bound for annuli, which we denote by
\[\A(a,\rho,\sigma) \doteq \{ x \in \mathbb R^2 : \rho \leq \lvert x - a\rvert \leq \sigma\} = \B(a,\rho)\setminus \B^{\circ}(a,\sigma)\]
where $a \in \mathbb R^2$ is the centre and the endpoints are given by $0 \leq \rho \leq \sigma$ and $\B^{\circ}$ denotes the open ball, recalling that balls $\B$ are always chosen to be closed.

\begin{lemma}\label{lem:annular_lowerbound}
  Let $u \in \WW^{1,2}(\A(a,\rho,\sigma); \mathcal N)$ with $a \in \mathbb R^2$ and $0 \leq \rho < \sigma < \frac{\sys}{2 \pi t_0}$ (with no restriction if $t_0 = 0$).
  If $E \doteq {\lambda(\Tr_{\partial \B_{\sigma}(a)}u(a + \sigma\,\cdot))^2}/{4\pi}$ is non-zero,  we have
  \begin{equation}\label{eq:lowerbound}
    \int_{\A(a,\rho,\sigma)} f(\lvert\D u\rvert) \,\d x 
    \geq E \left(\Lambda\left( \frac{\sigma}{E}\right) - \Lambda\left( \frac{\rho}{E} \right)\right) ,
  \end{equation} 
  where
  \begin{equation}\label{eq:lambda_general}
    \Lambda(t) = \Lambda_f(t) \doteq \frac12\sysN^2 \int_0^t f\left( \frac2{\sysN s} \right) s \d s,
  \end{equation} 
  which is finite by \ref{hyp:int}.
\end{lemma}
The quantity $\Lambda_f(\cdot)$ plays a similar role to the quantity $\Lambda_{\eps}(\cdot)$ considered by \textsc{Jerrard} in \cite[\S 3]{Jerrard1999} in the context of Ginzburg-Landau functionals.
Also setting $c_0 \doteq \frac12\sys$, we have
\begin{equation}\label{eq:defOfLambda}
  \Lambda(t) = 2c_0^2 \int_0^t f\left( \frac1{c_0s} \right) s\d s.
\end{equation} 
We will often use Lemma \ref{lem:annular_lowerbound} with $\rho=0$, however the general form is also used for instance in \eqref{eq:whereweusethelemma}.
The result also holds when $E=0$ if we interpret the right hand side of \eqref{eq:lowerbound} to be zero in that case.

\begin{Proof}{Lemma}{lem:annular_lowerbound}
  By converting to radial coordinates and using Jensen's inequality, we have
  \begin{equation}
    \begin{split}
      \int_{\mathbb A(a,\rho,\sigma)} f(\lvert\D u\rvert) \,\d x 
    &= 2\pi \int_{\rho}^{\sigma} \int_{\mathbb S^1} f\left(\lvert \D u(a + r \omega)\rvert \right) \d\mathcal H^1(\omega)\, r\d r \\
    &\geq 2\pi \int_{\rho}^{\sigma} f\left( \frac{\lambda}{ 2 \pi r} \right) r \d r,
    \end{split}
  \end{equation} 
  where $\lambda = \lambda([u,a])$ and  where we have used that $u \in \WW^{1,2} (\mathbb A(a,\rho,\sigma); \mathcal N)$ in this annulus, so $\lambda=\lambda(\Tr_{\mathbb S^1}u(a+r \,\cdot))$ for all $r \in [\rho,\sigma]$.
  Note that, since $\lambda \neq 0$, we have $\lambda \geq \sysN$.
  Then making the substitution $s = {r}/E$ we write
  \begin{equation}
    \begin{split}
      \int_{\mathbb A(a,\rho,\sigma)} f(\lvert\D u\rvert) \,\d x 
      &\geq 2\pi \int_{\rho}^{\sigma} f\left( \frac{\sqrt{E}}{\sqrt{\pi} r} \right) r \d r  \\
      &= 2\pi E^2 \int_{\frac{\rho}E}^{\frac{\sigma}E} f\left( \frac1{\sqrt{\pi E} s} \right) s \d s  \\
      &\geq \frac{\sysN^2}2 E \int_{\frac{\rho}E}^{\frac{\sigma}E} f\left( \frac{2}{\sysN} \frac1s \right) s\d s \\
      &= E \left( \Lambda(\frac{\sigma}E) - \Lambda(\frac{\rho}E)\right)
    \end{split}
  \end{equation} 
  where we used \ref{hyp:ndec2} noting that $\frac1{\sqrt{\pi E}s} \geq \frac{\sys}{2 \pi \sigma} \geq t_0$ for all $s \leq \frac{\sigma}E$, by the upper bound on $\sigma$.
\end{Proof}

\begin{lemma}\label{lem:Lambda_mean}
  Let $f$ be a Young function, and let $\Lambda(t)$ be as in \eqref{eq:lambda_general}.
  Then for $0 \leq \rho < \sigma$, we have
  \begin{equation}
    E  \mapsto  E\,\left(\Lambda\left(\frac{\sigma}E\right)-\Lambda\left(\frac{\rho}E\right)\right)
  \end{equation} 
  is non-decreasing in $E>0$.
\end{lemma}

Lemma \ref{lem:Lambda_mean} is analogous to the mean-value property observed by \textsc{Jerrard} \cite[Prop.~3.1, (3.3)]{Jerrard1999},
and by taking $\rho = 0$, we see that $t \mapsto \frac1t\,\Lambda(t)$ is also non-increasing in $t$.
In addition, by combining Lemmas \ref{lem:annular_lowerbound} and \ref{lem:Lambda_mean} we infer that \eqref{eq:lowerbound} holds for any $E \in [0, {\lambda(\Tr_{\mathbb S^1}u(a + \sigma\,\cdot))^2}/{4\pi}]$.  
\begin{Proof}{Lemma}{lem:Lambda_mean}
  For $E_1 < E_2$ we use the change of variable $r = \frac{E_1}{E_2}s$ in
  \begin{equation}
    \begin{split}
      E_1 \left( \Lambda\left(\frac{\sigma}E_1\right)- \Lambda\left(\frac{\rho}E_1\right)\right) 
      &= 2c_0^2 E_1 \int_{{\rho}/{E_1}}^{{\sigma}/{E_1}} f\left( \frac1{c_0s} \right) s\,\d s\\
      &= 2c_0^2 E_2 \int_{{\rho}/{E_2}}^{{\sigma}/{E_2}} f\left( \frac{E_1}{E_2}\frac1{c_0r} \right) \frac{E_2}{E_1} r\,\d r\\
      &\leq 2c_0^2 E_2 \int_{\rho/E_2}^{\sigma/E_2} f\left( \frac1{c_0r} \right) r \,\d r,
    \end{split}
  \end{equation} 
  where we used the fact that $f(st) \leq s f(t)$ with $s = \frac{E_1}{E_2} \in (0,1)$ in the last line, which follows from convexity of $f$ and the fact that $f(0)=0$.
\end{Proof}

We record the following ball merging construction, which can be found in \cite[Lem.~3.1]{Jerrard1999} (see also \cite[Fig.~1]{Sandier1998}).

\begin{lemma}[Merging lemma]\label{lem:merging_lemma}
  Given any finite collection $\mathcal{B}$ of closed balls in $\mathbb R^2$, we can associate
   a pairwise disjoint collection of balls $\widetilde{\mathcal{B}}$ such that
\(
  \bigcup \mathcal B \subset \bigcup \widetilde{\mathcal B}\) and \( \sum_{\B \in \mathcal B} r(\B) = \sum_{\tilde \B \in \widetilde{\mathcal B}} r(\tilde \B).
\) 
Moreover, the construction is such that for each $\tilde{\B} \in \widetilde{\mathcal B}$, there exist finitely many $\B_1,\dotsc,\B_k \in \mathcal B$ such that
\(
  \bigcup_{i=1}^k \B_i \subset \tilde{\B}\) and \(\sum_{i=1}^k r(\B_i) = r(\tilde{\B}).
\) 
In particular, if $k=1$ we have $\B_1 = \tilde \B$.
\end{lemma}

Informally, Lemma \ref{lem:merging_lemma} provides a way to systematically ``disjointify'' a given collection of balls, by merging any two balls which intersect.
In the sequel, we will say a ball $\tilde\B$ can be \emph{obtained by merging} balls $\B_1,\dotsc, \B_k$ if we have $\B_i \subset \tilde\B$ for all $i = 1, \dotsc, m$, and $r(\tilde\B) = \sum_{i=1}^m r(\B_i)$.

\begin{Proof}{Lemma}{lem:merging_lemma}
  We inductively merge any two balls $\B_1 = \B_{\sigma_1}(a_1)$ and $\B_2 = \B_{\sigma_2}(a_2)$ which intersect, by replacing them by $ \B_{\sigma}(a)$, where
  \( a = \frac{\sigma_1 a_1 + \sigma_2 a_2}{\sigma_1+\sigma_2} \) and \(  \sigma = \sigma_1+\sigma_2. \)
  We repeat this procedure for any two intersecting balls until they are all pairwise disjoint.
\end{Proof}
\begin{prop}[Ball merging for {$\mathrm{E}^{1,2}_{\sg}$}]\label{prop:genmergingballlemma} Let $f$ be a Young function satisfying \ref{hyp:ndec2} and \ref{hyp:int}, and let $\Lambda(t)$ be as in \eqref{eq:lambda_general}.
  For $\delta>0$, let $u \in \RR^{1,2}_{\delta}(\Omega;\mathcal N)$ have singularities at distinct $a_1,\dots,a_k \in \Omega $ with $k \geq 1$,
  and set $\eta_0 = \min\big\{\frac{\sys}{2\pi t_0},\delta\big\}>0$. 
  
  Then for each $\eta \in (0,\eta_0]$ there exists a finite collection $\mathcal{B}(\eta)$ of disjoint closed balls in $\Omega$ whose sum of radii is $\eta$, such that $\bigcup\mathcal B(\eta)$ contains $\{a_i\}_{i=1}^k$ and that
  \begin{equation}\label{eq:lowerboundinprop}
     \Big(\sum_{\B \in \mathcal B(\eta)} \Esgnu(\Tr_{\partial \B}u)\Big) \, \Lambda_f\Big(\frac{\eta}{\sum_{\B \in \mathcal B(\eta)}\Esgnu(\Tr_{\partial \B}u) } \Big) \leq \int_{\bigcup \mathcal B(\eta)} f(\lvert \D u\rvert ) \,\d x,
  \end{equation} 
  understanding the left-hand side to be zero when $u$ has no singularities in $\Omega$.
\end{prop}

Note that for the collection $\mathcal B(\eta)$ obtained by Proposition \ref{prop:genmergingballlemma},  $(\tr_{\partial \B}u)_{\B \in \mathcal B(\eta)}$ is topological resolution of $\tr_{\partial \Omega}u$ for each $\eta$, since $\{a_i\}_{i=1,\dots,k} \subset \bigcup\mathcal B(\eta)$ and the balls are disjoint.
The key step in proving Proposition \ref{prop:genmergingballlemma} will be to show the energy inequality \eqref{eq:lowerboundinprop} is preserved under a suitable merging process; we will record this as a separate lemma.

\begin{lemma}\label{lem:merging_energy_bound}
  Suppose $(\sigma_i)_{i=1}^k$, $(E_i)_{i=1}^k$ are positive numbers and let
  \begin{equation}\label{eq:useless}
    0<E_{\Sigma} \leq \sum_{i=1}^k E_i, \quad \sigma_{\Sigma} = \sum_{i=1}^k \sigma_i.
  \end{equation}
  Suppose there is a parameter $T>0$ such that
  \begin{align}
    \sigma_i &\leq E_i T \quad \text{ for all } 1 \leq i \leq k,\label{eq:sigmai_inequality}\\
    \sigma_{\Sigma} &\geq E_{\Sigma} T \label{eq:sigmasum_inequality},
  \end{align} 
  then we have
  \begin{equation}
    \sum_{i=1}^k E_i \Lambda\left(\frac{\sigma_i}{E_i}\right) \geq E_{\Sigma} \Lambda\left( \frac{\sigma}{E_{\Sigma} }\right) .
  \end{equation} 
\end{lemma}

In the proof of Lemma \ref{lem:merging_energy_bound}, we use Lemma \ref{lem:Lambda_mean} with $\rho = 0$ twice.

\begin{Proof}{Lemma}{lem:merging_energy_bound}
  Since \eqref{eq:sigmai_inequality} gives $\frac{\sigma_i}{E_i} \leq T$ for all $i$, and since $E \mapsto E \Lambda_f(\sigma/E)$ is non-increasing by Lemma~\ref{lem:Lambda_mean} (with $\rho=0$), we have
  \begin{equation}\label{eq:intermediate_lower}
    \sum_{i=1}^k E_i \Lambda_f\left( \frac{\sigma_i}{E_i} \right) \geq \sum_{i=1}^k \frac{\sigma_i}T \Lambda_f\left( T \right) = \frac{\sigma_{\Sigma}}{T} \Lambda_f\left( T \right),
  \end{equation} 
  where the last equality follows by $\eqref{eq:useless}_2$.
  Using  Lemma \ref{lem:Lambda_mean} and with \eqref{eq:sigmasum_inequality}, we see that
  \begin{equation}
    \frac{\sigma_{\Sigma}}{T} \Lambda_f\left( T \right) \geq E_{\Sigma} \Lambda_f\left( \frac{\sigma_{\Sigma}}{E_{\Sigma}} \right) ,
  \end{equation} 
  which we combine with \eqref{eq:intermediate_lower} to conclude.
\end{Proof}

\begin{Proof}{Proposition}{prop:genmergingballlemma}
  Fix $\eta \in (0,\eta_0]$.
  Given a parameter $t>0$, we will inductively define two collections of closed balls $\mathcal{B}^+(t), \mathcal{B}^-(t)$ which satisfy the following properties for $t>0$:
  \begin{enumerate}[(i)]
  	   \item\label{item:merging_disjoint} 
  	The collections $\mathcal B^{\pm}(t)$ are both finite collections of balls, such that balls in $\mathcal B^+(t)$ are pairwise disjoint, and balls in $\mathcal B^-(t)$ have pairwise disjoint interiors.
  	That is,
  	\begin{align*}
  		\overline \B_1 \cap \overline \B_2 &= \emptyset \quad \mbox{for all } \B_1, \B_2 \in \mathcal B^+(t),\\
  		\B_1^{\circ} \cap \B_2^{\circ} &= \emptyset \quad \mbox{for all } \B_1, \B_2 \in \mathcal B^-(t).
  	\end{align*} 
    Here, for clarity, we write $\overline \B$ and $\B^{\circ}$ for the closure and interior of $\B$ respectively.

  \item\label{item:merging_preservation} The collections $(c_{\B},\Tr_{\partial \B}u)_{\B \in \mathcal B^+(t)}$ and $(c_{\B},\Tr_{\partial \B}u)_{\B \in \mathcal B^-(t)}$ are topological resolutions of $\Tr_{\partial \Omega}u$, where $c_{\B}$ denotes the centre of the ball $\B$. We will say, by abuse of language, that both $\mathcal B^+(t), \mathcal B^-(t)$ are topological resolutions of $u$ in $\Omega$. In particular,
      \begin{equation*}
        \sum_{\B \in \mathcal B^{\pm}(t)} \mathrm{E}^{1,2}_{\sg}(\Tr_{\partial \B}u) \geq \mathrm{E}^{1,2}_{\sg}(\Tr_{\partial\Omega}u) \quad\mbox{for all $t>0$.}
      \end{equation*} 
      Here we use the notation $\mathcal B^{\pm}(t)$ to mean the property should hold with both $\mathcal B^+(t)$ and $\mathcal B^-(t)$ in place of $\mathcal B^{\pm}(t)$.

   \item\label{item:merging_containment} The collection $\mathcal{B}^+(t)$ can be obtained by applying Lemma \ref{lem:merging_lemma} to $\mathcal{B}^-(t)$.
     That is, to each $\tilde \B \in \mathcal B^+(t)$ there exists $\B_1, \dotsc, \B_m \in \mathcal B^-(t)$ such that $\tilde \B$ can be obtained by merging $\B_1, \dotsc, \B_m$.

    \item\label{item:merging_radius} We have the radius bounds
      \begin{align*}
        r(\B) &\leq  \mathrm E^{1,2}_{\sg}(\Tr_{\partial \B}u)\,t \quad \mbox{for all } \B \in \mathcal B^-(t),\\
        r(\B) &\geq  \mathrm E^{1,2}_{\sg}(\Tr_{\partial \B}u)\,t \quad \mbox{for all } \B \in \mathcal B^+(t),
      \end{align*} 
      where $r(\B)$ denotes the radius of the ball $\B$.

    \item\label{item:merging_energy} We have the lower bound
      \begin{equation*}
      \int_{\B}f(|\D u|) \d x \geq \mathrm{E}_{\sg}^{1,2}(\Tr_{\partial \B}u)\,\Lambda_f\left(\frac{r(\B)}{\mathrm{E}_{\sg}^{1,2}(\Tr_{\partial \B}u)}\right) \quad\mbox{
      for each $\B \in \mathcal B^{\pm}(t)$.}
      \end{equation*} 

    \item\label{item:merging_terminate}  The sum of radii satisfies
      \begin{equation*}
        r(\mathcal B^{\pm}(t)) \doteq \sum_{\B \in \mathcal B^{\pm}(t)} r(\B) \leq \eta.
      \end{equation*} 
  \end{enumerate}
  We will describe this process in Steps \hyperlink{step1}{1}--\hyperlink{step4}{4}.
  The idea will be to ``expand'' the balls with respect to the parameter $t$, where we consider concentric balls $\B(t)$ following the expansion law
  \begin{equation}\label{eq:growth_law}
    r(\B(t)) = \mathrm{E}_{\mathrm{sg}}^{1,2}(\Tr_{\partial \B(t)}u) \, t.
  \end{equation} 
  We cannot do this indefinitely without having balls intersect however, in which case we will need to merge balls using Lemma \ref{lem:merging_lemma}.
  To ensure the lower bound \eqref{item:merging_energy} is preserved under this process, we use Lemma \ref{lem:merging_energy_bound} on specific balls; for this reason, we have two collections of balls $\mathcal{B}^{\pm}(t)$ at each time-step satisfying \eqref{item:merging_radius}.

  In what follows we will abbreviate
  \begin{equation}
    \mathrm{E}^{1,2}_{\sg}( a,\rho) \doteq \mathrm{E}^{1,2}_{\sg}([\Tr_{\partial \B(a,\rho)}u]),
  \end{equation} 
  and also put $\mathrm E^{1,2}_{\sg}(a) \doteq \mathrm{E}^{1,2}_{\sg}([u,a])$, where $[u,a]$ denotes a limiting geodesic of $u$ at $a$.

  \subsubsection*{Step 1}\hypertarget{step1}{} The initial parametrisation.

  For $t>0$ sufficiently small, we define $\mathcal{B}^+(t)$ and $\mathcal{B}^-(t)$ to be equal as the collection of balls $\B(a_i,\rho_i(t))$ for each $1 \leq i \leq k$, where 
  \begin{equation}\label{eq:initial_radii_param}
    \rho_i(t) = \mathrm E^{1,2}_{\mathrm{sg}}(a_i) t.
  \end{equation} 
  Here $t>0$ is chosen so these balls do not intersect; that is $t<t_1$ with
  \begin{equation}\label{eq:first_time_step}
    t_{1} = \sup\{ t>0 : \mathcal{B}^{+}(t) \text{ is pairwise disjoint} \},
  \end{equation} 
  noting that $t_1>0$ since the starting points $a_1,\ldots,a_k$ were distinct.
  This defines $\mathcal{B}^{\pm}(t)$ for all $0<t<t_1$, which by definition satisfies \eqref{item:merging_preservation}--\eqref{item:merging_radius}.
  Since $u$ is only singular at each $a_i \in \B(a_i,\rho_i(t))$, by Lemma \ref{lem:annular_lowerbound} we have
  \begin{equation}
    \int_{\B(a_i,\rho_i(t))} f(\lvert \D u\rvert ) \d x \geq \mathrm{E}^{1,2}_{\sg}(a_i,\rho_i(t)) \Lambda_f \left( \frac{\rho_i(t)}{\mathrm{E}^{1,2}_{\sg}(a_i,\rho_i(t))} \right),
  \end{equation} 
  establishing \eqref{item:merging_energy}.
  Lastly, if we have that \eqref{item:merging_terminate} is satisfied with equality we will terminate the process and proceed to Step \hyperlink{step5}{5}; otherwise, our collection $\mathcal{B}^{\pm}(t)$ satisfies all the claimed properties for $0<t<t_{1}$, and so we can move to Step \hyperlink{step2}{2}.

  \subsubsection*{Step 2}\hypertarget{step2}{} Evolving and frozen balls.

  We will split $\mathcal{B}^{\pm}(t)$ further into disjoint unions
  \begin{equation}
    \mathcal{B}^{+}(t) = \mathcal{B}^{+}_{\mathrm{F}}(t) \cup \mathcal{B}^{+}_{\mathrm{G}}(t), 
    \quad \mathcal{B}^{-}(t) = \mathcal{B}^{-}_{\mathrm{F}}(t) \cup \mathcal{B}^{-}_{\mathrm{G}}(t)
  \end{equation} 
  distinguishing between the ``frozen'' and ``growing'' balls respectively.
  Note that none of the balls are initially frozen, so we define $\mathcal{B}^+_{\mathrm{G}}(t) = \mathcal{B}^-_{\mathrm{G}}(t) \doteq \mathcal{B}^+(t)$ and $\mathcal{B}^+_{\mathrm{G}}(t) = \mathcal{B}^-_{\mathrm{G}}(t) = \emptyset$ for all $t \in (0,t_1)$.

  We will inductively show there are finitely many  $0 = t_0 < t_1 < \cdots$ such that
  \begin{enumerate}[(a)]
    \item\label{item:step_2_case_a} $\mathcal{B}^+_{\mathrm{G}}(t) = \mathcal{B}^-_{\mathrm{G}}(t)$ for $t \in (t_{j-1},t_{j})$ and each ball $\B=\B(a,\sigma) \in \mathcal B^+_{\mathrm{G}}(t_{j-1})$ evolves as $\B(t)$ with respect to the expansion law \eqref{eq:growth_law}; that is, $\B(t) = \B(a, \sigma(t))$ where $\sigma(t) = \Esgnu(\tr_{\partial\B}u)t$.
      In particular, balls in $\mathcal{B}^{\pm}_{\mathrm{G}}$ satisfy \eqref{item:merging_radius} with equality.

    \item\label{item:step_2_case_b} both $\mathcal{B}^{\pm}_{\mathrm{F}}(t)$ are constant on each $(t_{j-1},t_{j})$ and $\mathcal{B}^+_{\mathrm{F}}(t)$ can be obtained by applying Lemma \ref{lem:merging_lemma} to $\mathcal{B}^-_{\mathrm{F}}(t)$.
    Moreover, we require that balls $\B \in \mathcal{B}^+_{\mathrm{F}}(t)$ satisfy \eqref{item:merging_radius} with a strict inequality.
	\end{enumerate}

  In Step \hyperlink{step1}{1}, we verified this holds up to $t \leq t_1$ defined in \eqref{eq:first_time_step}.
  Given $\mathcal B^{\pm}(t_{j-1})$, we can define $\mathcal{B}^{\pm}(t)$ via \eqref{item:step_2_case_a} and \eqref{item:step_2_case_b} until we reach $t \doteq t_{j} > t_{j-1}$ such that one of the following holds.

  \emph{Case 1}\hypertarget{case1}{} (collision): We have $\mathcal{B}^+(t)$ is no longer pairwise disjoint, in which case we define \eqref{eq:tj_minus} as below and move to Step \hyperlink{step3}{3}.

  \emph{Case 2}\hypertarget{case2}{} (saturation): There is $\B \in \mathcal{B}^+_{\mathrm{F}}(t)$ which satisfies \eqref{item:merging_radius} with equality, in which case we define \eqref{eq:tj_minus} as below and move to Step \hyperlink{step4}{4}.

  \emph{Case 3}\hypertarget{case3}{} (termination): We have $r(\mathcal B^+(t)) = \eta$, in which we terminate the process and move to Step \hyperlink{step5}{5}.

  If either Case \hyperlink{case1}{1} or \hyperlink{case2}{2} occurs, at this time step $t =t_j$ we let
  \begin{equation}\label{eq:tj_minus}
    \mathcal{B}^{\pm}(t_j-) = \mathcal{B}_{\mathrm{F}}^{\pm}(t_j-) \cup \mathcal{B}_{\mathrm{G}}^{\pm}(t_j-)
  \end{equation} 
  be the collections obtained by applying \eqref{item:step_2_case_a} and \eqref{item:step_2_case_b} at this time step, which we will suitably modify in the subsequent steps.
  If multiple cases occur at once, we apply the corresponding steps in order. 
  This procedure is shown in Figure \ref{fig:merging_diagram}.

\begin{figure}[h]
	\centering
	\begin{tikzcd}
		&& \begin{array}{c} \begin{matrix}\text{Collision} \\ \text{(Step \hyperlink{step3}{3})}\end{matrix} \end{array} \\
		\begin{array}{c} \begin{matrix}\text{Initialisation} \\ \text{(Step \hyperlink{step1}{1})}\end{matrix} \end{array} & \begin{array}{c} \begin{matrix}\text{Expansion} \\ \text{(Step \hyperlink{step2}{2})}\end{matrix} \end{array} \\
		&& \begin{array}{c} \begin{matrix}\text{Saturation} \\ \text{(Step \hyperlink{step4}{4})}\end{matrix} \end{array} \\
		 &\begin{array}{c} \begin{matrix}\text{Termination} \\ \text{(Step \hyperlink{step5}{5})}\end{matrix} \end{array}\\
		\arrow["{\small\substack{\text{Increment} \\ t_j \mapsto t_{j+1}}}", bend left = 10, from=1-3, to=2-2]
		\arrow["{\text{Case \hyperlink{case2}{2}}}",bend left = 10, dashed, from=1-3, to=3-3]
		\arrow["{\text{Case \hyperlink{case3}{3}}}",bend left = 10, dashed, from=3-3, to=4-2]
		\arrow[from=2-1, to=2-2]
		\arrow["{\text{Case \hyperlink{case1}{1}}}", bend left=10, from=2-2, to=1-3]
		\arrow["{\text{Case \hyperlink{case2}{2}}}"', bend right = 10, from=2-2, to=3-3]
		\arrow["{\text{Case \hyperlink{case3}{3}}}"',from=2-2, to=4-2]
		\arrow["{\small\substack{\text{Increment} \\ t_j \mapsto t_{j+1}}}"', bend right = 10, from=3-3, to=2-2]
	\end{tikzcd}
	\vspace{-1.2cm}
	\caption{Flowchart of the inductive process}
	\label{fig:merging_diagram}
	\vspace{-0.3cm}
\end{figure}
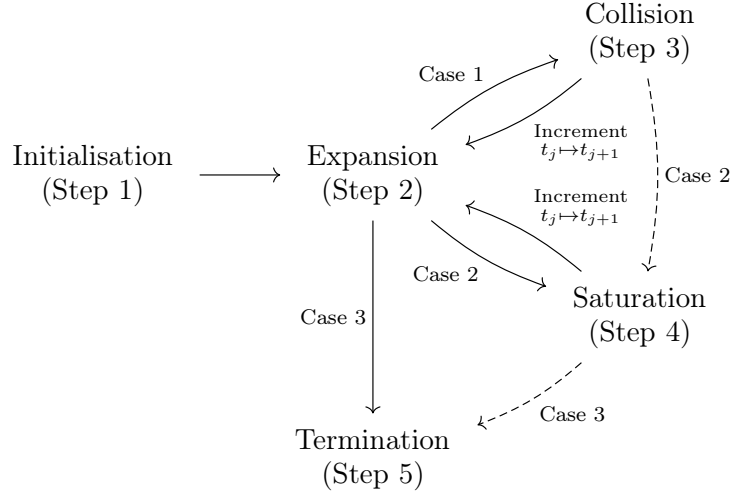

We claim that if $\mathcal{B}^{\pm}(t)$ satisfies \eqref{item:step_2_case_a}, \eqref{item:step_2_case_b} for all $t \in (t_{j-1},t_j)$, then the claimed properties \eqref{item:merging_disjoint}--\eqref{item:merging_terminate} hold for all $t$ in this range.
Indeed, property \eqref{item:merging_disjoint} holds for $t<t_j$ by Case \hyperlink{case1}{1} of $t_{j}$ and by the fact that the only expanding balls are those in $\mathcal{B}^+_{\mathrm{G}}(t) = \mathcal{B}^-_{\mathrm{G}}(t)$.
  Similarly Case \hyperlink{case2}{2} of $t_{j}$ ensures \eqref{item:merging_radius} holds provided $t<t_j$, noting the inequality for balls in $\mathcal{B}^-_{\mathrm{F}}(t)$ weakens as $t$ is increased.
Also since $\bigcup\mathcal{B}^{\pm}(t)$ still contains the initial singularities $a_1,\dotsc,a_k$, $\mathcal B^{\pm}(t)$ remains topological resolutions of $u$ in $\Omega$ and satisfies \eqref{item:merging_preservation}.
  By induction hypothesis \eqref{item:merging_containment} holds at $t=t_{j-1}$, and hence is preserved for $t \in (t_{j-1},t_j)$ since there are no new intersections.
  Also Case \hyperlink{case3}{3} ensures \eqref{item:merging_terminate} holds.

  For \eqref{item:merging_energy}, this is satisfied by balls in $\mathcal B^\pm_{\mathrm F}(t)$ by inductive assumption, so we only need to check the balls in $\mathcal B^{\pm}_{\mathrm{G}}(t)$. For this let $\B(t) = \B(a,\sigma(t)) \in \mathcal{B}^+_{\mathrm{G}}(t) = \mathcal{B}^-_{\mathrm{G}}(t)$ with $\sigma(t) = \Esgnu(\Tr_{\partial \B}u)t$, which corresponds to $\B=\B(a,\sigma(t_{j-1})) \in \mathcal B^+_{\mathrm{G}}(t_{j-1})$. 
  Then by applying Lemma \ref{lem:annular_lowerbound} (annular lower bound) and Lemma \ref{lem:Lambda_mean} (the mean value property), we can estimate 
  	\begin{IEEEeqnarray}{rCl}
      \label{eq:whereweusethelemma}  \int_{\B(a,\sigma(t))} f(\lvert \D u\rvert ) \d x &\geq& \int_{\B(a,\sigma(t))\setminus \B(a,\sigma(t_{j-1}))} f(\lvert \D u\rvert)\d x + \int_{\B(a,\sigma(t_{j-1}))} f(\lvert \D u\rvert)\d x\\
      &\geq&\notag \mathrm{E}^{1,2}_{\sg}(a,\sigma(t)) \left( \Lambda_f\left( \frac{\sigma(t)}{\mathrm{E}^{1,2}_{\sg}(a,\sigma(t))} \right) - \Lambda_f\left( \frac{\sigma(t_{j-1})}{\mathrm{E}^{1,2}_{\sg}(a,\sigma(t))} \right)  \right)  \\
      &&\notag+ \mathrm{E}^{1,2}_{\sg}(a,\sigma(t_{j-1})) \Lambda_f\left( \frac{\sigma(t)}{\mathrm{E}^{1,2}_{\sg}(a,\sigma(t_{j-1}))} \right)\\
      &=&\notag \mathrm{E}^{1,2}_{\sg}(a,\sigma(t)) \Lambda_f\left( \frac{\sigma(t)}{\mathrm{E}^{1,2}_{\sg}(a,\sigma(t))} \right) ,
  \end{IEEEeqnarray} 
  since $\mathrm{E}^{1,2}_{\sg}(a,\sigma(t)) = \mathrm{E}^{1,2}_{\sg}(a,\sigma(t_{j-1}))$. This verifies \eqref{item:merging_energy}, so we infer that \eqref{item:merging_disjoint}--\eqref{item:merging_terminate} holds for all $t \in (t_{j-1},t_j)$.

  \subsubsection*{Step 3}\hypertarget{step3}{} Ball merging I (Collision).

  Suppose we are in Case \hyperlink{case1}{1} of the previous step, so $\mathcal{B}^+(t_j-)$ is not pairwise disjoint.
  We will apply Lemma \ref{lem:merging_lemma} to $\mathcal{B}^+(t_j-)$ and let $\mathcal{B}^+_{\mathrm{New}}(t_j)$ denote the collection of \emph{new} balls obtained in this process, \textit{i.e.}\,those which strictly contains one of the existing balls.

  Then, each $\B \in \mathcal{B}^+_{\mathrm{New}}(t_j)$ arises from $\B_1^{\mathrm{F}},\dotsc,\B_{m}^{\mathrm{F}} \in \mathcal B^+_{\mathrm{F}}(t_j-)$ and $\B_1^{\mathrm{G}},\dotsc,\B_{\ell}^{\mathrm{G}} \in \mathcal B^+_{\mathrm{G}}(t_j-)$ which are contained in $\B$.
  Note that by \eqref{item:merging_containment}, there also exists $\widetilde\B_1^{\mathrm{F}}, \dotsc, \widetilde\B_{\tilde m}^{\mathrm{F}} \in \mathcal B^-_{\mathrm{F}}(t_j-)$ from which we can obtain $\B_1^{\mathrm{F}},\dotsc,\B_{m}^{\mathrm{F}}$ by merging.
  Then, since by \eqref{item:merging_radius} we have
  \begin{align}
    \label{eq:merging_fballs} r(\B_i^{\mathrm{F}}) &\geq \Esgnu(\Tr_{\partial \B_i^{\mathrm{F}}}u) t_j \quad\mbox{for all $i=1,\dotsc,m$,}\\
    \label{eq:merging_tildefballs} r(\widetilde\B_i^{\mathrm{F}}) &\leq \Esgnu(\Tr_{\partial \widetilde\B_i^{\mathrm{F}}}u) t_j \quad\mbox{for all $i=1,\dotsc,\tilde m$,}\\
    \label{eq:merging_gballs}r(\B_i^{\mathrm{G}}) &= \Esgnu(\Tr_{\partial \B_i^{\mathrm{G}}}u) t_j\quad\mbox{for all $i=1,\dotsc,\ell$,}
  \end{align} 
  it follows from \eqref{eq:merging_tildefballs} and \eqref{eq:merging_gballs} that 
  \begin{equation}
    \begin{split}\label{eq:merged_radius_energy}
      r(\B) &= \sum_{i=1}^m r(\B_i^{\mathrm{F}}) + \sum_{i=1}^{\ell} r(\B_i^{\mathrm{G}}) =\sum_{i=1}^{\tilde m} r(\widetilde\B_i^{\mathrm{F}}) + \sum_{i=1}^{\ell} r(\B_i^{\mathrm{G}})\\
      &\geq t_j\sum_{i=1}^m \Esgnu(\Tr_{\partial \B_i^{\mathrm{F}}}u) + t_j\sum_{i=1}^{\ell} \Esgnu(\Tr_{\partial \B_i^{\mathrm{G}}}u) \geq \Esgnu(\Tr_{\partial \B}u)t_j,
   \end{split}
  \end{equation} 
  where we used the subadditivity of the singular energy (namely \eqref{eq:esg_subadditive}) in the last line.
  Then since $\{\B_i^{\mathrm{F}}\}_{i=1}^m$, $\{\B_i^{\mathrm{G}}\}_{i=1}^\ell$ satisfies \eqref{item:merging_energy} by inductive assumption, by Lemma \ref{lem:merging_energy_bound} using \eqref{eq:merging_fballs} and \eqref{eq:merging_gballs}, we have
  \begin{equation}
    \begin{split}
    \int_{\B} f(\lvert \D u\rvert) \d x 
    &\geq \sum_{i=1}^m \int_{\B^{\mathrm{F}}_i} f(\lvert \D u\rvert) \d x + \sum_{i=1}^{\ell} \int_{\B^{\mathrm{G}}_i} f(\lvert \D u\rvert) \d x \\
    &\geq \sum_{i=1}^m \Esgnu(\tr_{\partial \B_i^{\mathrm{F}}}u) \Lambda_f\left( \frac{r(\B^{\mathrm{F}}_i)}{\Esgnu(\tr_{\partial \B_i^{\mathrm{F}}}u)} \right)  \\
    &\qquad+ \sum_{i=1}^\ell \Esgnu(\tr_{\partial \B_i^{\mathrm{G}}}u) \Lambda_f\left( \frac{r(\B^{\mathrm{G}}_i)}{\Esgnu(\tr_{\partial \B_i^{\mathrm{G}}}u)} \right) \\
    &\geq \Esgnu(\Tr_{\partial \B}u) \Lambda_f\left(\frac{r(\B)}{\Esgnu(\Tr_{\partial \B}u)}  \right),
    \end{split}
  \end{equation} 
  thereby establishing \eqref{item:merging_energy} for $\B$.

  We then define $\mathcal{B}^{\pm}_{\mathrm{F}}(t_j)$ and $\mathcal{B}^{\pm}_{\mathrm{G}}(t_j)$ as $\mathcal{B}^{\pm}_{\mathrm{F}}(t_j-)$ and $\mathcal{B}^{\pm}_{\mathrm{G}}(t_j-)$ respectively, except we make the following modifications corresponding to each $\B \in \mathcal{B}_{\mathrm{New}}^+(t_j)$:
  \begin{itemize}
    \item We remove $\B_1^{\mathrm{F}}, \dotsc, \B_m^{\mathrm{F}}$ from $\mathcal B_{\mathrm{F}}^+(t_j)$ and instead put $\B \in \mathcal B_{\mathrm{F}}^+(t_j)$.

    \item We remove $\B_1^{\mathrm{G}}, \dotsc, \B_{\ell}^{\mathrm{G}}$ from $\mathcal B_{\mathrm{G}}^{\pm}(t_j)$ and put them in $\mathcal{B}_{\mathrm{F}}^-(t_j)$.
  \end{itemize}
  The above shows that \eqref{item:merging_radius} and \eqref{item:merging_energy} hold, and \eqref{item:merging_containment} holds since each $\B \in \mathcal{B}_{\mathrm{new}}^+(t_j)$ can be obtained by merging $\widetilde\B_1^{\mathrm{F}}, \dotsc, \widetilde\B_{\tilde m}^{\mathrm{F}}, \B_1^{\mathrm{G}}, \dotsc \B_\ell^{\mathrm{G}}$, which now lie in $\mathcal B_{\mathrm{F}}^-(t_j)$.
  The remaining properties are a routine check.

  Thus we have defined $\mathcal{B}^{\pm}(t_j)$.
  Now if \eqref{eq:merged_radius_energy} holds with equality for any of the $\B \in \mathcal{B}^+_{\mathrm F}(t_j)$, since Case \hyperlink{case2}{2} would be violated, we rename our new collections as $\mathcal{B}^{\pm}(t_j-)$ and move to Step \hyperlink{step4}{4}.
  Otherwise, we return to Step \hyperlink{step2}{2} to continue the inductive process.

  \subsubsection*{Step 4}\hypertarget{step4}{} Ball merging II (Saturation).

  Suppose there is $\B \in \mathcal{B}^+_{\mathrm{F}}(t_j-)$ satisfying \eqref{item:merging_radius} with equality, that is
  \begin{equation}\label{eq:equality_case}
    r(\B) = \Esgnu(\Tr_{\partial \B}u) t_j.
  \end{equation} 
  If this holds, let $\B_1,\dotsc,\B_{\ell} \in \mathcal{B}^-_{\mathrm{F}}(t_j-)$ be the balls which merge to give $\B$, which exist by \eqref{item:merging_containment}.
  We then define $\mathcal{B}^{\pm}(t_j)$ as $\mathcal{B}^{\pm}(t_j-)$ except that for any such ball $\B$ we make the following modifications:
  \begin{itemize}
    \item We move $\B$ from $\mathcal{B}_{\mathrm{F}}^+(t_j)$ to $\mathcal{B}_{\mathrm{G}}^+(t_j)$.
    \item We remove $\B_1,\dotsc,\B_{\ell}$ from $\mathcal{B}_{\mathrm{F}}^-(t_j)$ and add $\B$ to $\mathcal{B}_{\mathrm{G}}^-(t_j)$.
  \end{itemize}

	We show that \eqref{item:merging_disjoint}--\eqref{item:merging_terminate} holds. 
  Since we apply Step \hyperlink{step3}{3} first if both Cases \hyperlink{case1}{1} and \hyperlink{case2}{2} occur, \eqref{item:merging_disjoint} holds true. 
  Since we have replaced $\B_1,\dotsc,\B_{\ell}$ by the merged ball $\B$ whenever $\B$ satisfies \eqref{eq:equality_case}, items \eqref{item:merging_preservation}, \eqref{item:merging_containment} and \eqref{item:merging_radius} still hold.
  For \eqref{item:merging_energy}, we use in order that $\bigcup_{i = 1}^\ell \B_i \subset \B$, \eqref{item:merging_energy} and \eqref{item:merging_radius} holds for each $\B_i\in \mathcal{B}^-_{\mathrm{F}}(t_j-)$, Lemma \ref{lem:Lambda_mean}, that $r(\B) = \sum_{i=1}^{\ell}r(\B_i)$ and \eqref{eq:equality_case} to bound 
  \begin{align*}
    \int_{\B} f(\lvert \D u\rvert) \d x 
    &\geq  \sum_{i = 1}^\ell \int_{\B_i} f(\lvert \D u\rvert) \d x \\
		&\geq  \sum_{i = 1}^\ell\mathrm{E}_{\sg}^{1,2}(\Tr_{\partial \B_i}u)\,\Lambda_f\left(\frac{r(\B_i)}{\mathrm{E}_{\sg}^{1,2}(\Tr_{\partial \B_i}u)}\right) \\
		&\geq \sum_{i = 1}^\ell\frac{r(\B_i)}{t_j} \Lambda_f(t_j) \\
		&= \frac{r(\B)}{t_j} \Lambda_f(t_j) \\
		&=\mathrm{E}_{\sg}^{1,2}(\Tr_{\partial \B}u)\,\Lambda_f\left(\frac{r(\B)}{\mathrm{E}_{\sg}^{1,2}(\Tr_{\partial \B}u)}\right).
	\end{align*}

  \subsubsection*{Step 5}\hypertarget{step5}{} Termination.

  Observe the sum of radii $r(\mathcal B^-(t))$ is continuous and non-decreasing in $t$, and by \eqref{item:merging_radius} satisfies
  \begin{equation}
    r(\mathcal B^-(t)) \geq \sum_{\B \in \mathcal B(t)}\mathrm E^{1,2}_{\sg}(\Tr_{\partial \B}u) t \geq \mathrm E^{1,2}_{\sg}(\Tr_{\partial\Omega}u) t,
  \end{equation} 
  where we used the subadditivity inequality \eqref{eq:esg_subadditive} of the singular energy.
  Since the right hand side is unbounded as $t \to \infty$, there exists $t= T>0$ where
  \begin{equation}
    r(\mathcal B(T)) = \eta,
  \end{equation} 
  and hence we must eventually reach Case \hyperlink{case3}{3} in the iteration of Step \hyperlink{step2}{2}.
  When this occurs, by construction the collection $\mathcal B \doteq \mathcal B^-(T)$ is a collection of disjoint balls $\{\B(c_i,\sigma_i)\}_{i=1}^N$ for which
  \begin{equation}
    \int_{\B(c_i,\sigma_i)} f(\lvert \D u\rvert) \d x  \geq \mathrm E^{1,2}_{\sg}( c_i,\sigma_i)\Lambda_f \left(\frac{\sigma_i}{\mathrm E^{1,2}_{\sg}( c_i,\sigma_i)}\right),
  \end{equation} 
  by \eqref{item:merging_energy} and
  \begin{equation}\label{eq:final_T_inequality}
    \sigma_i \leq \mathrm{E}^{1,2}_{\sg}( c_i,\sigma_i) T 
  \end{equation} 
  by \eqref{item:merging_radius}.
  Moreover we have
  \begin{equation}
    \eta = r(\mathcal B^-(T)) = r(\mathcal B^+(T)) \geq \sum_{\B \in \mathcal B^+(T)} \mathrm E^{1,2}_{\sg}(\Tr_{\partial \B}u)  T \geq \mathrm{E}^{1,2}_{\sg}(\Tr_{\partial\Omega}u) T,
  \end{equation} 
  so applying Lemma \ref{lem:merging_energy_bound} we deduce that
  \begin{equation}
    \int_{\bigcup \mathcal B} f(\lvert \D u\rvert) \d x \geq \left( \sum_{\B \in \mathcal B} \mathrm{E}^{1,2}_{\sg}(\Tr_{\partial \B}u) \right) \Lambda_f\left( \frac{\eta}{\sum_{\B \in \mathcal B} \mathrm{E}^{1,2}_{\sg}(\Tr_{\partial \B}u)} \right) ,
  \end{equation} 
  establishing our desired lower bound for the collection $\mathcal B$.
\end{Proof}

\begin{Rmk}\label{rem:merging_moreover}
  The proof of Proposition \eqref{prop:genmergingballlemma} moreover shows that, analogously to what is proven in \cite[Prop.~4.3]{vanschaftingen2023asymptotic}, for each $\eta \in (0,\eta_0]$ there exists a collection $\mathcal S$ of circles $\mathbb S^1(a,\rho) = \partial \B(a,\rho)$ (understanding $\mathbb S^1(c,0)=\{c\}$) for which the following holds.
  \begin{enumerate}[(a)]
    \item We have $\bigcup \mathcal S$ is a disjoint union of finitely many annuli; that is, there exists $N \in \mathbb N$ and centres $c_j \in \Omega$, radii $0 \leq \underline{r}_j \leq \overline{r}_j \leq \eta$ for each $1 \leq j \leq N$ such that $\mathbb S^1(c_j,r) \in \mathcal S$ for each $\underline r_j \leq r < \overline r_j$.

    \item for each $\eps \in (0,\eta]$ there exists finitely many $\{\mathbb S^1(c_{j_i},\rho_i)\}_{i=1}^{\ell} \subset \mathcal S$ such that $\sum_{i=1}^{\ell} \rho_i = \eps$ and $\{a_i\}_{i=1}^{\kappa} \subset \bigcup_{i=1}^{\ell} \B(c_{j_i},\rho_i)$.
      If $\eps = \eta,$ this collection corresponds precisely to $\{ \partial \B : \B \in \mathcal B(\eta)\} \subset \mathcal S$.

    \item For each $1 \leq j \leq N$ we have the localised estimate
    \[
      \Esgnu(\Tr_{\partial \B(c_j,r)}u) \Lambda_f\Big(\frac{r}{\Esgnu(\Tr_{\partial \B(c_j,r)}u) } \Big) 	\leq \int_{\B(c_j,r)}f(|\D u|)\,\d x \quad \text{for each } r \in [\underline r_i,\overline r_i).
   \]
  \end{enumerate}
  This can be seen by defining
  \begin{equation}
    \mathcal S = \{ \partial \B : \B \in \mathcal{B}_{\mathrm{G}}^+(t),\  0 \leq t \leq T\},
  \end{equation} 
  and verifying that the claimed properties follow from the items \eqref{item:merging_disjoint}--\eqref{item:merging_terminate} in the proof of Proposition \ref{prop:genmergingballlemma}.
  We will omit the details, as we will not need this in our subsequent analysis.
\end{Rmk}

\subsection{Asymptotic lower bounds}

We will show how this merging ball construction can be used to obtain bounds for (almost-)minimising sequences in the class $\mathrm{R}^{1,2}_{\delta}(\Omega;\mathcal N)$.
For this we will first establish a general lower bound valid for all sequences, namely \eqref{eq:firstorderbndndess}, and then we will prove refined bounds along sequences attaining this bound in the sense of \eqref{eq:secondorder}.

\begin{proposition}\label{prop:firstorderconve} Let $g \in \WW^{\sfrac{1}{2},2}(\partial \Omega;\mathcal N)$ and $\delta>0$, and suppose $(u_n)_n \subset \RR^{1,2}_{\delta}(\Omega;\mathcal N)$ is a sequence such that $\Tr_{\partial\Omega}u_n  = g$ for all $n$.
  Then given a sequence $f_n \colon \mathbb R_+ \to \mathbb R$ of approximating integrands, the following holds:
	\begin{enumerate}[\rm(i)]
    \item\label{item:firstorder} We always have
      \begin{equation}\label{eq:firstorderbndndess}
        \Esgnu([g]) \leq \varlimsup_{n \to \infty} \frac{1}{\mathcal V(f_n)} \int_{\Omega} f_n(\lvert\D u_n\rvert) \d x,
      \end{equation} 
      where the right-hand side may in general be infinite.

			\item \label{item:2ndorder}   Assume that
		\begin{equation}\label{eq:secondorder}
      \varlimsup_{n\to\infty}\left[\int_{\Omega}f_n(|\D u_n|) - \mathcal V(f_n)\Esgnu(g)\right] \text{ is finite}.
		\end{equation}
		Then for each $\eta \in (0, \eta_0)$ and $n \in \mathbb N$, there exists a finite collection of balls $\mathcal B_n^{\mathrm{Top}}(\eta)$ whose sum of radii is at most $\eta$, 
    $(\Tr_{\partial \B}u_n)_{\B \in \mathcal B_n^{\mathrm{Top}}(\eta)}$ forms a topological resolution of $g$, we have $\Esgnu(\Tr_{\partial \B}u)>0$ for all $\B \in \mathcal B_n^{\mathrm{Top}}(\eta)$, and the following asymptotic lower bound
    \begin{equation}\label{eq:limsupesg}
			\varlimsup_{n\to \infty}\sum_{\B \in \mathcal B_{n}^{\mathrm{Top}}(\eta)} \Esgnu(\Tr_{\partial \B}u_n) \leq \varlimsup_{n\to\infty}\mathcal V(f_n)^{-1}\int_{\bigcup\mathcal B_{n}^{\mathrm{Top}}(\eta)} f_n(|\D u_n|)
		\end{equation} 
    holds.
    Moreover we have the remainder estimate
		\begin{equation}\label{eq:bndawayfrmSing}
      \begin{split}
      &\varlimsup_{n\to\infty}\int_{\Omega\setminus \bigcup \mathcal B_n^{\mathrm{Top}}(\eta)}f_n(|\D u_n|) \\
      &\quad\leq \varlimsup_{n\to\infty}\left[\int_{\Omega}f_n(|\D u_n|) - \mathcal V(f_n)\Esgnu(g)\right]  +\Esgnu(g)\log \frac{\sys\Esgnu(g)}{2\eta}.
    \end{split}
		\end{equation}
	\end{enumerate}
\end{proposition}

Under the assumptions of \eqref{item:2ndorder}, note that \eqref{eq:secondorder} implies that 
\begin{equation}\label{eq:singularenergyaslim}
	\Esgnu([g]) = \lim_{n \to \infty}\frac{1}{\mathcal V(f_n)}\int_\Omega f_n(|\D u_n|),
\end{equation}
so \eqref{item:firstorder} holds with equality.
We also point out that $\mathcal B^{\mathrm{Top}}_n(\eta)$ is not the collection of balls $\mathcal B_n(\eta)$ obtained by applying Proposition \ref{prop:genmergingballlemma} to $u_n$ and $f_n$, as we will discard any balls carrying zero singular energy (\emph{i.e.}\,balls $\B$ for which $\Esgnu(\Tr_{\partial \B}u_n)=0$).

\begin{lemma}\label{lemma:asympLambda}
	For a family of approximating integrands $(f_n)_n$, we have
	\begin{equation}\label{eq:asympLambda}
    \lim_{n\to\infty}\left[\Lambda_{f_n}(t) - \mathcal V(f_n)\right] = \log \frac{\sys t}{2} \quad\mbox{for all $t>0$}.
	\end{equation}
  where $\Lambda_{f_n}(t)$ was defined in \eqref{eq:lambda_general}.
	In particular,
	\begin{equation}\label{eq:asympLambdaFirstorder}
    \lim_{n\to\infty}\frac{\Lambda_{f_n}(t)}{\mathcal V(f_n)} = 1.
	\end{equation}
\end{lemma}

This result follows by applying Lemma \ref{lemma:entropy_error}, with $\rho = t$ and $\lambda = \frac{4\pi}{\sys}$.

\begin{Proof}{Proposition}{prop:firstorderconve}
	  For \eqref{item:firstorder}, we set
  \begin{equation}\label{eq:firstorderbndndessbar}
    \mathrm{L}_1 \doteq \varlimsup_{n \to \infty} \mathcal{V}(f_n)^{-1} \int_{\Omega} f_n(\lvert \D u_n\rvert) \d x,
  \end{equation} 
  which we assume is finite, since \eqref{eq:firstorderbndndess} is vacuously true otherwise.
  Let $\eta_0 = \min\big\{\frac{\sys}{2\pi t_0},\delta\big\}$, then by Proposition \ref{prop:genmergingballlemma} applied to $u_{n}$ and $f_n$, 
	there exists a collection $\mathcal{B}_n(\eta)$ of disjoint closed balls in $\Omega$ whose sum of radii is $\eta$, such that $\bigcup\mathcal B_n(\eta)$ contains singular points of $u_n$ and that
	\begin{equation}\label{eq:lowerboundinprop2}
    \begin{split}
    \mathrm{E}(u_n,\mathcal B_n(\eta))\, \Lambda_{f_n}\Big(\frac{\eta}{\mathrm{E}(u_n,\mathcal B_n(\eta))} \Big)
    \leq \int_{\bigcup \mathcal B_{n}(\eta)} f_n(\lvert \D u_{n}\rvert ) \,\d x ,
    \end{split}
	\end{equation} 
  where we abbreviate
  \begin{equation}
    \mathrm{E}(u_n,\mathcal B_n(\eta)) \doteq \sum_{\B \in \mathcal B_n(\eta)}\Esgnu(\Tr_{\partial \B} u_{n}).
  \end{equation} 
  Let $0 \leq \mathrm{M} < \limsup_{n \to \infty}\mathrm{E}(u_n,\mathcal B_n(\eta))$, so that we have $\mathrm{M} \leq \mathrm{E}(u_n,\mathcal B_n(\eta))$ for infinitely many $n$.
  Thus for any such $n$, combining \eqref{eq:lowerboundinprop2} with the mean value property of $\Lambda_{f_n}$ (Lemma \ref{lem:Lambda_mean}) we have
	\begin{equation}\label{eq:importantlimsup}
		\mathrm{M}\,\Lambda_{f_n}({\eta}/{\mathrm{M}}) \leq \int_{\bigcup\mathcal B_{n}(\eta)}f_n(|\D u_{n}|).
	\end{equation}
  Multiplying both sides by $\mathcal V(f_n)^{-1}$ we obtain, by \eqref{eq:asympLambdaFirstorder} from Lemma \ref{lemma:asympLambda},
	\begin{align*}
		\mathrm{M} &= \lim_{n \to \infty}\mathcal V(f_n)^{-1} \mathrm{M}\,\Lambda_{f_n}({\eta}/{\mathrm{M}}) 
       \leq \limsup_{n \to \infty}\mathcal V(f_n)^{-1}\int_{\bigcup\mathcal B_{n}(\eta)}f_n(|\D u_{n}|) \leq \mathrm{L}_1.
	\end{align*}
  Now increasing $\mathrm{M}$ to equal $\limsup_{n \to \infty} \mathrm E(u_n,\mathcal B_n(\eta))$ we infer that
  \begin{equation}
    \limsup_{n \to \infty} \sum_{\B \in \mathcal B_n(\eta)}\Esgnu(\Tr_{\partial \B} u_{n}) \leq \mathrm{L}_1,
  \end{equation} 
  from which we deduce \eqref{eq:limsupesg}, and combining this with the subadditivity property \eqref{eq:esg_subadditive} of the singular energy, \eqref{eq:firstorderbndndess} also follows.
  We now let
  \begin{equation}
    \begin{split}
      \mathcal B_{n}^{\mathrm{Top}}(\eta) 
      &\doteq \{\B \in \mathcal B_n(\eta) : \Esgnu(\Trace_{\partial \B}\tilde u_n) > 0\}.
    \end{split}
  \end{equation}
  Then \eqref{eq:limsupesg} and \eqref{eq:lowerboundinprop2} remains satisfied with $\mathcal B^{\mathrm{Top}}_n(\eta)$ in place of $\mathcal B_n(\eta)$, since we have a smaller collection of balls.
  Also since for any $\B \in \mathcal B_n(\eta) \setminus \mathcal B_n^{\mathrm{Top}}(\eta)$ we have $\Esgnu(\Tr_{\partial \B}u_n)=0$, it follows that $\Tr_{\partial \B}u_n$ is null-homotopic in $\mathcal N$, and hence we can find $\tilde u_{n,\B} \in \WW^{1,2}(\B,\mathcal N)$ satisfying $\Tr_{\partial \B}u_n   = \Tr_{\partial \B}\tilde u_{n,\B}$.
  Then we can define $\tilde u_n \in \WW^{1,2}(\Omega \setminus \bigcup\mathcal B_n^{\mathrm{Top}}(\eta);\mathcal N)$ by setting $\tilde u_n = \tilde u_{u,\B}$ in each $\B \in \mathcal B_n(\eta) \setminus \B_n^{\mathrm{Top}}(\eta)$, and setting $\tilde u_n = u_n$ otherwise.
  Using this construction, it follows that $(\Tr_{\partial \B}u_n)_{\B \in \mathcal B_n^{\mathrm{Top}}(\eta)}$ remains a topological resolution of $g$.

  For \eqref{item:2ndorder} we assume that
  \begin{equation}
    \mathrm{L}_2 \doteq \varlimsup_{n\to\infty}\left[\int_{\Omega}f_n(|\D u_n|)\d x - \mathcal V(f_n)\Esgnu(g)\right]
  \end{equation} 
  is finite,
  and dividing both sides by $\mathcal V(f_n)$ we see that
  \begin{equation}
    \lim_{n \to \infty} \frac1{\mathcal V(f_n)}\left[\int_{\Omega} f_n(\lvert \D u_n\rvert) \d x - \Esgnu(g)\right] = \lim_{n \to \infty} \frac{\mathrm{L}_2}{\mathcal V(f_n)}  = 0.
  \end{equation} 
  In particular we can apply \eqref{item:firstorder}, so for $\eta \in (0,\eta_0)$ let $\mathcal B_n^{\mathrm{Top}}(\eta)$ be the collection of balls as constructed above, noting that \eqref{eq:lowerboundinprop2} remains satisfied.
  Then by subadditivity of the singular energy and Lemma \ref{lem:Lambda_mean},
  \begin{equation}\label{eq:esg_lowerbound}
    \Esgnu(g) \Lambda_{f_n}\left( \frac{\eta}{\Esgnu(g)} \right) \leq \mathrm{E}(u_n,\mathcal B_n^{\mathrm{Top}}(\eta)) \Lambda_{f_n}\left( \frac{\eta}{\mathrm{E}(u_n,\mathcal B_n^{\mathrm{Top}}(\eta))} \right) \leq \int_{\bigcup \mathcal B_n^{\mathrm{Top}}(\eta)} f_n(\lvert \D u_n\rvert) \d x
  \end{equation} 
  holds for all $n$.
  Since by Lemma~\ref{lemma:asympLambda} we have
  \begin{equation}
    \lim_{n \to \infty} \left[ \Lambda_{f_n}\left( \frac{\eta}{\Esgnu(g)} \right) - \mathcal V(f_n) \right]  = \log \Big( \frac{\sys \Esgnu(g)}{2\eta} \Big) ,
  \end{equation} 
  it follows that
  \begin{equation*}
    \begin{split}
      \varlimsup_{n \to \infty} \int_{\Omega \setminus \bigcup \mathcal B_n^{\mathrm{Top}}(\eta)} f_n(\lvert \D u_n\rvert) \d x 
      &\leq \varlimsup_{n \to \infty} \left[ \int_{\Omega} f_n(\lvert \D u_n\rvert) \d x - \mathcal V(f_n) \right] \\
      &\quad- \varliminf_{n \to \infty} \Big[ \int_{\bigcup \mathcal B_n^{\mathrm{Top}}(\eta)} f_n(\lvert \D u_n\rvert) - \Esgnu(g) \Lambda_{f_n}\left( {\eta}/{\Esgnu(g)} \right)  \Big] \\
      &\quad+ \Esgnu(g) \lim_{n \to \infty} \left[ \Lambda_{f_n}\left. ( {\eta}/{\Esgnu(g)}) \right. - \mathcal V(f_n) \right] \\
      &\leq \mathrm{L}_2 + \Esgnu(g) \log \Big( \frac{\sys \Esgnu(g)}{2\eta} \Big),
    \end{split}
  \end{equation*} 
  where we used \eqref{eq:esg_lowerbound} to bound the second term by zero, thereby establishing the result.
\end{Proof}

\section{Compactness results}\label{sec:compactness-results}

\subsection{Compactness theorem}
Using the asymptotic lower bounds established in Proposition \ref{prop:firstorderconve}, we can infer compactness for sequences $(u_n)_n$ satisfying the renormalised bound \eqref{eq:asymptotic_finiteness}.
More precisely, we will show the following.

\begin{theorem}\label{thm:convsub}
  Let $g \in \WW^{\sfrac12,2}(\partial\Omega;\mathcal N)$ and $(f_n)_n$ be a sequence of approximating integrands.
  Given a sequence $(u_n)_n \subset \WW^{1,1}(\Omega;\mathcal N)$  satisfying $\Tr_{\partial\Omega}u_n = g$ for all $n \in \mathbb N$, assume that
    \begin{equation}\label{eq:asymptotic_finiteness}
      \mathrm L \doteq \varliminf_{n \to \infty} \left[ \int_{\Omega} f_n(\lvert \D u_n\rvert) \d x- \mathcal V(f_n) \Esgnu(g) \right]  \quad\text{is finite}.
    \end{equation} 
    Then there exists an unrelabelled subsequence and a limit map $u_{\ast} \in \WW^{1,1}(\Omega;\mathcal N)$ satisfying $\tr_{\partial\Omega}u_{\ast}=g$ and for which the following holds:
    \begin{enumerate}[\rm(i)]
      \item\label{item:convsubseq} The limit map $u_{\ast} \in \WW^{1,2}_{\ren}(\Omega;\mathcal N)$ with $\kappa \in \mathbb N$ distinct singular points $a_1,\dotsc,a_{\kappa} \in \Omega$, such that for each $0<\rho<\rho_{\Omega}(a_1,\dotsc,a_{\kappa})$, $(\tr_{\mathbb S^1}u_{\ast}(a_i+\rho\,\cdot))_{i=1}^{\kappa}$ is a minimal topological resolution of $g$, in that
        \begin{equation*}
          \Esgnu(g) = \sum_{i=1}^{\kappa} \frac{\lambda([u_{\ast},a_i])^2}{4\pi}.
        \end{equation*} 
      In particular, $\kappa \leq {4\pi \Esgnu(g)}/{\sys^2}$.
      \item\label{item:limit_renormalisable} The renormalised energy can be controlled:
        \(
          \mathrm E^{1,2}_{\ren}(u_{\ast}) + \mathrm H([u_{\ast},a_i])_{i=1}^{\kappa} \leq \mathrm L.
        \) 
\item\label{item:isboundaway} For $0 < \rho < \rho_{\Omega}(a_1,\dotsc,a_{\kappa})$,
        \begin{align*}
          \varliminf_{n \to \infty} \int_{\Omega \setminus \bigcup_{i=1}^{\kappa}\B(a_i,\rho)} f_n(\lvert \D u_n\rvert) \d x &\geq \int_{\Omega \setminus \bigcup_{i=1}^{\kappa} \B(a_i,\rho)} \frac12\lvert \D u\rvert^2 \d x,  \\
          \varlimsup_{n \to \infty} \int_{\Omega \setminus \bigcup_{i=1}^{\kappa}\B(a_i,\rho)} f_n(\lvert \D u_n\rvert) \d x
          &\leq\int_{\Omega \setminus \bigcup_{i=1}^{\kappa} \B(a_i,\rho)} \frac12\lvert \D u\rvert^2 \d x
          + \left(\mathrm{L} - \mathrm E_{\ren}^{1,2}(u_{\ast}) - \mathrm H([u_{\ast},a_i])_{i=1}^{\kappa}\right).
          \end{align*}
        \item\label{item:isboundloc} For $0 < \rho< \rho_{\Omega}(a_1,\dotsc,a_{\kappa})$
          and all $i = 1,\dotsc,\kappa$,
          \begin{align*}
          \varliminf_{n \to \infty} \left[\int_{\B(a_i,\rho)} f_n(\lvert \D u_n\rvert) \d x - \frac{\lambda([u_{\ast},a_i])^2}{4\pi} \mathcal V(f_n)\right]
          &\geq \mathrm{E}^{1,2}_{\ren}(u_{\ast};\B(a_i,\rho)) \\
          &+ \frac{\lambda([u_{\ast},a_i])^2}{4\pi} \log\frac{2\pi}{\lambda([u_{\ast},a_i])^2} \\
          \varlimsup_{n \to \infty} \left[\int_{\B(a_i,\rho)} f_n(\lvert \D u_n\rvert) \d x - \frac{\lambda([u_{\ast},a_i])^2}{4\pi} \mathcal V(f_n)\right]
          &\leq \mathrm{L} - \mathrm{E}_{\ren}^{1,2}(u_{\ast})- \mathrm{H}([u_{\ast},a_i])_{i=1}^{\kappa} \\
          &\hspace{-10mm}+ \mathrm{E}^{1,2}_{\ren}(u_{\ast};\B(a_i,\rho)) + \frac{\lambda([u_{\ast},a_i])^2}{4\pi} \log \frac{2\pi}{\lambda([u_{\ast},a_i])}.
        \end{align*}      
        recalling the localised renormalised energy on balls was defined in \eqref{eq:localised_renormalised}.

      \item\label{item:vaguemeasures} We have narrow  convergence of measures
        \begin{equation*}
          \frac{f_n(\lvert \D u_n\rvert)}{\mathcal V(f_n)}  \Leb^2 \mres \Omega \rightharpoonup \sum_{i=1}^{\kappa} \frac{\lambda([u_{\ast},a_i])^2}{4\pi} \delta_{a_i},
        \end{equation*} 
    in that the measures converge \weaklystar and the total mass converges.
      \item\label{item:L1bound} For $0 < \rho < \rho_{\Omega}(a_1,\dotsc,a_{\kappa})$ and each $i = 1,\dotsc,\kappa$ we have
      \begin{equation*}
        \begin{split}
          \varliminf_{n \to \infty} \int_{\B(a_i,\rho)} \lvert \D u_n\rvert \d x 
          \leq \frac{2\pi \rho}{\lambda([u_{\ast},a_i])} &\bigg[\mathrm{L} - \mathrm{E}_{\ren}^{1,2}(u_{\ast}) - \mathrm{H}([u_{\ast},a_i])_{i=1}^{\kappa} + \lambda([u_{\ast},a_i])^2  \\
          &\quad+ \int_0^{\rho}\int_{\partial\B(a_i,r)} \frac12\lvert \D u\rvert^2\d\mathcal H^1 - \frac{\lambda([u_{\ast},a_i])^2}{4\pi r} \d r \bigg].
          \end{split}
      \end{equation*} 
      \item\label{item:weakconvergence} The sequence $(\D u_n)_n$  is equi-integrable and converges weakly to $\D u_{\ast}$ in $\LL^1(\Omega)$.
    \end{enumerate}
\end{theorem}

Observe that if equality is attained in \eqref{item:limit_renormalisable}, that is, we have
\begin{equation}
  \lim_{n \to \infty}\left[ \int_{\Omega} f_n(\lvert \D u_n\rvert) \d x - \mathcal{V}(f_n) \Esgnu(g)\right] = \Eren^{1,2}(u_{\ast}) + \mathrm{H}([u_{\ast},a_i])_{i=1}^{\kappa},
\end{equation} 
then the upper and lower limits in \eqref{item:isboundaway} and \eqref{item:isboundloc} coincide.
We will analyse this equality case in further detail in Theorem \ref{thm:strongconv}.
Additionally as the proof will reveal, this result also applies in the case $\mathrm{L} = -\infty$; more precisely such a possibility is ruled out by \eqref{item:limit_renormalisable}.

We will present the proof of Theorem \ref{thm:convsub} in Section \ref{sec:proof_constrained}.
For this we will need results concerning weak $\LL^1$-compactness and strong convergence of traces, which we first collect in Sections \ref{sec:weak-convergence-of-bounded-sequences} and \ref{sec:strong_traces} respectively.

\subsection{Weak convergence of bounded sequences}\label{sec:weak-convergence-of-bounded-sequences}
Our first ingredient will be a compactness result away from the singularities.
More precisely, we prove a weak compactness result in $\LL^1$ under the presence of the uniform bound \eqref{eq:boundedseq}.

\begin{proposition}\label{prop:compacitnoprblem}
	Let $\Omega \subset \R^d$ be a bounded open set, and $f_n : \R_+ \to \R_+$ be a sequence of Young functions converging pointwise to $\frac12 t^2$.
  Let $u_n  \in \WW^{1,1}(\Omega;\mathbb R^{\nu})$. Assume that
	\begin{equation}\label{eq:boundedseq}
    \mathrm M \doteq \varliminf_{n\to \infty} \int_\Omega f_n (|\D u_n|) \d x  \quad\text{is finite}.
	\end{equation}
  Then up to an unrelabelled subsequence, the following holds.
  \begin{enumerate}[\rm(a)]
    \item\label{item:unconstrained_bound} We have
	\begin{equation}\label{eq:asympequicont}
    \varlimsup_{n \to \infty}\int_{A}|\D u_n| \d x \leq \sqrt{2\mathrm{M} \Leb^d(A)}\quad\mbox{for all $A \subset \Omega$ measurable},
	\end{equation}
	and $(|\D u_n|)_n$ is uniformly integrable.

\item\label{item:unconstrained_conv} If we additionally have $u_n$ converges strongly to a limit map $u$ in $\LL^1(\Omega;\mathbb R^{\nu})$, then $u \in \WW^{1,2}(\Omega;\mathbb R^{\nu})$, $\D u_n$ converges weakly to $\D u$ in $\LL^1(\Omega)$, and 
	\begin{equation}\label{eq:fatoucompactplb}
	\int_\Omega \frac{1}{2}|\D u|^2 \d x \leq \varliminf_{n\to \infty} \int_{\Omega}f_n (|\D u_n|) \d x. 
	\end{equation}
	Furthermore, if $u_n(x) \in \mathcal N$ for all $n$ and almost every $x \in \Omega$, then the limit $u$  satisfies $u(x) \in \mathcal N$ for almost every $x \in \Omega$.
  \end{enumerate}
\end{proposition}

Proposition \ref{prop:compacitnoprblem} is a general compactness result for Young functions, we do not assume \ref{hyp:cdec2} nor \ref{hyp:int}.
We will primarily use Proposition~\ref{prop:compacitnoprblem} in the case $d=2$, however we have stated it in this generality as we will also apply it in $\LL^1(\mathbb S^1,\mathcal H^1)$, which we can identify with $\LL^1((0,1))$.
We will use the latter case in the proof of Theorem \ref{thm:convsub}, and we will omit the details of the necessary modifications.

\begin{Rmk}[Weak and \weakstar convergence]\label{rmk:weakstar}
We recall the difference between \weakstar and weak convergence for our uses. We say a sequence of functions $(U_n)_n \subset \LL^1(\Omega)$ converges weakly in $\LL^1$ to a limit map $U \in \LL^1(\Omega)$ if 
\begin{equation}\label{eq:weakconvegrnece}
  \lim_{n \to \infty} \int_\Omega U_n \phi \d x = \int_\Omega U \phi \d x \quad\mbox{for all $\phi \in \LL^\infty(\Omega)$}.
\end{equation}
We say a sequence $(U_n)_n \subset \LL^1(\Omega)$ converges \weaklystar as measures to $U\in \LL^1(\Omega)$ if \eqref{eq:weakconvegrnece} only holds for all $\phi \in \CC_0(\Omega)$. 
Note that for mappings $(u_n)_n$ in $\WW^{1,1}$, weak${}^{\ast}$ convergence of $(\D u_n)_n$ corresponds to \weakstar convergence of $(u_n)$ in $\BV$.
\end{Rmk}

The following variant of the Dunford-Pettis theorem \cite[Thm.~2.4.5]{attouch2014variational} will be of use.
\begin{lemma}\label{lemma:ourinstanceofDunfordPettis}
  Let $\Omega \subset \mathbb R^d$ be a bounded open set, and consider a sequence of mappings $(u_n)_n \subset \WW^{1,1}(\Omega; \mathcal N)$ converging strongly in $\LL^1(\Omega)$.
  Then $(\D u_n)_n$ converges weakly in $\LL^1$ if and only if $\{|\D u_n| : n \in \mathbb N\}$ is uniformly integrable.
\end{lemma}

Classically, the Dunford-Pettis theorem asserts that weak compactness of a sequence $(U_n)_n$ in $\LL^1(\Omega)$ is equivalent to uniform integrability of this family. 
Recall that equi-integrability and uniform integrability are synonymous.
In our setting we apply this with $U_n = \D u_n$, noting the limiting map is necessarily unique by strong convergence of $(u_n)_n$ in $\LL^1(\Omega)$.

In the proof of Proposition \ref{prop:compacitnoprblem}, we will use the following result (Lemma \ref{lemma:asympunifint}), which allows us to convert ``asymptotic'' uniform integrability into the well--known statement of uniform integrability.
While we will state the result for the Lebesgue measure in $\mathbb R^d$, we note the following also holds in more general measure spaces.

\begin{lemma} \label{lemma:asympunifint} For a sequence of functions $(U_n)_n \subset \LL^1(\Omega)$, the following are equivalent:
	\begin{enumerate}[\rm(i)]
		\item \label{item:supunifint} the sequence $(U_n)$ is uniformly integrable; that is for all $\eps>0$ there exists $\delta >0$ such that for every measurable  $A \subset \Omega$ such that $\Leb^d(A) \leq \delta$, we have $\ds \sup_{n}\int_{A}|U_n|\d x \leq \epsilon$,
		\item \label{item:limsupunifint} for all $\eps>0$ there exists $\delta >0$ such that for every measurable $A \subset \Omega$ such that $\Leb^d(A) \leq \delta$, we have  $\ds \varlimsup_{n \to \infty}\int_{A}|U_n|\d x \leq \epsilon$.
	\end{enumerate}
\end{lemma}

This result follows by applying the Vitali-Hahn-Saks Theorem (see \emph{e.g.}\,\cite[Thm.~2.53]{FonsecaLeoni2007}), considering a countable family of measurable sets on which the superior limit in \eqref{item:limsupunifint} is a limit.

The Vitali-Hahn-Saks is usually proven using Baire's category theorem, and it  also proves Lemma \ref{lemma:asympunifint} directly along similar lines as follows: 
 given $\eps>0$, let $\delta>0$ as in \eqref{item:limsupunifint} and consider the complete metric space $\mathcal A_{\delta}$ of measurable subsets $A \subset \Omega$ satisfying $\Leb^d(A) \leq \delta$, with respect to the metric $\Leb^d(A \triangle B) = \int_{\Omega} \lvert \chi_A - \chi_B \rvert \d x$.
Then consider the closed subsets
\begin{equation}
  \mathcal M(m) = \bigcap_{n \geq m}\left\{A \in \mathcal A_\delta : \int_A |U_n| \d x\leq 2\epsilon\right\},
\end{equation} 
which by \eqref{item:limsupunifint} satisfies $\bigcup_m \mathcal M(m) = \mathcal A_{\delta}$.
Therefore by Baire's theorem, there is some $\mathcal M(m_0)$ with non-empty interior, and hence contains an $r$-ball for some $r>0$.
One can then verify that setting $\delta_0 = \min\{\delta,r\}$, we have $\Leb^d(A) \leq \delta_0$ implies that $\sup_{n \geq m_0} \int_A U_n \d x < 2\eps$.
Since the finite collection $\{U_1,\dotsc,U_{m_0-1}\}$ is uniformly integrable, shrinking $\delta_0$ further if necessary \eqref{item:supunifint} follows.

\begin{Proof}{Proposition}{prop:compacitnoprblem} 
  By passing to a subsequence, we can assume that 
  \begin{equation}
    \mathrm{M} = \varlimsup_{n \to \infty} \int_{\Omega} f_n(\lvert \D u_n\rvert) \d x < \infty.
  \end{equation} 
  We will first show \eqref{item:unconstrained_bound}, namely that \eqref{eq:asympequicont} holds. 
  For this let $A \subset \Omega$ be measurable, and assume that $\Leb^d(A) > 0$ since there is otherwise nothing to show.
  Then Jensen's inequality gives
	\begin{equation*}
    \mathrm{M} \geq \Leb^d(A)\varlimsup_{n \to \infty} f_n\Big(\frac{1}{\Leb^d(A)}\int_A |\D u_n| \d x\Big).
	\end{equation*}
  Let
  \begin{equation}
    0 < \mathrm{L} < \frac{1}{\Leb^d(A)}\varlimsup_{n \to \infty} \int_{A} \lvert \D u_{n}\rvert \d x,
  \end{equation} 
  assuming this limit inferior is positive, as otherwise \eqref{eq:asympequicont} holds since the left-hand side is zero.
  Then passing to a further subsequence where the limit supremum is attained, we can assume that $\mathrm{L}\, \Leb^d(A) \leq \int_A \lvert \D u_{n}\rvert \d x$ for all $n$ sufficiency large.
  For such $n$, since each $f_n(t)$ is non-decreasing,
  \(
  \mathrm{M} \geq \Leb^d(A) \varlimsup_{n \to \infty} f_n(\mathrm{L}) = \Leb^d(A) {\mathrm{L}^2}/2,
  \) 
  and hence that
  \begin{equation}
    \mathrm{M} \geq \frac1{2\Leb^d(A)} \Big(\varlimsup_{n \to \infty} \int_A \lvert \D u_{n}\rvert\d x\Big)^2,
  \end{equation} 
  which rearranges to give \eqref{eq:asympequicont}.
  Therefore we have $(\lvert\D u_n\rvert)_n$ is uniformly integrable by Lemma \ref{lemma:asympunifint}, thereby establishing \eqref{item:unconstrained_bound}.

  Now for \eqref{item:unconstrained_conv}, by the Dunford-Pettis theorem (see Lemma \ref{lemma:ourinstanceofDunfordPettis}) we have the strong $\LL^1$ convergence and uniform integrability of the gradients implies that $\D u_n$ converges weakly to $\D u$ in $\LL^1(\Omega)$.
  To show \eqref{eq:fatoucompactplb},
  let $(\kappa_m)_m \subset (0,\infty)$ be any sequence such that $\kappa_m \nearrow \infty$.
  Then by Lemma \ref{lem:fn_conv}\eqref{item:derivative_pointwise} we have $\lim_{n \to \infty} f_{n,+}'(\kappa_m) \to \kappa_m$ for each $m \in \mathbb N$,
  and writing $f(t) = \frac12 t^2$ we consider the truncations $T_{\kappa_m}f_n$ and $T_{\kappa_m}f$ from Definition \ref{defn:kappa_truncation}.
  Note that for each $m$, $T_{\kappa_m}f$ is the pointwise limit of $T_{\kappa_m}f_{n}(t)$ as $n \to \infty$,
  and moreover this limit is uniform on $[0,\kappa_m]$ by Lemma \ref{lem:fn_conv}\eqref{item:local_uniform}.
  Also for $t>\kappa_m$ we have 
  \begin{equation}
    \lvert T_{\kappa_m}f_n(t) - T_{\kappa_m}f(t) \rvert \leq \lvert f'_{n,+}(\kappa_m) - \kappa_m\rvert \lvert t\rvert + \left\lvert f_n(\kappa_m) - f'_{n,+}(\kappa_m)\kappa_m + \frac12 \kappa_m^2\right\rvert,
  \end{equation} 
  where both coefficients vanish as $n \to \infty$.
  Hence setting $A_{m,n} = \{ x \in \Omega : \lvert \D u_n\rvert(x) \leq \kappa_m\}$, 
  for any $m$ we have
  \begin{IEEEeqnarray*}{rCl}
  	   \mathrm{M}&=&\varlimsup_{n \to \infty} \int_{\Omega} f_n(\lvert \D u_n\rvert) \d x 
    \\
    &\geq&  \varlimsup_{n \to \infty} \int_{\Omega} T_{\kappa_m}f_n(\lvert \D u_{n}\rvert) \d x \\
    &\geq& \varlimsup_{n \to \infty} \int_{\Omega} T_{\kappa_m}f(\lvert \D u_{n}\rvert) \d x -  \varlimsup_{n \to \infty}\left\lvert \int_{\Omega} T_{\kappa_m}f_n(\lvert \D u_{n}\rvert) - T_{\kappa_m}f(\lvert \D u_{n}\rvert) \d x\right\rvert\\
    &\geq& \varlimsup_{n \to \infty} \int_{\Omega} T_{\kappa_m}f(\lvert \D u_{n}\rvert) \d x - \varlimsup_{n \to \infty} \int_{A_{m,n}} \left\lvert f_{n}(\lvert \D u_{n}\rvert) - \frac12\lvert \D u_{n}\rvert^2 \right\rvert \d x \\
    && - \varlimsup_{n \to \infty}\lvert f_{n,+}'(\kappa_m) - \kappa_m\rvert \int_{\Omega \setminus A_{m,n}} \lvert \D u_{n}\rvert \d x \\
    && - \varlimsup_{n \to \infty} \left\lvert f_{n}(\kappa_m) - f'_{n,+}(\kappa_m)\kappa_m + \frac12 \kappa_m^2\right\rvert\Leb^d(\Omega \setminus A_{m,n}) \\
    &=& \varlimsup_{n \to \infty} \int_{\Omega} T_{\kappa_m}f(\lvert \D u_{n}\rvert) \d x.
  \end{IEEEeqnarray*}  
  Note that in the last line, we used the $\LL^1$-boundedness of $\lvert \D u_n\rvert$ to show that $\int_{\Omega \setminus A_{m,n}} \lvert \D u\rvert \d x$ is uniformly bounded.
  Since $\D u_n \rightharpoonup \D u$ weakly in $\LL^1$, using the lower-semicontinuity from Lemma \ref{lemma:bv_relaxation} we have
  \begin{equation}
    \mathrm{M} \geq \varlimsup_{n \to \infty} \int_{\Omega} T_{\kappa_m}f(\lvert \D u_n\rvert) \d x \geq \int_{\Omega} T_{\kappa_m}f(\lvert \D u \rvert) \,\d x,
  \end{equation} 
  so sending $m \to \infty$ and applying the monotone convergence theorem we have
  \begin{equation}
    \varliminf_{n \to \infty} \int_{\Omega} f_n(\lvert \D u_n\rvert) \d x = \mathrm{M} \geq \lim_{m \to \infty} \int_{\Omega} T_{\kappa_m}f(\lvert \D u\rvert)\d x = \int_{\Omega} \frac12 \lvert \D u\rvert^2 \d x,
  \end{equation} 
	
  Finally since $(u_n)_n$ converges strongly in $\LL^1$, passing to a subsequence we can assume $u_n$ converges almost everywhere to $u$ in $\Omega$.
  In particular, if for all $n$ we have $u_n(x) \in \mathcal N$ for almost every $x \in \Omega$, then the same holds for the limit map $u$.
\end{Proof}

\subsection{Strong convergence of traces}\label{sec:strong_traces}

As the example $(\chi_{\B(0,1 - \sfrac{1}{n})})_n$ illustrates, the trace operator  $\tr_{\partial \B(0,1)} : \BV(\R^2) \to  \LL^1(\mathbb S^1)$ is not \weaklystar continuous. We therefore need refined convergence results of traces, which we establish in this section. 

\begin{proposition}\label{prop:trace_conv}
	Let $\Omega \subset \mathbb R^2$ be a Lipschitz domain, $(f_n)_n$ a sequence of approximating integrands, and $(u_n)_n \subset \WW^{1,1}(\Omega;\R^\nu)$ be a sequence converging \weaklystar in $\BV$ to $u \in \WW^{1,2}(\Omega;\mathbb R^\nu)$. If
	\begin{equation}\label{eq:boundfortrace}
    	\varlimsup_{n\to \infty} \int_\Omega f_n (|\D u_n|)\d x \quad\text{is finite,}
	\end{equation}
  then passing to a subsequence we have $\Tr_{\partial \Omega}u_n \to \Tr_{\partial\Omega}u$ strongly in $\LL^1(\partial \Omega;\mathbb R^\nu)$.  Moreover
  \begin{equation}\label{eq:limsuptrace}
  	\varlimsup_{n\to \infty}\int_{\partial \Omega}f_n(|\tr_{\partial \Omega}u_n|) \d \mathcal H^1 \text{ is finite}
  \end{equation}
  and
  \begin{equation}\label{eq:trace_lsc}
		\int_{\partial\Omega} \frac{\lvert\Tr_{\partial \Omega}u\rvert^2}{2}\d \mathcal H^1 \leq \varliminf_{n\to \infty} \int_{\partial\Omega} f_n (|\Tr_{\partial \Omega}u_n|)\d \mathcal H^1.
  \end{equation}
\end{proposition}

\begin{lemma}\label{lemma:tracebound}
  Let $\Omega \subset \mathbb R^2$ be convex and $f \colon \R_+ \to \R$ be a Young function satisfying \ref{hyp:ndec2}.
  Then there exists $C>0$ such that for all $v \in \WW^{1,1}(\Omega;\R^{\nu})$,
	\[
		\int_{\partial \Omega} f(|v|)\d \Hau^1 \leq C\int_\Omega 1+f(|v|) +f(|\D v|) \d x.
	\]
  Moreover if $(f_n)_n$ is a sequence of approximating integrands, then the constant $C>0$ may be chosen to be independent of $n$.
\end{lemma}
We will present a proof relying on Gagliardo's $\LL^1$-trace theorem.

\begin{Proof}{Lemma}{lemma:tracebound}
  By considering $\lvert v\rvert$ in place of $v$, we may assume that $v \geq 0$.
  Then by the chain rule and Lemma \ref{lemma:dec2andlegendre} below, we have
  \begin{equation}
    \D (f(v)) = f_+'(v) \D v \leq c_1 (1+f(\lvert v\rvert) + f(\lvert \D v\rvert)),
  \end{equation} 
  so $f(\lvert v\rvert) \in \WW^{1,1}(\Omega)$.
  For a sequence of approximating integrands, Lemma \ref{lemma:dec2andlegendre} also asserts that this constant $c_1>0$ can be chosen independently of $n$.
  Hence by Gagliardo's trace theorem (see \emph{e.g.}\,\cite[Thm.~18.18]{Leoni2017}) there exists $c_2>0$, depending only on $\Omega$, such that
  \begin{equation}
    \int_{\partial\Omega} f(\Tr_{\partial\Omega}v) \d \mathcal H^1 \leq c_2 \int_{\Omega} f(v) + \lvert \D(f(v)) \rvert \d x \leq c_2(1+c_1) \int_{\Omega} 1+ f(\lvert v\rvert) + f(\lvert \D u\rvert) \d x,
  \end{equation} 
  as required.
\end{Proof}

\begin{lemma}\label{lemma:dec2andlegendre}
	If $f : \R_+ \to \R$ is a Young function satisfying \ref{hyp:ndec2}, there is $C>0$ such that
  \begin{equation}
    f'_+(s) t \leq C(1+f(s) + f(t)) \quad\mbox{for all $s,t \geq 0$}.
 \end{equation}
 Moreover if $(f_n)$ is a sequence of approximating integrands, the constant $C>0$ may be chosen to be independent of $n$.
\end{lemma}
\begin{Proof}{Lemma}{lemma:dec2andlegendre}
  By \eqref{eq:convex_linearisation} we know that $f_+'(s) \leq f(2s)/s$, and by \eqref{eq:slope_inequality} we know that $f(t)/t$ is non-decreasing on $(0,\infty)$. Hence we have
  \begin{equation}
    f_+'(s)t  \leq \frac{f(2s)t}{s} \leq 
    \begin{cases} f(2s) &\mbox{if } t \leq s, \\ f(2t) &\mbox{if } t > s,
    \end{cases}
  \end{equation} 
  By \ref{hyp:ndec2} we have $f(2t) \leq 4f(t)$ for all $t \geq t_0$, and for $t< t_0$ we have $f(2t) \leq f(2t_0)$, thereby giving $f(2t) \leq C (1 + f(t))$ for all $t \geq 0$, from which the desired inequality follows.
  Observe that for an approximating family $(f_n)_n$, since $\sup_{n\in\mathbb N} f_n(2t_0)$ is bounded, this constant can be chosen uniformly in $n$.
\end{Proof}

\begin{Proof}{Proposition}{prop:trace_conv}
  Define
  \begin{equation}
    \mathrm L = \varlimsup_{n\to \infty} \int_\Omega f_n (|\D u_n|) \d x,
  \end{equation} 
  which by \eqref{eq:boundfortrace} is finite.
  By \cite[Thm.~18.21(iii)]{Leoni2017} there is $C>0$ for which
  \begin{equation}\label{eq:l1_trace_bound}
    \int_{\partial\Omega} \lvert \Tr_{\partial\Omega}u - \Tr_{\partial\Omega}u_n\rvert\d\mathcal{H}^{1} \leq \frac{C}{\eps} \int_{A_{\eps}} \lvert u - u_n\rvert \d x + C \lvert \D u - \D u_n\rvert(A_{\eps}),
  \end{equation} 
  where $A_{\eps} = \{ x \in \Omega : \dist(x,\partial\Omega) <\eps\}$.
  Applying \eqref{eq:asympequicont} from Proposition \ref{prop:compacitnoprblem} and using that $u_n \to u$ strongly in $\LL^1(\Omega;\mathbb R^{\nu})$ (for instance by \cite[Thm.~5.5]{EvansGariepy2015}),
  we can send $n \to \infty$ in \eqref{eq:l1_trace_bound} to obtain
  \begin{equation}
    \varlimsup_{n \to \infty}\int_{\partial\Omega} \lvert \Tr_{\partial\Omega}u - \Tr_{\partial\Omega}u_n\rvert\d\mathcal{H}^{1} \leq 2 C \sqrt{ 2\mathrm L  \Leb^2(A_{\eps}}),
  \end{equation} 
  which vanishes as $\eps \to 0$.
  In particular, passing to a subsequence we have $\Tr_{\partial\Omega}u_n \to \Tr_{\partial \Omega}u$ pointwise $\mathcal H^1$-almost everywhere, from which \eqref{eq:trace_lsc} follows by Fatou's lemma.
  Finally since $\mathcal N$ is compact, $\sup_{n \in \mathbb N} \|u_n\|_{\LL^\infty(\Omega)}$ is finite (see Section \ref{sec:existence_minima}), so by Lemma \ref{lemma:tracebound} we can estimate
  \begin{equation}
    \limsup_{n \to \infty} \int_{\partial\Omega} f_n(\lvert \tr_{\partial\Omega} u_n\rvert) \d \mathcal H^1 \leq C (1 + \sup_{n \in \mathbb N} \|u_n\|_{\LL^\infty(\Omega)}) \Leb^2(\Omega) + C\,\mathrm{L},
  \end{equation} 
  establishing\eqref{eq:limsuptrace}.
\end{Proof}

The following Lemma \ref{lemma:topologyonthesphere} will be of use in the proof of Theorem \ref{thm:convsub},
in order to associate a well-defined homotopy class to each singularity of the limiting map.
We will use the notation $\D_{\tau}$ for the tangential derivative on any circle $\mathbb S^1(a,r) =\partial \B(a,r)$.

\begin{lemma}\label{lemma:topologyonthesphere}
	Let $a \in \mathbb R^2$ and $r>0$.
	Suppose $(u_n)_n \subset \WW^{1,1}(\partial \B(a,r);\mathcal N)$ and $(f_n)_n$ is a sequence of approximating integrands such that
	\begin{equation}\label{eq:sphere_limsup}
		\varlimsup_{n \to \infty} \int_{\partial \B(a,r)} f_{n}(\lvert \D_{\tau} u_n\rvert) \,\d \mathcal H^1 \quad\text{is finite},
	\end{equation} 
	and $\Tr_{\partial\B(a,r)}u_n \to \Tr_{\partial\B(a,r)}u$ strongly in $\LL^1(\partial \B(a,r); \mathcal N)$ to some limit $u \in \WW^{1,1}(\partial \B(a,r);\mathcal N)$.
	Then we have $u_n$ and $u$ are homotopic on $\partial \B(a,r)$ for all $n$ sufficiently large, and in particular
	\(
		\lambda(\Tr_{\partial\B(a,r)}u_n) = 	\lambda(\Tr_{\partial\B(a,r)}u)\) and \( \Esgnu(\Tr_{\partial\B(a,r)}u_n) = \Esgnu(\Tr_{\partial\B(a,r)}u)
	\) for all such $n$.
\end{lemma}

\begin{Proof}{Lemma}{lemma:topologyonthesphere}
 It suffices to show the traces of $u_n$ and $u$ on $\partial\B(a,r)$ are homotopic for sufficiently large $n$, which will imply the equality of the minimal lengths and singular energies.
	Taking the truncations $T_{\kappa}f_n$ of $f_n$, by \eqref{eq:sphere_limsup} there is $M>0$ such that
	\begin{equation}
		\int_{\partial \B(a,r)} T_{\kappa}f_n(\lvert \D_{\tau}u_n\rvert) \d x \leq \int_{\partial \B(a,r)} f_n(\lvert \D_{\tau}u_n \rvert) \d x \leq M
	\end{equation} 
	for all $n \in \mathbb N$ and $\kappa>0$.
  Since the traces $\Tr_{\partial \B(a,r)}u_n$, $\Tr_{\partial \B(a,r)}u$ are assumed to lie in $\WW^{1,1}(\partial \B(a,r);\mathcal N)$, they admit a continuous representative which we denote simply by $u_n$, $u$ respectively.
	By parametrising $\partial \B(a,r)$ by the map $\theta \mapsto a+(r\cos\theta,r\sin\theta)$, we mollify in $\theta$ to define $u_{n,\eps} \doteq u_n \ast \rho_{\eps}$, which satisfies
  \begin{equation}\label{eq:mollification_linfty}
		\lVert u_{n,\eps} - u_{\eps}\rVert_{L^{\infty}(\partial\B(a,r))} \leq \eps^{-1} \lVert u_n - u \rVert_{L^1(\partial\B(a,r))}.
  \end{equation} 
	
	We will first estimate $\lvert u_{n,\eps} - u_n\rvert$ pointwise. For this let $x,y \in \partial \B(a,r)$, and let $S_{x,y} \subset\partial\B(a,r)$ denote an arc connecting $x$ and $y$ with minimal length.
	Then by the fundamental theorem of calculus we have
	\begin{equation}
		\begin{split}
			\lvert u_n(x) - u_n(y)\rvert 
			&\leq \int_{S_{x,y}} \lvert \D_{\tau} u_n\rvert \,\d \mathcal H^1 \\
			&= \int_{S_{x,y}} \chi_{\{\lvert \D_{\tau} u_n\rvert \geq \kappa\}} \lvert \D_{\tau} u_n\rvert \d\mathcal H^1 
			+ \int_{S_{x,y}} \chi_{\{\lvert \D_{\tau} u_n\rvert < \kappa\}} \lvert \D_{\tau} u_n\rvert \d\mathcal H^1\\
			&\leq \frac1{f_{n,+}'(\kappa)} \int_{S_{x,y}} T_{\kappa}f_n(\lvert \D_{\tau} u_n\rvert) \d\mathcal H^1\\
      &\quad+ \mathcal H^1(S_{x,y}) f_n^{-1}\left( \frac1{\mathcal H^1(S_{x,y})} \int_{S_{x,y}}\chi_{\{\lvert \D_{\tau} u_n\rvert<\kappa\}} f_n(\lvert\D_{\tau} u_n\rvert) \d\mathcal H^1 \right) \\ 
			&\leq \frac{M}{f_{n,+}'(\kappa)}  + \mathcal H^1(S_{x,y}) f_n^{-1}\left( \frac{M}{\mathcal H^1(S_{x,y})} \right).
		\end{split}
	\end{equation} 
	Then noting that
	\(
		\lvert x- y \rvert \leq \mathcal H^1(S_{x,y}) = \dist_{\mathbb S^1}(x,y) = 2r\arcsin\left( {\lvert x- y\rvert}/{2r} \right) \leq 2 \lvert x-y\rvert,
	\) 
	we have
  \begin{equation}\label{eq:trace_diff_mollif}
		\lVert u_{n,\eps} - u_{n}\rVert_{\LL^{\infty}(\partial \B(a,r))} \leq \frac{M}{f'_{n,+}(\kappa)} + 2\eps f_n^{-1}\left( \frac{M}{\eps} \right) .
	\end{equation} 
  Since $u_n \to u$ almost everywhere on $\partial \B(a,r)$, passing to the limit in this estimate (noting by Lemma \ref{lem:fn_conv} that $f_{n,+}'(\kappa)\to \kappa$ and $f_n^{-1}(t) \to \sqrt{2t}$ as $n \to \infty$) we have
  \begin{equation}\label{eq:trace_diff_mollif2}
		\lVert u_{\eps} - u \rVert_{L^{\infty}(\partial\B(a,r))} \leq {M}/{\kappa} + \sqrt{8 M\eps}.
  \end{equation} 
	Now consider the tubular neighbourhood projection $\Pi_{\mathcal N} \colon \mathcal U \to \mathcal N$, so there is $\delta>0$ such that $\mathcal N + \B(0,\delta) \subset \mathcal U$. 
	We then define
	\(
	\tilde H_n(t,x) = (1-t)u_n(x) + t u(x),
	\) 
	and we will show that $\tilde H_n$ maps into $\mathcal U$ for $n$ sufficiently large.
	Indeed for $(t,x) \in [0,1] \times \partial\B(a,r)$ we have by \eqref{eq:mollification_linfty}, \eqref{eq:trace_diff_mollif} and \eqref{eq:trace_diff_mollif2},
	\begin{equation}
		\begin{split}
			\dist(\tilde H_n(x,t); \mathcal N) 
			&\leq \lvert \tilde H_n(t,x) - u_n(x)\rvert \\
			&\leq \lvert u_n - u_{n,\eps}\rvert(x) + \lvert u_{n,\eps} - u_{\eps}\rvert(x) + \lvert u_{\eps}-u\rvert(x)\\
			&\leq \frac{M}{f'_{n,+}(\kappa)} + 2\eps f_n^{-1}\left( \frac{M}{\eps} \right) + \eps^{-2} \lVert u_n - u \rVert_{L^1(\partial \B(a,r))} + \frac{M}{\kappa}+ \sqrt{8M\eps}.
		\end{split}
	\end{equation} 
	Now we choose $\eps>0$ sufficiently small so that
	\(
	\sqrt{8M \eps}< {\delta}/9.
	\) 
	We also take $\kappa > {9M}/{\delta}$, then for $n$ sufficiently large we have 
	\begin{equation}
		\frac{M}{f'_{n,+}(\kappa)} + 2\eps f_n^{-1}\left( \frac{M}{\eps} \right) \leq \frac{2M}{\kappa} + 2 \sqrt{8M \eps} < \frac{4\delta}9,
	\end{equation} 
	and since $u_n \to u$ in $L^1$ on $\partial\B(a,r)$, increasing $n$ if necessary we have a
	\(
		\eps^{-1} \lVert u_n - u \rVert_{L^1(\partial \B(a,r))} < {\delta}/3.
	\) 
	Combining our estimates, we have
	\begin{equation}
		\dist(\tilde H_n(t,x); \mathcal N) \leq \frac{4\delta}9+ \eps^{-1} \lVert u_n - u\rVert_{\LL^1(\partial \B(a,r)} + \frac{M}{\kappa} + \sqrt{8M\eps} < \delta
	\end{equation} 
	for all $(t,x) \in [0,1] \times \partial \B(a,r)$,
	\( \dist(\tilde H_n(t,x); \mathcal N) < \delta \text{ for all } (t,x) \in [0,1] \times \partial \B(a,r).
	\) 
	Thus we see that $H_n(t,x) = \Pi_{\mathcal N}\circ\tilde H_n(t,x)$ is a well-defined homotopy between $u_n$ and $u$ for sufficiently large $n$, thereby establishing the result.
\end{Proof}

\subsection{Constrained compactness}\label{sec:proof_constrained}

\begin{Proof}{Theorem}{thm:convsub}
  
We will break up the proof into a series of steps.
Steps \hyperref[sec:step-1-extension-and-approximation]{1}--\hyperref[sec:step-4-limit-map-is-renormalisable]{4} are devoted to obtaining a limiting map $u_{\ast}$ which is renormalisable in $\Omega$, along with \eqref{item:convsubseq} and \eqref{item:limit_renormalisable}.
This involves approximating $(u_n)_n$ by mappings  in the $\mathrm{R}^{1,2}$-class (see Definition \ref{defn:Rclass}) in Step \hyperref[sec:step-1-extension-and-approximation]{1}, to which we can apply Proposition \ref{prop:firstorderconve}, and from which a weakly${}^{\ast}$-convergent limit can be extracted in Step \hyperref[sec:step-2-existence-of-the-limiting-map]{2}.
It is then shown in Step \hyperref[sec:step-3-minimal-topological-resolution]{3} that the limiting map has singularities forming a minimal topological resolution of $g$, and furthermore in Step \hyperref[sec:step-4-limit-map-is-renormalisable]{4} it is shown that $u_{\ast}$ is renormalisable and that singularities are bounded away from $\partial\Omega$.
In Step \hyperref[sec:step-5-localized-asymptotic-estimates]{5}, we prove the asymptotic estimates \eqref{item:isboundaway} and \eqref{item:isboundloc} for $f_n(\lvert \D u_n\rvert)$ , from which we can infer convergence of the defect measures, namely \eqref{item:vaguemeasures}, which we establish in Step \hyperref[sec:step-6-strict-measure-convergence]{6}.
Finally in Step \hyperref[sec:bounds-in-ll1]{7}, we use \eqref{item:isboundloc} to prove the $\LL^1$-estimates \eqref{item:L1bound}, which is used to establish weak $\LL^1$ convergence \eqref{item:weakconvergence}.

\subsubsection{Step 1: Extension and approximation}\label{sec:step-1-extension-and-approximation}

\begin{lemma}\label{lemma:density_of_the_R_class}
	Fix an open bounded Lipschitz domain $\Omega \subset \R^2$, a compact Riemannian manifold $\mathcal N$ and a map $u \in \WW^{1,1}(\Omega;\mathcal N)$. Assume that there exists an open set $\omega \subset\subset \Omega$ with $u \in \WW^{1,2}(\Omega\setminus  \bar\omega,\mathcal N)$.
	Then, for every $\varepsilon> 0$ and every $\theta \in (0,1)$, setting
	\[\omega_\theta \doteq \{x\in \Omega : \dist(\omega,x) < \theta\, \dist(\omega,\Omega)\},\] there exists a map $v \in \WW^{1,1}(\Omega,\mathcal N)$  such that 
	\begin{enumerate}[\rm(i)]
		\item for some points   $\{a_i\}_{i = 1,\dots,k}\subset \omega_\theta$ ($k\in \N$) such that $v \in C^\infty(\Omega\setminus\{a_i\}_{i = 1,\dots,k},\mathcal N)$,
		\item $\|u - v\|_{\WW^{1,1}(\Omega)} \leq \varepsilon$,
		\item $\Trace_{\partial \Omega}u$ and $\Trace_{\partial \Omega}v$ are homotopic.
	\end{enumerate}
\end{lemma}

This follows by adapting the arguments in \cite[Thm.~2]{bethuel1991approximation}, noting that in the regions where $u$ lies in $\WW^{1,2}$, we can locally approximate by composing mollifications with a nearest point projection to $\mathcal N$ (see \cite[\S 4]{SchoenUhlenbeck1983}). This can be done up to the boundary by working with a local extension, which exists by \cite[Lem.~5.5]{MonteilEtAl2022}. Note that the analogous result for $p \in (1,2)$ was proven in \cite[Prop.~4.2]{vanschaftingen2023asymptotic}; in this case one does not need to assume the compactness of $\mathcal N$.

Since $g \in \WW^{\sfrac{1}{2},2}(\partial \Omega;\mathcal N)$, by \cite[Lem.~5.5]{MonteilEtAl2022} there exists a Lipschitz domain $\Omega' \subset \R^2$ such that $\Omega\Subset \Omega'$ and a map $U \in \WW^{1,2}(\Omega'\setminus \Omega;\mathcal N)$ such that $\Tr_{\partial \Omega}U = g$ and such that $\Esgnu(\Tr_{\partial\Omega'}U) = \Esgnu(g)$. 
    For each $n$, we define
		\begin{equation}\label{eq:extensionofthesequence}
			\bar u_n = \begin{cases}
				u_n & \text{ on }\Omega \\
				U &  \text{ on }\Omega'\setminus \Omega,
			\end{cases}
		\end{equation}
    then each $\bar u_n \in \WW^{1,1}(\Omega';\mathcal N)$, using \cite[Thm.~18.1(ii)]{Leoni2017} to verify the weak differentiability.
    Observe by the dominated convergence theorem and Lemma \ref{lemma:quadgrowth} (as in the proof of Proposition \ref{prop:upperBound}) that,
    \begin{equation}\label{eq:DU_conv}
			\lim_{n \to \infty}\int_{\Omega'\setminus \Omega}f_n(|\D \bar u_n|) = \lim_{n \to \infty}\int_{\Omega'\setminus \Omega}f_n(|\D  U|) =\int_{\Omega'\setminus \Omega}\frac{|\D U|^2}{2}.
    \end{equation}
  Now recalling the $\kappa$-truncation operator from Definition \ref{defn:kappa_truncation},
  noting that $T_{\kappa}f_n \leq f_n$ for all $n$ and $T_{\kappa}f_n \to f_n$ pointwise as $\kappa \to \infty$.
  Thus we can find a sequence $(\kappa_n)_n$ such that
  \begin{align}
    \lim_{n \to \infty} \left(\mathcal V(f_{n}) - \mathcal V(T_{\kappa_n}f_{n})\right)&=0 \label{eq:truncation_vf}\\
    \lim_{n \to \infty} \int_{\Omega'}\lvert f_n(\D \bar u_n) - T_{\kappa_n}f_n(\D \bar u_n) \rvert \d x &=0.\label{eq:trunctation_firstorder}
  \end{align} 
  Then combining with \eqref{eq:asymptotic_finiteness} and \eqref{eq:DU_conv} we deduce that
  \begin{equation}
    \begin{split}
      \varlimsup_{n \to \infty} \left[ \int_{\Omega'} T_{\kappa_n}f_n(\lvert \D \bar u_n\rvert) - \mathcal V(T_{\kappa_n}f_n) \Esgnu(g) \right]
                   &=  \mathrm{L} + \int_{\Omega' \setminus \Omega} \frac12 \lvert \D U\rvert^2 \d x,
    \end{split}
  \end{equation} 
  recalling $\mathrm L$ was defined in\eqref{eq:asymptotic_finiteness}.
  In particular,
	\begin{equation}\label{eq:firstorderbndndessbarlin}
    \varlimsup_{n\to\infty}\mathcal V(T_{\kappa_n}f_{n})^{-1}\int_{\Omega'} T_{\kappa_n}f_{n}(|\D \bar u_n|) = \Esgnu(g).
	\end{equation}
  Now putting $\delta \doteq \dist(\Omega,\partial\Omega')>0$, for each $n \in \N$ by Lemma \ref{lemma:density_of_the_R_class} we can find a sequence $(u_{n,m})_m \subset \RR^{1,2}_{\delta}(\Omega'; \mathcal N)$ such that $u_{n,m} \to \bar u_n$ strongly in $\WW^{1,1}(\Omega'; \mathcal N)$ as $m \to \infty$ and such that $\Esgnu(\Tr_{\partial\Omega'}u_{n,m}) = \Esgnu(g)$ for all $m$.
  Then for each $n$, by passing to a subsequence in $m$ if necessary, by the partial converse of the dominated convergence theorem  (see for instance \cite[Thm.~4.9]{Brezis2011}), we can find $g_n \in \LL^1(\Omega',\mathbb R)$ non-negative such that $\lvert \D u_{n,m}\rvert \leq g_n$ for all $m$.

  Therefore, by the linear growing behaviour of $T_{\kappa_n}f_n$ and the dominated convergence theorem (using the majorant $g_n$),
	\[
    \lim_{m\to \infty}\int_{\Omega}|T_{\kappa_n}f_n(|\D u_{n,m}|)  -T_{\kappa_n}f_n(|\D u_{n}|)| \d x  = 0 \quad\mbox{for all $n$}.
	\] 
	Thus we choose a sequence $(m_n)_n$ such that 
  \begin{equation}\label{eq:densityinfirstorder}
    \lim_{n \to \infty} \int_{\Omega} \lvert T_{\kappa_n}f_n( \lvert\D u_{n,m_n}\rvert )  -T_{\kappa_n}f_n(\lvert\D \bar u_{n}\rvert)\rvert \d x = 0,
  \end{equation} 
  so we have
  \begin{equation}\label{eq:truncation_2nd_error}
    \varlimsup_{n \to \infty} \left[ \int_{\Omega'} T_{\kappa_n}f_n(\lvert \D u_{n,m_n}\rvert) - \mathcal V(T_{\kappa_n}f_n) \Esgnu(g) \right]
                   =  \mathrm{L} + \int_{\Omega' \setminus \Omega} \frac12 \lvert \D U\rvert^2 \d x.
  \end{equation} 
  We can also ensure that $\lvert \D \bar u_n - \D u_{n,m_n} \rvert \to 0$ strongly in $\LL^1(\Omega')$.
  Hence by applying Proposition \ref{prop:firstorderconve}\eqref{item:2ndorder} to $(u_{n,m_n})_n$ and $(T_{\kappa_n}f_n)_n$, for each $\eta \in (0, \min\{\frac{\sys}{2\pi t_0},\delta\})$ there exists a collection $\mathcal B_n^{\mathrm{Top}}(\eta)$ such that putting
  \begin{equation}\label{eq:sum_energy_notation}
    \mathrm{E}(u_{n,m_n},\mathcal B_n^{\mathrm{Top}}(\eta)) \doteq \sum_{\B \in \mathcal B_n^{\mathrm{Top}}(\eta)} \Esgnu(\Tr_{\partial \B}u_{n,m_n})
  \end{equation} 
  we have
  \begin{equation}\label{eq:ballmerging_fixed_n}
    \mathrm{E}(u_{n,m_n},\mathcal B_n^{\mathrm{Top}}(\eta)) \Lambda_{T_{\kappa_n}f_n}\left( \frac{\eta}{\mathrm{E}(u_{n,m_n},\mathcal B_n^{\mathrm{Top}}(\eta))} \right)  \leq \int_{\bigcup \mathcal B_n^{\mathrm{Top}}(\eta)} f_n(\lvert \D u_n\rvert) \d x,
  \end{equation} 
  \begin{equation}\label{eq:step1_esg_bound}
      \varlimsup_{n\to \infty}\mathrm{E}(u_{n,m_n},\mathcal B_n^{\mathrm{Top}}(\eta)) \leq \varlimsup_{n\to\infty}\mathcal V(T_{\kappa_n}f_n)^{-1}\int_{\bigcup\mathcal B_{n}(\eta)} f_n(|\D u_{n,m_n}|),
		\end{equation} 
    and combining with \eqref{eq:truncation_2nd_error},
    \begin{equation}\label{eq:step1_remainder_bound}
      \begin{split}
        \varlimsup_{n\to\infty}\int_{\Omega'\setminus \bigcup \mathcal B_n^{\mathrm{Top}}(\eta)}T_{\kappa_n}f_n(|\D u_{n,m_n}|) 
      \leq \mathrm{L} + \int_{\Omega' \setminus \Omega} \frac12 \lvert \D U\rvert^2 \d x  +\Esgnu(g)\log \frac{\sys\Esgnu(g)}{2\eta}.
    \end{split}
		\end{equation}
    We also have $\Esgnu(\Tr_{\partial B}u_{n,m_n}) > 0$ for each $\B \in \mathcal B_n^{\mathrm{Top}}(\eta)$.

\subsubsection{Step 2: Existence of a limiting map}\label{sec:step-2-existence-of-the-limiting-map}

  To simplify notation we will denote
  \begin{equation}\label{eq:notationxcvvbb}
    \tilde u_n \doteq u_{n,m_n} \in \RR^{1,2}_{\delta}(\Omega';\mathcal N), \quad \tildef_n(t) \doteq T_{\kappa_n}f_n(t).
  \end{equation} 
  For each $\B \in \mathcal B_n^{\mathrm{Top}}(\eta_n)$, we have $\Esgnu(\Tr_{\partial \B}u_n)>\frac{\sys^2}{4\pi}$ since it is non-zero (using the systole defined in \eqref{eq:systole}).
  Combining this with \eqref{eq:step1_esg_bound} we can bound the number of balls as 
	\begin{equation}\label{eq:bornesup}
    \limsup_{n \to \infty}\#\mathcal B_{n}^\mathrm{Top}(\eta)  \leq \frac{4\pi}{\sysN^2}\varlimsup_{n\to\infty}\mathcal V(\tildef_n)^{-1}\int_{\bigcup\mathcal B_{n}(\eta)} \tildef_n(|\D \tilde u_n|) = \frac{4\pi\Esgnu(g)}{\sysN^2}.
	\end{equation}
  Hence by passing to a subsequence in $n$ (which depends on $\eta>0$), we can assume $\mathcal B^{\mathrm{Top}}_n(\eta)$ converges to a limit collection $\mathcal B^{\mathrm{Top}}(\eta)$ in the following sense: we have $\#\mathcal B^{\mathrm{Top}}_n(\eta) \equiv \kappa_1$ is constant and writing $\mathcal B^{\mathrm{Top}}_n(\eta) = \{ \B(a_i^n,r_i^n)\}_{i=1}^{\kappa_1}$, we have $a_i^n \to \tilde a_i$ and $r_i^n \to \tilde r_i$ for all $1 \leq i \leq \kappa_1$.
  The limit collection $\mathcal B_{\ast}^{\mathrm{Top}}(\eta) = \{\B(a_i,r_i)\}_{i=1}^{\kappa}$ is then defined as these balls $\B(\tilde a_i,\tilde r_i)$, discarding any balls with $\tilde r_i = 0$, so in particular $\kappa \leq \kappa_1$.
  Note that if $\tilde a_i = \tilde a_j$ for some $i \neq j$, the disjointness of $\mathcal B^{\mathrm{Top}}_n(\eta)$ implies that either $\tilde r_i^n \to 0$ or $\tilde r_j^n \to 0$; thus we infer the centre points $\{a_i\}_{i=1}^{\kappa}$ are distinct, and that the limiting collection has pairwise disjoint interiors.

  Now let $(\eta_m)_m$ be a positive sequence such that $\eta_m \searrow 0$, then applying the above to each $\eta_m$ combined with a Cantor diagonal argument, we can find a common subsequence in $n$ such that for each $m$, $\mathcal B^{\mathrm{Top}}_n(\eta_m)$ converges to $\mathcal B_{\ast}^{\mathrm{Top}}(\eta_m)$ in the above sense. 
  Passing to a subsequence in $m$, we have the centres of $\mathcal B_{\ast}^{\mathrm{Top}}(\eta_m)$ converges to a set of finitely many points $\{a_i\}_{i=1}^{\kappa} \subset \overline\Omega$, such that $\kappa \leq \frac{4\pi \Esgnu(g)}{\sys^2}$ by \eqref{eq:bornesup}.

  We will let $\bar \rho_{\Omega'}  \doteq \rho_{\Omega'}(a_1,\dotsc,a_{\kappa})$ be the non-intersection radius as defined in \eqref{eq:nonintersection}.
  Then let $0 < \eps < \bar\rho_{\Omega'}$, and choose $m_0 \geq 0$ such that $\eta_m < \eps/2$ and $\bigcup \mathcal B_{\ast}^{\mathrm{Top}}(\eta_m) \subset \bigcup_{i=1}^{\kappa} \B(a_i,\eps/2)$ whenever $m \geq m_0$. 
  Fixing any such $n$, the above convergence we have $\bigcup \mathcal B_n^{\mathrm{Top}}(\eta_m) \subset \bigcup_{i = 1}^\kappa \B(a_i,\eps)$ for all $n$ sufficiently large.
  Thus from \eqref{eq:step1_remainder_bound}, for all $m\geq m_0$ we have,
	\begin{equation}\label{eq:bndawayfrmSingwithrho}
    \varlimsup_{n\to\infty}\int_{\Omega'\setminus \bigcup_{i = 1}^\kappa \B(a_i,\epsilon)}\tildef_n(|\D \tilde u_n|) \leq \tilde{\mathrm{L}}_2 + \Esgnu(g) \log\frac1{\eta_m},
  \end{equation}
  where we set
  \begin{equation}\label{eq:tildeE2}
    \tilde{\mathrm{L}}_2  \doteq\mathrm{L} + \Esgnu(g) \log \frac{\sys \Esgnu(g)}{2} + \int_{\Omega'\setminus\Omega} \frac12 \lvert \D U\rvert^2 \d x
  \end{equation} 
  Therefore by Proposition \ref{prop:compacitnoprblem}\eqref{item:unconstrained_bound} and a Cantor diagonal argument in $m$, we have $(\lvert \D \tilde u_n\rvert)_n$ is bounded in each $\LL^1(\Omega' \setminus \bigcup_{i=1}^{\kappa}\B(a_i,\eps))$. Moreover, since $\tilde u_n$ takes values in the compact manifold $\mathcal N$, the sequence is bounded in $\WW^{1,1}(\Omega' \setminus \bigcup_{i=1}^{\kappa}\B(a_i,\eps))$.
  Hence by taking a further diagonal sequence, we obtain a limit map $u_{\ast}$ such that passing to a subsequence, $\tilde u_n$ converges \weaklystar in $\BV$ to $u_{\ast}$ in $\Omega' \setminus \bigcup_{i=1}^{\kappa} \B(a_i,\eps)$. 
  Moreover \weakstar convergence in $\BV$ implies strong $\LL^1$ convergence of $\tilde u_n$, so Proposition \ref{prop:compacitnoprblem}\eqref{item:unconstrained_conv} we also have the estimate
  \begin{equation}\label{eq:G_lsc_away}
    \int_{\Omega' \setminus \bigcup_{i=1}^{\kappa} \B(a_i,\eps)} \frac12\lvert \D u_{\ast}\rvert^2 \d x \leq \varlimsup_{n \to \infty} \int_{\Omega' \setminus \bigcup_{i=1}^{\kappa} \B(a_i,\eps)} \tildef_n(\lvert \D \tilde u_n\rvert) \d x \leq \tilde{\mathrm{L}}_2 + \Esgnu(g) \log\frac1{\eta_m}
  \end{equation} 
  using \eqref{eq:bndawayfrmSingwithrho}.
  Since this holds for all $\eps>0$ sufficiently small, this defines a map in $\WW^{1,2}_{\mathrm{loc}}(\Omega' \setminus \{a_i\}_{i=1}^{\kappa};\mathcal N)$, noting the $\LL^1$ convergence ensures the limit map $u_{\ast}$ takes values in $\mathcal N$ for almost every $x \in \Omega$.
  Moreover since $\D \bar u_n - \D\tilde u_n \to 0$ strongly in $\LL^1$ and $\bar u_n \rvert_{\Omega' \setminus \Omega} \equiv U$ for all $n$, it follows that $u_{\ast} \rvert_{\Omega'\setminus\Omega } \equiv U$.

  Now combining \eqref{eq:G_lsc_away} with \eqref{eq:trunctation_firstorder} and \eqref{eq:densityinfirstorder}, we have
  \begin{equation}\label{eq:bndawayfrmSingwithrhoIII}
    \begin{split}
      & \varlimsup_{n \to \infty} \int_{\Omega \setminus \bigcup_{i=1}^{\kappa} \B(a_i,\eps)}  f_n(\lvert \D u_n\rvert) \d x\\
      &\quad= \varlimsup_{n \to \infty} \int_{\Omega' \setminus \bigcup_{i=1}^{\kappa} \B(a_i,\eps)} \tildef_n(\lvert \D \tilde u_n\rvert) \d x - \int_{\Omega' \setminus \left(\Omega \cup\bigcup_{i=1}^{\kappa} \B(a_i,\eps)\right)} \frac12 \lvert \D U\rvert^2 \d x\\
      &\quad\leq \tilde{\mathrm{L}}_{2}+\log\frac{1}{\eta_m} + \int_{(\Omega' \setminus \Omega) \cap \left(\bigcup_{i=1}^{\kappa} \B(a_i,\eps)\right)} \frac12 \lvert\D U\rvert^2 \d x.
    \end{split}
  \end{equation}

  \subsubsection{Step 3: Minimal topological resolution}\label{sec:step-3-minimal-topological-resolution}
  To show the singularities of $\tilde u_n$ are transferred to the limiting map $u_{\ast}$, we will rely on Lemma \ref{lemma:topologyonthesphere}.

  Given $0 < \eps < \bar\rho_{\Omega'}$, by converting to polar coordinates and applying Fatou's lemma, we have
  \begin{equation}
    \int_{\eps}^{\bar\rho_{\Omega'}} \varliminf_{n \to \infty} \int_{\partial\B(a_i,r)} \tildef_n(\lvert \D \tilde u_n\rvert) \d \mathcal H^1 \d r \leq \varliminf_{n \to \infty}\int_{\B(a_i,\bar\rho_{\Omega'}) \setminus \B(a_i,\eps)} \tildef_n(\lvert \D \tilde u_n\rvert) \d x
  \end{equation} 
  is finite for each $i = 1,\dotsc,\kappa$ by \eqref{eq:bndawayfrmSingwithrhoIII}.
  Thus for almost every $r \in (0,\bar\rho_{\Omega'})$ we have
	\begin{equation}\label{eq:G_trace_finite}
		\varliminf_{n \to \infty}\int_{\partial \B(a_i,r)}\tildef_n(|\D \tilde u_n|)\d \mathcal H^1\quad \text{ is finite for each $i = 1,\dotsc,\kappa$.}
	\end{equation}
  By Proposition \ref{prop:trace_conv}, there exists a subsequence $\tilde u_{n_k}$ (depending on $r$) attaining the limit inferior and along which $\tr_{\partial \B(a_i,r)}\tilde u_{n_k} \to \tr_{\partial \B(a_i,r)}u_*$ in $\LL^1(\partial \B(a_i,r);\mathcal N)$. By applying Proposition \ref{prop:compacitnoprblem} on each $\mathbb S^1(a_i,r)$, we deduce that for almost every $r \in (0, \bar \rho)$ and each $i = 1,\dots,\kappa$,
	\begin{equation}\label{eq:lsc_on_spheres}
		\int_{\partial \B(a_i,r)}\frac{|\D u_*|^2}{2} \d \mathcal H^1\leq \varliminf_{n \to \infty}\int_{\partial \B(a_i,r)}\tildef_n(|\D \tilde u_n|) \d \mathcal H^1.
	\end{equation}
  Now applying Lemma \ref{lemma:topologyonthesphere} to $\Tr_{\partial \B(a,r)}\tilde u_{n_k}$ and $\tildef_{n_k}$ we have
  \begin{equation}\label{eq:lambda_constant2}
    \lambda(\Tr_{\partial\B(a_i,r)}\tilde u_{n_k}) = \lambda(\Tr_{\partial\B(a_i,r)}u_{\ast})
  \end{equation}
  for all $1 \leq i \leq \kappa$ and $k$ sufficiently large.

  Given $\eps>0$ such that \eqref{eq:G_trace_finite} holds with $r = \eps$, along the subsequence such that \eqref{eq:lambda_constant2} holds, for $k$ sufficiently large by Lemma \ref{lem:annular_lowerbound} we have
  \begin{equation}
    \sum_{i=1}^{\kappa} \frac{\lambda(\Tr_{\partial \B(a_i,\eps)} u_{\ast})^2}{4\pi}\Lambda_{\tildef_{n_k}}\left( \frac{4 \pi \eps}{\lambda(\Tr_{\partial \B(a_i,\eps)} u_{\ast})^2} \right) \leq \int_{\bigcup_{i=1}^{\kappa} \B(a_i,\eps)} \tildef_{n_k}(\lvert \D \tilde u_{n_k}\rvert)\d x.
  \end{equation} 
  Dividing both sides by $\mathcal V(\tildef_{n_k})$ and sending $k \to \infty$, using Lemma \ref{lemma:asympLambda} and \eqref{eq:truncation_2nd_error} we have
  \begin{equation}
    \begin{split}
      \sum_{i=1}^{\kappa} \frac{\lambda(\Tr_{\partial \B(a_i,\eps)} u_{\ast})^2}{4\pi}
      &= \sum_{i=1}^{\kappa} \frac{\lambda(\Tr_{\partial \B(a_i,\eps)} u_{\ast})^2}{4\pi}\lim_{k \to \infty}  \mathcal V(\tildef_{n_k})^{-1}\Lambda_{\tildef_{n_k}}\left( \frac{4 \pi \eps}{\lambda(\Tr_{\partial \B(a_i,\eps)} u_{\ast})^2} \right) \\
      &\leq \varlimsup_{n \to \infty}   \frac{1}{\mathcal V(\tildef_n)}\int_{\Omega} \tildef_n(\lvert \D \tilde u_n\rvert) \d x
      = \Esgnu(g).
    \end{split}
  \end{equation} 
  By subadditivity of the singular energy (namely \eqref{eq:esg_subadditive}) and noting that $u_{\ast} \in \WW^{1,2}_{\mathrm{loc}}(\Omega' \setminus \{a_i\}_{i=1}^{\kappa};\mathcal N)$, the above holds with equality for all $\eps \in (0,\bar\rho_{\Omega'})$. 
  Thus
  \begin{equation}\label{eq:uast_minimal_top}
    \sum_{i=1}^{\kappa} \frac{\lambda([u_{\ast},a_i])^2}{4\pi} = \sum_{i=1}^{\kappa} \frac{\lambda(\Tr_{\partial \B(a_i,\eps)} u_{\ast})^2}{4\pi} = \Esgnu(g) = \Esgnu(\Tr_{\partial \Omega'}U),
  \end{equation} 
  and so $(\Tr_{\mathbb S^1}(a_i+ \eps \,\cdot))_{i=1}^{\kappa}$ forms a minimal topological resolution for $\Tr_{\partial\Omega'}U$ in $\Omega'$ for all $\eps \in (0,\bar\rho_{\Omega'})$ sufficiently small.
  Since $\Tr_{\partial\Omega}u_{\ast} = g$, \emph{if} each $a_i \in \Omega$, then for $\eps>0$ we will also have that $(\Tr_{\mathbb S^1}(a_i+ \eps \,\cdot))_{i=1}^{\kappa}$ forms a minimal topological resolution for $g$ in $\Omega$.

  \subsubsection{Step 4: Limit map is renormalisable}\label{sec:step-4-limit-map-is-renormalisable}

  Next, we will show that the limit map is renormalisable in $\Omega'$.
  By \eqref{eq:lsc_on_spheres}, for almost all $r \in (0,\bar \rho_{\Omega'})$ we have
  \begin{equation}\label{eq:renom_slicewise}
    \int_{\partial \B(a_i,r)} \frac1{2}{|\D u_*|^2}\d \mathcal H^1 - \frac{\lambda([u_{\ast},a_i])^2}{4\pi r} \leq \varliminf_{n\to\infty}\int_{\partial \B(a_i,r)}\tildef_n(|\D \tilde u_n|) - \tildef_n\left(\frac{\lambda([u_{\ast},a_i])}{2\pi r}\right)\d \mathcal H^1,
  \end{equation}
  where we used the dominated convergence theorem and Lemma \ref{lemma:quadgrowth} to pass to the limit in the second term.
  Now for $0 <  \rho < \bar\rho_{\Omega'}$, integrating \eqref{eq:renom_slicewise} over $r \in (0,\rho)$, summing over $i = 1, \dotsc, \kappa$ and applying Fatou's lemma, we have
  \begin{equation}\label{eq:u_ast_renormalisable}
    \begin{split}
    &\sum_{i=1}^{\kappa}\int_{0}^{\rho}\left[\int_{\partial\B(a_i,r)} \frac{\lvert \D u_{\ast}\rvert^2}2 \d\mathcal H^1 - \frac{\lambda([u_{\ast},a_i])^2}{4\pi r} \right] \d r\\
    &\leq \sum_{i=1}^{\kappa}\int_{0}^{\rho}\varliminf_{n \to \infty} \left[ \int_{\partial \B(a_i,r)}\tildef_n(|\D \tilde u_n|) - \tildef_n\left(\frac{\lambda([u_{\ast},a_i])}{2\pi r}\right)\d \mathcal H^1\right] \d r \\
    &\leq \varliminf_{n \to \infty} \left[\int_{\bigcup_{i=1}^{\kappa}\B(a_i,\rho)} \tildef_n(\lvert \D \tilde u_n\rvert)\d x - \sum_{i=1}^{\kappa}\frac{\lambda([u_{\ast},a_i])^2}{4\pi}\int_0^{\frac{2\pi \rho}{\lambda([u_{\ast},a_i])}} \tildef_n(1/r)r \d r \right]\\
    &= \varliminf_{n \to \infty} \left[\int_{\bigcup_{i=1}^{\kappa}\B(a_i,\rho)} \tildef_n(\lvert \D \tilde u_n\rvert)\d x - \Esgnu(g)\mathcal V(\tildef_n) \right] 
    - \sum_{i=1}^{\kappa} \frac{\lambda([u_{\ast},a_i])^2}{4\pi} \log\frac{2\pi\rho}{\lambda([u_{\ast},a_i])}
    \end{split}
  \end{equation} 
  where we used Lemma \ref{lemma:entropy_error} and \eqref{eq:uast_minimal_top} in the last line, which is moreover finite by \eqref{eq:truncation_2nd_error}.
  Hence $u_{\ast}$ verifies property \eqref{item:actualuseofren} of Proposition \ref{prop:renormalisedenergy}, from which we infer that $u_{\ast} \in \WW^{1,2}_{\ren}(\Omega';\mathcal N)$.

  To show the singular points $a_i$ actually lie in $\Omega$,
  we will make use of the following lemma from \cite[Lem.~6.2]{MonteilEtAl2022}.
\begin{lemma}\label{lemma:ainomega}
	Let $\Omega \subset \R^2$ be an open set, $a \in \R^2$ and $0 <  \sigma < \tau$. If $u  \in \WW^{1,2}(\B(a,\tau)\setminus \B(a,\sigma); \mathcal N)$ is such that $\lambda(\Tr_{\partial \B(a,\tau)}u)>0$,
  then
	\begin{equation*}
    \Gamma_{\tau}^{\sigma}(a,\D u)	\frac{\lambda(\Tr_{\partial \B(a,\tau)}u)^2}{4\pi \nu_{\tau}^\sigma(a)}\log \frac{\tau}{\sigma}  \leq \int_{\Omega \cap \B(a,\tau) \setminus \B(a,\sigma)} \frac{|\D u|^2}{2} ,
	\end{equation*}
	where we have set
	\begin{equation*}
    \Gamma_{\tau}^{\sigma}(a,\D u) \doteq \Big [1 - \frac{\sqrt{2\pi}}{\lambda(\Tr_{\partial \B(a,\tau)}u)}\Big (\frac{1}{\log \frac{\tau}{\sigma}} \int_{\B(a,\tau) \setminus (\Omega \cup \B(a,\sigma))} |\D u|^2 \Big)^{\frac{1}{2}} \Big]^2
	\end{equation*}
and
	\begin{equation}\label{eq:densityboundary}
		\nu_{\tau}^\sigma(a)  \doteq \frac{1}{2\pi \log \frac{\tau}{\sigma}}\int_{\Omega \cap \B(a,\tau) \setminus \B(a,\sigma)}\frac{\d x}{|x - a|^2} \leq 1.
	\end{equation}
\end{lemma}
This will be used in conjunction with the following observation from \cite[\S 6]{MonteilEtAl2022}.
\begin{lemma}\label{lemma:point_in_omega}
  Let $\Omega \subset \mathbb R^2$ be a bounded Lipschitz domain, and let $a \in \mathbb R^2$. Then for any $\tau>0$, we have
  \(
    \lim_{\sigma \to 0} \nu_{\tau}^{\sigma}(a) = 1\) if and only if \(a \in \Omega.
  \) 
\end{lemma}
\begin{Proof}{Lemma}{lemma:point_in_omega}
  Observe that for any $\tau_1, \tau_2 >0$ we have
  \(
    \lim_{\sigma \to 0} (\nu_{\tau_1}^{\sigma}(a) - \nu_{\tau_2}^{\sigma}(a)) = 0,
  \) 
  so it suffices to prove the result for some $\tau>0$.
  If $a \in \Omega$ or $a \notin \overline\Omega$, by taking $\tau$ sufficiently small so $\B(a,\tau) \subset \Omega$ or $\B(a,\tau)\subset\mathbb R^2\setminus\overline\Omega$, we have $\lim_{\sigma \to 0} \nu_{\tau}^{\sigma}(a) = 1$ or $0$ respectively.
  If $a \in \partial\Omega$, by the exterior cone condition (see \cite[Thm.~1.2.2.2]{Grisvard2011}) there exists $\sigma_0>0$, $v \in \mathbb S^1$ and $\theta\in (0,\pi/2)$ such that
  \begin{equation}
    C_{\sigma_0}(a,\theta,v) \doteq \{ a + w : 0 < \lvert w\rvert< \sigma_0, \lvert\measuredangle(v,w)\rvert <\theta/2  \} \subset \B(a,\sigma_0) \setminus \Omega,
  \end{equation} 
  where $\measuredangle(v,w) = \arccos{((v\cdot w)/(|v||w|))}$ denotes the angle between the vectors $v$ and $w$ in $\mathbb R^2$.
  Then choosing $\tau = \sigma_0$ we have
  \begin{equation}
    \varlimsup_{\sigma \to 0} \nu_{\sigma_0}^{\sigma}(a) \leq \frac1{2\pi \log(\sigma_0/\sigma)} \int_{\B(a,\sigma_0) \setminus ((\B(a,\sigma) \cup C_{\sigma_0}(a,\theta,v))} \frac{\dd x}{\lvert x- a\rvert^2} =1 - \frac{\theta}{2\pi} < 1,
  \end{equation} 
  thereby establishing the result.
\end{Proof}

  By \eqref{eq:u_ast_renormalisable}, for $0 < \sigma<\rho < \bar\rho_{\Omega'}$ and each $1 \leq i \leq \kappa$ 
  \begin{equation}
    \frac{\lambda([u_{\ast},a_i])^2}{4\pi} \log \frac{\rho}{\sigma} \leq \int_{\B(a_i,\rho) \setminus \B(a_i,\sigma)} \frac{\lvert \D u_{\ast}\rvert^2}2 \d x \leq \frac{\lambda([u_{\ast},a_i])^2}{4\pi} \log \frac{\rho}{\sigma} + C_{\rho},
  \end{equation} 
  where $C_{\rho}$ is an upper bound for \eqref{eq:u_ast_renormalisable}, and the lower bound follows by converting to polar coordinates and recalling the definition of $\lambda([u_{\ast},a_i])$.
  From this we deduce that
  \begin{equation}\label{eq:log_uast_concentration}
    \lim_{\sigma \to 0} \frac1{\log(\rho/\sigma)} \int_{\B(a_i,\rho)\setminus\B(a_i,\sigma)} \frac{\lvert \D u_{\ast}\rvert^2}2 \,\d x = \frac{\lambda([u_{\ast},a_i])^2}{4\pi}.
  \end{equation} 
  Moreover since $u_{\ast}\rvert_{\Omega'\setminus\Omega} \equiv U \in \WW^{1,2}(\Omega'\setminus\Omega;\mathcal N)$, we have
  \begin{equation}\label{eq:log_boundary_concentration}
    \lim_{\sigma \to 0} \frac1{\log(\rho/\sigma)} \int_{\B(a_i,\rho)\setminus(\Omega \cup \B(a_i,\sigma))} \frac{\lvert \D u_{\ast}\rvert^2}2 \,\d x = 0,
  \end{equation} 
  so in particular, $\lim_{\sigma \to 0} \Gamma_{\rho}^{\sigma}(a_i,\D u_{\ast}) = 1$.
  Combining \eqref{eq:log_uast_concentration} and \eqref{eq:log_boundary_concentration} we obtain
  \begin{align}
    \lim_{\sigma \to 0} \frac1{\log(\rho/\sigma)} \int_{\Omega \cap \B(a_i,\rho)\setminus\B(a_i,\sigma)} \frac{\lvert \D u_{\ast}\rvert^2}2 \,\d x = \frac{\lambda([u_{\ast},a_i])^2}{4\pi},
  \end{align} 
  so by Lemma \ref{lemma:ainomega} we have
  \begin{equation}
    \varliminf_{\sigma \to 0} \nu_{\rho}^{\sigma}(a_i) \geq \lim_{\sigma \to 0}(\Gamma_{\rho}^{\sigma}(a_i,\D u_{\ast}))^{-1} \frac{4\pi}{\lambda([u_{\ast},a_i])^2} \int_{\Omega \cap \B(a_i,\rho) \setminus \B(a_i,\sigma)} \frac{\lvert \D u_{\ast}\rvert^2}2 \,\d x = 1
  \end{equation} 
  for each $i = 1 ,\dotsc, \kappa$.
  It follows that $\lim_{\sigma \to 0} \nu_{\rho}^{\sigma}(a_i) = 1$, and hence each $a_i \in \Omega$ by Lemma \ref{lemma:point_in_omega}.
  Therefore the restriction of $u_{\ast}$ to $\Omega$ is also a renormalisable map, establishing \eqref{item:convsubseq}.

  We can also compute the renormalised energy of $u_{\ast}$ in $\Omega$ using \eqref{eq:erenwithoutlimit}, \eqref{eq:uast_minimal_top} and \eqref{eq:u_ast_renormalisable} as, for any $\rho < \rho_{\Omega}(a_1,\dotsc,a_{\kappa})$,
  \begin{equation*}
    \begin{split}
      \mathrm{E}_{\ren}^{1,2}(u_{\ast};\Omega) 
      &= \int_{\Omega \setminus \bigcup_{i=1}^{\kappa}\B(a_i,\rho)} \frac{\lvert \D u_{\ast}\rvert^2}2 \,\d x - \sum_{i=1}^{\kappa} \frac{\lambda([u_{\ast},a_i])^2}{4 \pi} \log\frac1{\rho} \\
      &\quad+ \sum_{i=1}^{\kappa} \int_0^\rho \left[\int_{\partial\B(a_i,r)} \frac{\lvert \D u_{\ast}\rvert^2}2 \d\mathcal H^1 - \frac{\lambda([u_{\ast},a_i])^2}{4\pi r}\right] \d r \\
      &\leq \varliminf_{n \to \infty} \int_{\Omega \setminus \bigcup_{i=1}^{\kappa}\B(a_i,\rho)} \tildef_n(\D \tilde u_n) \d x - \sum_{i=1}^{\kappa} \frac{\lambda([u_{\ast},a_i])^2}{4\pi} \log\frac{2\pi\rho}{\lambda([u_{\ast},a_i])}\\
      &\quad+ \varliminf_{n \to \infty} \left[\int_{\bigcup_{i=1}^{\kappa}\B(a_i,\rho)} \tildef_n(\lvert \D \tilde u_n\rvert)\d x - \Esgnu(g)\mathcal V(\tildef_n) \right] - \Esgnu(g)\log\frac1{\rho}  \\
      &\leq \varliminf_{n \to \infty} \left[\int_{\Omega} \tildef_n(\lvert \D \tilde u_n\rvert)\d x - \Esgnu(g)\mathcal V(\tildef_n) \right] - \sum_{i=1}^{\kappa} \frac{\lambda([u_{\ast},a_i])^2}{4\pi} \log\frac{2\pi}{\lambda([u_{\ast},a_i])}\\
      &\leq \varliminf_{n \to \infty} \left[\int_{\Omega} f_n(\lvert \D u_n\rvert)\d x - \Esgnu(g)\mathcal V(f_n) \right] - \sum_{i=1}^{\kappa} \frac{\lambda([u_{\ast},a_i])^2}{4\pi} \log\frac{2\pi}{\lambda([u_{\ast},a_i])},
    \end{split}
  \end{equation*} 
  where the last line follows by \eqref{eq:truncation_2nd_error}, thereby establishing \eqref{item:limit_renormalisable}.

\subsubsection{Step 5: Localized asymptotic estimates}\label{sec:step-5-localized-asymptotic-estimates}

For $0<\rho < \bar\rho_{\Omega} \doteq \rho_{\Omega}(a_1,\dotsc,a_{\kappa})$, by integrating \eqref{eq:renom_slicewise} over $(0,\rho)$ and applying Fatou's lemma we have
\begin{equation}
  \begin{split}
  &\varliminf_{n\to\infty} \left[\int_{\B(a_i,\rho)} \tildef_n(\lvert\D \tilde u_n\rvert)\d x - 2\pi\int_0^{\rho} \tildef_n\left( \frac{\lambda([u_{\ast},a_i])}{2\pi r} \right)r \d r \right] \\
  &\quad\geq \int_{0}^{\rho}\left[\int_{\partial\B(a_i,r)} \frac12 \lvert \D u_{\ast}\rvert^2 \d\mathcal H^1 - \frac{\lambda([u_{\ast},a_i])^2}{4\pi r}\right] \d r \\
  &\quad = \mathrm{E}^{1,2}_{\ren}(u_{\ast};\B(a_i,\rho)) + \frac{\lambda([u_{\ast},a_i])^2}{4\pi} \log \frac1{\rho},
\end{split}
\end{equation} 
recalling the localised renormalised energy was defined in \eqref{eq:localised_renormalised}.
Then Lemma \ref{lemma:entropy_error} gives
\begin{equation}\label{eq:lowerbound_perturb_vfn}
  \lim_{n \to \infty} \left[2\pi \int_0^{\rho} \tildef_n\left( \frac{\lambda([u_{\ast},a_i])^2}{2\pi r} \right) r \d r
  - \frac{\lambda([u_{\ast},a_i])^2}{4\pi} \mathcal V(\tildef_n)\right] =  \frac{\lambda([u_{\ast},a_i])^2}{4\pi} \log\frac{2\pi \rho}{\lambda([u_{\ast},a_i])},
\end{equation} 
which we combine with \eqref{eq:truncation_vf}, \eqref{eq:trunctation_firstorder} and \eqref{eq:densityinfirstorder} to obtain
\begin{equation}\label{eq:ball_concerntration_bound}
  \begin{split}
    &\varliminf_{n \to \infty} \left[\int_{\B(a_i,\rho)} f_n(\lvert \D u_n\rvert) \d x - \frac{\lambda([u_{\ast},a_i])^2}{4\pi} \mathcal V(f_n)\right] \\
    &\qquad=\varliminf_{n \to \infty} \left[\int_{\B(a_i,\rho)} \tildef_n(\lvert \D \tilde u_n\rvert) \d x - \frac{\lambda([u_{\ast},a_i])^2}{4\pi} \mathcal V(\tildef_n)\right] \\
     &\qquad\geq \mathrm{E}^{1,2}_{\ren}(u_{\ast};\B(a_i,\rho)) + \frac{\lambda([u_{\ast},a_i])^2}{4\pi} \log\frac{2\pi}{\lambda([u_{\ast},a_i])}.
\end{split}
\end{equation} 
This establishes the lower bound in \eqref{item:isboundloc}.
Using this we can estimate
\begin{equation}
  \begin{split}
    &\varlimsup_{n \to \infty} \int_{\Omega \setminus \bigcup_{i=1}^{\kappa} \B(a_i,\rho)} f_n(\lvert \D u_n\rvert) \d x\\
    &\quad\leq \mathrm L - \sum_{i=1}^{\kappa}\varliminf_{n \to \infty} \left[\int_{\B(a_i,\rho)} f_n(\lvert \D u_n\rvert) \d x - \frac{\lambda([u_{\ast},a_i])^2}{4\pi} \mathcal V(f_n)\right] \\
    &\quad\leq \mathrm L - \sum_{i=1}^{\kappa}\frac{\lambda([u_{\ast},a_i])^2}{4\pi}\log\frac{2\pi}{\lambda([u_{\ast},a_i])} - \sum_{i=1}^{\kappa} \mathrm{E}_{\ren}^{1,2}(u_{\ast};\B(a_i,\rho)) \\
    &\quad= \int_{\Omega \setminus \bigcup_{i=1}^{\kappa} \B(a_i,\rho)} \frac12\lvert \D u_{\ast}\rvert^2 \d x+ \left(\mathrm{L} - \mathrm E_{\ren}^{1,2}(u_{\ast}) - \mathrm H([u_{\ast},a_i])_{i=1}^{\kappa}\right),
  \end{split}
\end{equation} 
where we used \eqref{eq:erenwithoutlimit} to write 
\begin{equation}\label{eq:renormalisable_additivity}
  \mathrm{E}_{\ren}^{1,2}(u_{\ast}) = \sum_{i=1}^{\kappa} \mathrm{E}_{\ren}^{1,2}(u_{\ast};\B(a_i,\rho)) + \int_{\Omega \setminus \bigcup_{i=1}^{\kappa} \B(a_i,\rho)} \frac12 \lvert \D u_{\ast}\rvert^2 \d x.
\end{equation} 
By Proposition \ref{prop:compacitnoprblem} we also have the lower bound
\begin{equation}\label{eq:away_lower_bound}
  \int_{\Omega \setminus \bigcup_{i=1}^{\kappa}\B(a_i,\rho)} \frac12 \lvert \D u_{\ast}\rvert^2 \,\d x \leq \liminf_{n \to \infty}\int_{\Omega \setminus \bigcup_{i=1}^{\kappa}\B(a_i,\rho)} f_n( \lvert \D u_{\ast}\rvert) \,\d x ,
\end{equation} 
which establishes \eqref{item:isboundaway}.
To obtain asymptotic upper bounds on $\B(a_i,\rho)$, we use \eqref{eq:renormalisable_additivity}, \eqref{eq:away_lower_bound} and \eqref{eq:ball_concerntration_bound} applied for each $j \neq i$ to estimate
\begin{equation*}
  \begin{split}
    &\varlimsup_{n \to \infty} \left[\int_{\B(a_i,\rho)} f_n(\lvert \D u_n\rvert) \d x - \frac{\lambda([u_{\ast},a_i])^2}{4\pi} \mathcal V(f_n)\right]\\
    &\quad\leq \mathrm{L} - \sum_{j \neq i}\varliminf_{n \to \infty}  \left[\int_{\B(a_j,\rho)} f_n(\lvert \D u_n\rvert) \d x- \frac{\lambda([u_{\ast},a_i])^2}{4\pi} \mathcal V(f_n)\right] - \varliminf_{n \to \infty} \int_{\Omega \setminus \bigcup_{i=1}^{\kappa} \B(a_i,\rho)} f_n(\lvert \D u_{\ast}\rvert) \d x  \\
    &\quad\leq \mathrm{L} - \sum_{j\neq i} \mathrm{E}_{\ren}^{1,2}(u_{\ast},\B(a_j,\rho)) - \sum_{j \neq i} \frac{\lambda([u_{\ast},a_j])^2}{4\pi} \log \frac{2\pi}{\lambda([u_{\ast},a_j])} - \int_{\Omega \setminus \bigcup_{i=1}^{\kappa} \B(a_i,\rho)} \frac12\lvert \D u_{\ast}\rvert^2 \d x \\
    &\quad = \mathrm{L} - \Eren^{1,2}(u_{\ast}) - \mathrm{H}([u_{\ast},a_i]) + \mathrm{E}^{1,2}_{\ren}(u_{\ast};\B(a_i,\rho)) + \frac{\lambda([u_{\ast},a_i])^2}{4\pi} \log \frac{2\pi}{\lambda([u_{\ast},a_i])},
  \end{split}
\end{equation*} 
thereby establishing \eqref{item:isboundloc}.

\subsubsection{Step 6: Strict measure convergence}\label{sec:step-6-strict-measure-convergence}

Multiplying the inequalities in \eqref{item:isboundloc} by $\mathcal V(f_n)^{-1}$, for each $0 < \rho < \bar\rho_{\Omega}$ and $i=1,\dotsc,\kappa$ we have
\begin{equation}\label{eq:localised_asymptoptic}
  \lim_{n \to \infty} \frac1{\mathcal V(f_n)} \int_{\B(a_i,\rho)} f_n(\lvert \D u_n\rvert) \d x = \frac{\lambda([u_{\ast},a_i])^2}{4 \pi}.
\end{equation} 
Also by \eqref{eq:asymptotic_finiteness} and \eqref{item:convsubseq} we know that
\begin{equation}\label{eq:global_asymptoptic}
  \lim_{n \to \infty} \frac1{\mathcal V(f_n)} \int_{\Omega} f_n(\lvert \D u_n\rvert) \d x = \Esgnu(g) = \sum_{i=1}^{\kappa} \frac{\lambda([u_{\ast},a_i])^2}{4\pi}.
\end{equation} 
Hence if we define the measures
\begin{equation}
  \mu_n \doteq \frac{f_n(\lvert \D u_n\rvert)}{\mathcal V(f_n)}  \Leb^2 \mres \Omega,
\end{equation} 
by \eqref{eq:localised_asymptoptic}, for all $0<\rho<\bar\rho_{\Omega}$ and $i=1,\dotsc,\kappa$ we have
\begin{equation}
  \lim_{n \to \infty} \mu_n(\B(a_i,\rho)) = \frac{\lambda([u_{\ast},a_i])^2}{4\pi},
\end{equation} 
whereas by \eqref{eq:global_asymptoptic} we have
\begin{equation}
  \lim_{n \to \infty} \mu_n(\Omega) = \sum_{i = 1}^\kappa\frac{\lambda([u_{\ast},a_i])^2}{4 \pi}.
\end{equation} 
Hence we infer that $\mu_n$ converges narrowly to a measure $\mu$ that is supported on $\{a_i\}_{i=1}^{\kappa}$, taking the form
\begin{equation}
  \mu = \sum_{i=1}^{\kappa} \frac{\lambda([u_{\ast},a_i])^2}{4\pi} \delta_{a_i},
\end{equation} 
thereby establishing \eqref{item:vaguemeasures}.

\subsubsection{Step 7 : Bounds in $\LL^1$}\label{sec:bounds-in-ll1}

For $0 < \rho < \rho_{\Omega}(a_1,\dotsc,a_{\kappa})$ and each $i = 1,\dotsc,\kappa$, set
\begin{equation}
  \mathrm{L}_{i,\rho} \doteq \varlimsup_{n \to \infty}\left[\int_{\B(a_i,\rho)} f_n(\lvert \D u_n\rvert) \d x - \frac{\lambda([u_{\ast},a_i])^2}{4\pi} \mathcal V(f_n)\right] ,
\end{equation} 
which we know is finite by \eqref{item:isboundloc}.
Then by convexity of $f_n$, using \eqref{eq:convex_linearisation} we have 
  \begin{equation}\label{eq:conv_firstorder}
    f_n(|\D u_n|) - f_n\Big(\frac{\lambda([u_*,a_i])}{2\pi r}\Big)\geq f_{n,+}'\Big(\frac{\lambda([u_*,a_i])}{2\pi r }\Big) \Big(|\D u_n| - \frac{\lambda([u_*,a_i])}{2\pi r}\Big)
	\end{equation}
  for all $x \in \B(a_i,\rho)$ with $r= \lvert x\rvert$.
	Then integrating \eqref{eq:conv_firstorder} over $x \in \B(a_i,\rho)$,
	\begin{align*}
    &\int_{\B(a_i,\rho)}f_n(|\D u_n|) - 2\pi \int_{0}^\rho f_n\Big(\frac{\lambda([u_*,a_i])}{2\pi r}\Big) r\d r \\
    &\quad\geq \int_0^\rho f_{n,+}'\Big(\frac{\lambda([u_*,a_i])}{2\pi r}\Big) \Big(\int_{\partial\B(a_i,r)}|\D u_n|\d\mathcal H^1 - {\lambda([u_*,a_i])}\Big) \d r\\
    &\quad\geq f_{n,+}'\Big(\frac{\lambda([u_*,a_i])}{2\pi \rho}\Big)\int_0^\rho  \Big(\int_{\partial\B(a_i,r)}|\D u_n| \d\mathcal H^1 - {\lambda([u_*,a_i])}\Big) \d r
	\end{align*}
  since each $f_{n,+}'(t)$ is non-decreasing in $t$.
  Noting that $\lim_{n \to \infty}f_{n,+}'(t) = t$ by Lemma \ref{lem:fn_conv}\eqref{item:derivative_pointwise}, taking the limit in $n$,
  \begin{equation}
    \begin{split}
     &\frac{\lambda([u_{\ast},a_i])}{2\pi \rho} \varlimsup_{n \to \infty}\int_{\B(a_i,\rho)} \lvert \D u_n\rvert \d x \\ 
    &\quad= \mathrm{L}_{i,\rho} + \frac{\lambda([u,a_i])^2}{2 \pi} + \varlimsup_{n \to \infty} \left[ \frac{\lambda([u_{\ast},a_i])^2}{4\pi} \mathcal V(f_n) - 2\pi \int_{0}^\rho f_n\Big(\frac{\lambda([u_*,a_i])}{2\pi r}\Big) r\d r \right] \\
    &\quad= \mathrm{L}_{i,\rho} + \frac{\lambda([u,a_i])^2}{2 \pi} + \frac{\lambda([u_{\ast},a_i])^2}{4\pi} \log\frac{\lambda([u_{\ast},a_i])^2}{4\pi\rho} 
  \end{split}
  \end{equation} 
  using Lemma \ref{lemma:entropy_error} in the last line.
  Then using the estimate \eqref{item:isboundloc},
  \begin{equation}
    \begin{split}
    \mathrm{L}_{i,\rho} +\frac{\lambda([u_{\ast},a_i])^2}{4\pi} \log\frac{\lambda([u_{\ast},a_i])^2}{4\pi\rho} 
    &\leq \mathrm{D} + \mathrm{E}^{1,2}_{\ren}(u_{\ast};\B(a_i,\rho)) + \frac{\lambda([u_{\ast},a_i])^2}{4\pi} \log \frac1{\rho} \\
    &= \mathrm{D} + \int_0^{\rho} \int_{\partial \B(a_i,r)} \frac12\lvert \D u_{\ast}\rvert^2 \d\mathcal H^1 - \frac{\lambda([u_{\ast},a_i])^2}{4\pi r} \d r,
    \end{split}
  \end{equation} 
  where we have set $\mathrm{D} \doteq \mathrm{L} - \Eren^{1,2}(u_{\ast}) - \mathrm{H}([u_{\ast},a_i]).$
  Therefore we have
  \begin{equation}\label{eq:l1bound_proof}
    \begin{split}
    \varliminf_{n \to \infty} \int_{\B(a_i,\rho)} \lvert \D u_n\rvert \d x 
    \leq \frac{2\pi \rho}{\lambda([u_{\ast},a_i])} &\bigg[ \mathrm D + \lambda([u_{\ast},a_i])^2  \\
    &\quad+ \int_0^{\rho}\int_{\partial\B(a_i,r)} \frac12\lvert \D u_{\ast}\rvert^2\d\mathcal H^1 - \frac{\lambda([u_{\ast},a_i])^2}{4\pi r} \d r \bigg] ,
    \end{split}
  \end{equation} 
  establishing \eqref{item:L1bound}.

  To show equi-integrability, for each $i=1,\dotsc,\kappa$ we will set
  \begin{equation}
    \mathrm{M}_i \doteq \frac{2\pi}{\lambda([u_{\ast},a_i])}\bigg[ \mathrm D + \lambda([u_{\ast},a_i])^2
    +\int_0^{\bar\rho_{\Omega}}\int_{\partial\B(a_i,r)} \frac12\lvert \D u_{\ast}\rvert^2\d\mathcal H^1 - \frac{\lambda([u_{\ast},a_i])^2}{4\pi r} \d r \bigg].
  \end{equation} 
  We will assume $\kappa \geq 1$ in what follows, as otherwise the result follows from Proposition \ref{prop:compacitnoprblem}\eqref{item:unconstrained_bound}.
  Then recall from \eqref{item:isboundaway}, using \eqref{eq:erenwithoutlimit} for the renormalised energy, that for $0 < \rho < \bar\rho_{\Omega}$ we have
  \begin{equation}
    \begin{split}
      \varlimsup_{n \to \infty} \int_{\Omega \setminus \bigcup_{i=1}^{\kappa} \B(a_i,\rho)} f_n(\lvert \D u_n\rvert) \d x 
      &\leq \mathrm{D} + \int_{\Omega \setminus \bigcup_{i=1}^{\kappa}\B(a_i,\rho)} \frac12 \lvert \D u_{\ast}\rvert^2 \d x \\
      &= \mathrm{L} - \mathrm{H}([u_{\ast},a_i])_{i=1}^{\kappa} + \Esgnu(g) \log\frac1{\rho}.
    \end{split}
  \end{equation} 
  Combining with \eqref{eq:asympequicont} from Proposition \ref{prop:compacitnoprblem}\eqref{item:unconstrained_bound}
, for each $0 < \rho < \bar\rho_{\Omega}$  and any measurable subset $A \subset \Omega \setminus \bigcup_{i=1}^{\kappa} \B(a_i,\rho)$ we have 
  \begin{equation}\label{eq:l1_asympint}
    \begin{split}
      \varlimsup_{n \to \infty} \int_A \lvert \D u_n\rvert \d x \leq \left( \mathrm{L} - \mathrm{H}([u_{\ast},a_i])_{i=1}^{\kappa} + \Esgnu(g) \log\frac1{\rho} \right)^{\frac12}  \sqrt{2 \Leb^2(A)}.
    \end{split}
  \end{equation} 
  Let $\eps>0$, and choose $\rho >0$ such that $\mathrm{M}_i\rho < \frac{\eps}{2\kappa}$
  for all $i$. 
  For this choice of $\rho$ we then take $\delta>0$ such that
  \begin{equation}
    \left( \mathrm{L} - \mathrm{H}([u_{\ast},a_i])_{i=1}^{\kappa}+ \Esgnu(g) \log\frac1{\rho_0} \right)^{\frac12}  \sqrt{2 \delta} < \frac{\eps}2.
  \end{equation} 
  Then given $A \subset \Omega$ measurable such that $\Leb^2(A) \leq \delta$, using \eqref{eq:l1bound_proof}, \eqref{eq:l1_asympint} we have
  \begin{equation}
    \begin{split}
    \varlimsup_{n \to \infty} \int_A \lvert \D u_n\rvert \,\d x \leq \varlimsup_{n \to \infty} \int_{A \setminus \bigcup_{i=1}^{\kappa} \B(a_i,\rho_0)}\lvert \D u_n\rvert \,\d x + \sum_{i=1}^{\kappa}\varlimsup_{n \to \infty} \int_{\B(a_i,\rho)} \lvert \D u_n\rvert \,\d x < \eps.
  \end{split}
  \end{equation} 
  Thus by Lemma \ref{lemma:asympunifint} we have $(\D u_n)_n$ is uniformly integrable, from which weak $\LL^1$-convergence also follows by the Dunford--Pettis theorem (Lemma \ref{lemma:ourinstanceofDunfordPettis}), establishing \eqref{item:weakconvergence}.

This concludes the proof of Theorem \ref{thm:convsub}.
\end{Proof}

\section{Convergence of almost minimisers}\label{sec:convergence-of-almost-minimisers}

\subsection{Convergence theorem}
We will show the convergence of Theorem \ref{thm:convsub} can be improved for sequences of almost minimisers in the sense of \eqref{eq:uscren} below.
This will also imply Theorem \ref{thm:2ndorder} from the introduction.

\begin{theorem}\label{thm:strongconv}
  Given an approximating family $(f_n)_n$ and $g \in \WW^{\sfrac12,2}(\partial\Omega;\mathcal N)$ define
  \begin{equation}\label{eq:min_problem}
    \mathrm{I}_n = \inf\left\{ \int_{\Omega} f_n(\lvert\D v\rvert)\d x : \begin{matrix} v \in \WW^{1,1}(\Omega;\mathcal N)\\ \Tr_{\partial\Omega}v=g\end{matrix}\right\}.
  \end{equation} 
  Then suppose $(u_n)_n \subset \WW^{1,1}(\Omega;\mathcal N)$ is an asymptotically minimising in the sense that $\Tr_{\partial\Omega}u_n = g$ for all $n$ and
  \begin{equation}\label{eq:uscren}
    \lim_{n \to \infty} \left(\int_{\Omega} f_n(\lvert \D u_n\rvert) \d x - \mathrm{I}_n\right) = 0.
  \end{equation} 
  Then passing to a subsequence, we have $u_n \to u_{\ast}$ strongly in $\LL^1(\Omega)$ to a limit map $u_{\ast} \in \WW^{1,2}_{\ren}(\Omega;\mathcal N)$ satisfying $\Tr_{\partial\Omega}u_{\ast}=g$ with $\kappa \in \mathbb N$ singularities $a_1,\dotsc,a_{\kappa} \in \Omega$, such that the following holds:
  \begin{enumerate}[\rm(i)]
    \item \label{item:approxoferenmin}
        $\displaystyle\lim_{n \to \infty} \left[\int_{\Omega} f_n(\lvert\D u_n\rvert) \d x - \mathcal{V}(f_n) \Esgnu(g)\right] = \mathrm{E}_{\ren}^{1,2}(u_{\ast}) + \mathrm H([u_{\ast},a_i])_{i=1}^{\kappa}.$
    \item \label{item:mineren} The limit map satisfies
      \begin{equation*}
      \mathrm{E}^{1,2}_{\ren}(u_{\ast}) + \mathrm{H}([u_{\ast},a_i])_{i=1}^{\kappa} = \inf\left\{\mathrm{E}^{1,2}_{\ren}(\bar u) + \mathrm{H}([\bar u,\bar a_i])_{i=1}^{\bar\kappa} \right\},
      \end{equation*}
      where the infimum is taken over all $\bar u\in \WW^{1,2}_{\ren}(\Omega;\mathcal N)$ with $\Tr_{\partial\Omega}\bar u = g$, where $\bar u$ has $\bar\kappa\in\mathbb N$ singularities $\bar a_1,\dotsc,\bar a_{\bar\kappa}$, forming a minimal topological resolution in that $\Esgnu(g)=\sum_{i=1}^{\bar\kappa} \frac{\lambda([\bar u,\bar a_i])^2}{4\pi}$.

    \item \label{item:convofmass} $\displaystyle\lim_{n \to \infty} \int_{\Omega} f_n(\lvert \D u_{\ast}\rvert)  - f_n(\lvert \D u_n\rvert) \d x = 0$.

    \item \label{item:L1convergence} $\D u_n \to \D u_{\ast}$ strongly in $\LL^1(\Omega)$ and pointwise almost everywhere.
    \item \label{item:strongconvawayfromsing} For each $\rho < \rho_{\Omega}(a_1,\dotsc,a_{\kappa})$ we have,
      \begin{equation*}
        \lim_{n \to \infty} \int_{\Omega \setminus \bigcup_{i=1}^{\kappa} \B(a_i,\rho)} f_n(\lvert \D u_{\ast} - \D u_n\rvert) \d x = 0
      \end{equation*} 
    \item \label{item:localised_conv} For each $\rho<\rho_{\Omega}(a_1,\dotsc,a_{\kappa})$ and each $i = 1, \dotsc, \kappa$ ,
      \begin{equation*}
        \lim_{n \to \infty} \left[\int_{\B(a_i,\rho)} f_n(\lvert \D u_n\rvert) - f_n\left( \frac{\lambda([u_{\ast},a_i])}{2\pi \lvert x - a_i\rvert} \right) \d x \right] = \int_0^{\rho}\int_{\partial\B(a_i,r)} \frac{\lvert \D u_{\ast}\rvert^2}2 \d\mathcal H^1 - \frac{\lambda([u_{\ast},a_i])^2}{4 \pi r} \d r.
      \end{equation*} 
  \end{enumerate}
\end{theorem}

By \eqref{item:mineren}, if $v \in \WW^{1,2}_{\ren}(\Omega;\mathcal N)$ also has $\kappa$ singularities at $b_1,\dotsc,b_{\kappa} \in \Omega$ such that $[u_{\ast},a_i] = [v,b_i]$ for all $1 \leq i \leq \kappa$ (\emph{i.e.}\,the limiting geodesics are homotopic), then $\Eren^{1,2}(u_{\ast}) \leq \Eren^{1,2}(v)$.
Therefore we have $u_{\ast}$ is a \emph{minimal renormalisable singular harmonic map} in the sense of \cite[Def.~7.8]{MonteilEtAl2022}.
In particular, by a direct comparison argument, for all $\rho>0$ sufficiently small we have $u_{\ast}$  is a minimising harmonic map in $\Omega \setminus \bigcup_{i=1}^{\kappa} \B(a_i,\rho)$ with respect to its own boundary conditions, thereby proving Theorem \ref{thm:2ndorder}\eqref{item:minimisingaway}.

\subsection{Asymptotic uniform convexity}

The classical notion of uniform convexity of a normed space (see \emph{e.g.}\,\cite[Def.~5.2.2]{willem2013functional}) implies that weak convergence, combined the convergence of the norm, implies strong convergence (known as Radon–Riesz property). 
In this section we will prove an asymptotic version, valid along sequences of approximating integrands $(f_n)_n$, by exploiting the uniform convexity of the limiting space $\LL^2$. 
More precisely, we will establish the following.

\begin{proposition}\label{prop:strong_conv}
	Consider a sequence of mappings $u_n \in \WW^{1,1}(\Omega;\mathcal N)$ and a sequence of approximating integrands $(f_n)_n$.
	Assume that $(\D u_n)_n$ converges \weaklystar to $\D u$, and that
	\begin{equation}\label{eq:normconvlemma}
		\int_{\Omega}\frac12{|\D u|^2}\d x= \lim_{n\to\infty}\int_{\Omega}f_n(|\D u_n|)\d x.
	\end{equation}
	Then
	\begin{equation}\label{eq:strongconvlim}
		\lim_{n \to \infty} \int_{\Omega} f_n(\lvert \D u_n - \D u\rvert )\d x = 0.
	\end{equation} 
	In particular, $\D u_n \to \D u$ in $\LL^1$. 
\end{proposition}

Here we will work at the level of gradients for the convenience of the reader, however one could more generally with a sequence $(U_n)_n \subset \LL^1(\Omega)$ in place of $(\D u_n)_n$.
We will adapt the strategy described in \cite{brezis1993convergence}, and as a key step we will establish almost everywhere convergence of the gradients, which we specify as a separate lemma.

\begin{lemma}\label{prop:brezisapproachtopointwise} Under the assumptions of Proposition \ref{prop:strong_conv}, up to a subsequence, $\D u_n$ converges pointwise almost everywhere to $\D u$ in $\Omega$.
\end{lemma}
We rely on the following result (Lemma \ref{lemma:poiutwesbrezis}) concerning sequences in $\mathbb R^d$.

\begin{lemma}\label{lemma:poiutwesbrezis}
	Consider a sequence of Young function $(f_n)_n$ converging to $|z|^2/2$. Suppose $a \in \mathbb R$ and $(b_n)_n \subset \mathbb R^d$ is a sequence satisfying 
	\[
	\frac{1}{2}f_n(|a|) + \frac{1}{2}f_n(|b_n|) - f_n\Big(\Big|\frac{a +b_n}{2}\Big|\Big) \xrightarrow{n\to \infty}0.
	\] Then $\ds\lim_{n\to \infty}b_n = a$.
\end{lemma}
\begin{Proof}{Lemma}{lemma:poiutwesbrezis}
	It suffices to prove that the sequence \((b_n)_n\) is bounded. Indeed if this holds, for any convergent subsequence \(b_{n_k} \to b\), since \(f_n\) converges locally uniformly to \(t^2/2\) by Lemma \ref{lem:fn_conv}, we have
	\[
	0 = \frac14 \lvert a \rvert^2 + \frac14 \lvert b \rvert^2 - \frac12 \left\lvert \frac{a + b}{2} \right\rvert^2 = \frac{1}{8} \lvert b - a \rvert^2,
	\]
	which implies that \(b = a\). 
	Therefore every subsequence has a further subsequence converging to \(a\), from which it follows that \(b_n \to a\).

	To show the boundedness of $(b_n)_n$, we suppose otherwise. Then, passing to a subsequence, we can assume that $\lvert b_n\rvert$ is monotone and divergent.
	Letting
	\begin{equation*}
		g_n(t) = \frac12 f_n(\lvert a\rvert) + \frac12 f_n(t) - f_n\left( \frac{\lvert a\rvert+t}2 \right),
	\end{equation*} 
	since $f_n$ is increasing and convex we have
	\begin{equation*}
		0 \leq g_n(\lvert b_n\rvert) \leq \frac12 f_n(\lvert a\rvert) + \frac12 f_n(\lvert b_n\rvert) - f_n\left( \left\lvert \frac{a+b_n}2\right\rvert \right) 
	\end{equation*} 
	and so
	\(
	\lim_{n \to \infty} g_n(\lvert b_n\rvert) \to 0.
	\) 
	On the other hand we have $g_n$ is non-decreasing by noting that
  \[g_n'(t) = \frac12 \left( f_n'(t) - f_n'\left(\frac{\lvert a\rvert+t}2\right) \right) \geq 0\] 
  for almost all $t \geq \lvert a\rvert$,
  where we used that $t \mapsto f_+'(t)$ is right-continuous and non-decreasing, and that $g$ is absolutely continuous since $f$ is (see \emph{e.g.}\,\cite[Thm.~3.20, Ex.~3.21]{Leoni2017}). 
  Now fixing $m$ sufficiently large such that $\lvert b_m\rvert > \lvert a\rvert$, which exists by unboundedness of $(b_n)_n$, we have
	\begin{equation*}
		0 = \lim_{n \to \infty} g_n(\lvert b_n\rvert) \geq \lim_{n \to \infty} g_n(\lvert b_m\rvert) = \frac14 \lvert a\rvert^2 + \frac14 \lvert b_m\rvert^2 - \frac12\left( \frac{\lvert a\rvert+\lvert b_m\rvert}2 \right)^2 \geq 0.
	\end{equation*} 
  Since $\lvert a\rvert \neq \lvert b_m\rvert$, we have the right-hand side is strictly positive, thereby contradicting the unboundedness of $(b_n)_n$.
\end{Proof}

\begin{Proof}{Lemma}{prop:brezisapproachtopointwise}
  Applying Proposition \ref{prop:compacitnoprblem} to $(\D u_n + \D u)/2$, we have \eqref{eq:fatoucompactplb} reads as
  \begin{equation}\label{eq:half_lsc}
	\int_{\Omega}\frac12{|\D u|^2} \d x \leq \varliminf_{n\to \infty}\int_\Omega f_n\Big(\Big|\frac{\D u_n + \D u}{2}\Big|\Big) \d x,
\end{equation}
and by the dominated convergence theorem and Lemma \ref{lemma:quadgrowth} we also have
\begin{equation}\label{eq:frozen_conv}
  \lim_{n \to \infty} \int_{\Omega} f_n(\lvert \D u\rvert) \d x = \int_{\Omega} \frac12 \lvert \D u\rvert^2 \d x.
\end{equation} 
	For each $n$, on $\Omega$ we define
	\[
	g_n \doteq \frac{1}{2}f_n(|\D u|) + \frac{1}{2}f_n(|\D u_n|) - f_n\Big(\Big|\frac{\D u_n + \D u}{2}\Big|\Big).
	\]
	By convexity of $f_n$ we have $g_n \geq 0$, and so we have 
	\begin{align*}
		\varlimsup_{n \to \infty}\int_\Omega g_n \d x &\leq \frac{1}{2}	\varlimsup_{n \to \infty}\int_\Omega f_n(|\D u|)\d x  + \frac{1}{2}	\varlimsup_{n \to \infty}\int_\Omega f_n(|\D u_n|) \d x- \varliminf_{n\to \infty}\int_\Omega f_n\Big(\Big|\frac{\D u_n + \D u}{2}\Big|\Big) \d x \\
		&\leq \left(\frac{1}{2} + \frac12 - 1\right)	\int_{\Omega}\frac12{|\D u|^2}\d x = 0
	\end{align*}
  where we used \eqref{eq:normconvlemma}, \eqref{eq:half_lsc} and \eqref{eq:frozen_conv}. We deduce that $g_n \to 0$ in $\LL^1(\Omega)$, 
  and hence (using \emph{e.g.}\,\cite[Prop.~4.2.10]{willem2013functional}), up to a subsequence we have $g_n \to 0$ almost everywhere. Then by Lemma \ref{lemma:poiutwesbrezis}, we infer that $\D u_n \to \D u$ almost everywhere.
\end{Proof}
\begin{Proof}{Proposition}{prop:strong_conv}
	By Proposition \ref{prop:brezisapproachtopointwise}, passing to a subsequence we have $\D u_n \to \D u$ almost everywhere.
	Then by Scheff\'e's lemma (asserting that pointwise convergence and convergence of the norm implies strong $\LL^1$ convergence, see \textit{e.g.}\,\cite[Prop.~2.4.6]{willem2013functional}) we have
	\begin{equation}
		\lim_{n\to\infty}\int_{\Omega} \Big\lvert f_n(\lvert \D u_n\rvert) - \frac12 \lvert \D u\rvert^2 \Big\rvert \d x = 0. 
	\end{equation} 
	In particular, $(f_n(\lvert \D u_n\rvert))_n$ is equi-integrable.
	By \ref{hyp:ndec2} we have the almost-doubling estimate
	\begin{equation}
		f_n(2t) \leq f_n(2t_0) + 4 f(t)
	\end{equation} 
	for all $n$ and $t \geq 0$, so we obtain the pointwise estimate
	\begin{equation}
		\lvert f_n(\lvert \D u - \D u_n\rvert)\rvert \leq 2f_n(2t_0) + 4 \left( f_n(\lvert \D u\rvert) + f_n(\lvert \D u_n\rvert) \right).
	\end{equation} 
	Since $f_n(\lvert \D u\rvert) \to \frac12 \lvert \D u\rvert^2$ strongly in $\LL^1(\Omega)$ and $f_n(2t_0) \to 2t_0^2$ as $n \to \infty$,
  we infer the left hand side is equi-integrable. Therefore \eqref{eq:strongconvlim} follows by Vitali's convergence theorem (see \emph{e.g.}\,\cite[Thm.~3.1.12]{willem2013functional}), from which convergence in $\LL^1$ follows by Jensen's inequality.
\end{Proof}

\subsection[Strong convergence of almost minimisers]{Strong convergence of almost minimisers --- Proof of Theorem \ref{thm:strongconv}}

\begin{Proof}{Theorem}{thm:strongconv}
  Let $\bar u \in \WW^{1,2}_{\ren}(\Omega;\mathcal N)$ be such that $\Tr_{\partial\Omega}\bar u = g$, which always exists by \cite[Prop.~2.4, Prop.~8.1]{MonteilEtAl2022}, and suppose it has $\bar{\kappa}$ singularities at $\bar a_1,\dotsc, \bar a_{\bar\kappa}$ such that $\Esgnu(g) = \sum_{i=1}^{\bar\kappa} \frac{\lambda([\bar u,\bar a_i])^2}{4\pi}$.
  Then by the upper bound (Proposition \ref{prop:upperBound}) we have
  \begin{equation}
    \begin{split}
      \varlimsup_{n \to \infty} \left(\mathrm{I}_n - \mathcal V(f_n) \Esgnu(g)\right) 
      &\leq \left[\lim_{n \to \infty} \int_{\Omega} f_n(\lvert \D \bar u\rvert) - \mathcal V(f_n) \Esgnu(g)\right]  \\
      &= \mathrm{E}_{\ren}^{1,2}(\bar u) + \mathrm{H}([\bar u,\bar a_i])_{i=1}^{\bar\kappa},
    \end{split}
  \end{equation} 
  noting $\bar u$ is admissible for the minimisation problem \eqref{eq:min_problem} for each $n$.
  Therefore if $u_n \in \WW^{1,1}_g(\Omega;\mathcal N)$ satisfies \eqref{eq:uscren}, it follows that
  \begin{equation}\label{eq:almostmin_ineq}
    \varlimsup_{n \to \infty} \left[ \int_{\Omega} f_n(\lvert \D u_n\rvert) \,\d x - \mathcal V(f_n)\Esgnu(g) \right] \leq \mathrm{E}_{\ren}^{1,2}(\bar u) + \mathrm{H}([\bar u, \bar a_i])_{i=1}^{\bar\kappa}.
  \end{equation} 
  Therefore by Theorem \ref{thm:convsub}, noting \eqref{eq:asymptotic_finiteness} is satisfied, by passing to a subsequence we infer that $u_n$ converges to a limiting map $u_{\ast} \in \WW^{1,2}_{\ren}(\Omega;\mathcal N)$ with $\kappa$ singularities $a_1,\dotsc,a_{\kappa}$, which by Theorem \ref{thm:convsub}\eqref{item:limit_renormalisable} satisfies
  \begin{equation}\label{eq:lsc_renormalisable}
    \varliminf_{n \to \infty} \left[\int_{\Omega} f_n(\lvert \D u_n\rvert) - \mathcal V(f_n) \Esgnu(g)\right] \geq \mathrm{E}^{1,2}_{\ren}(u_{\ast}) + \mathrm{H}([u_{\ast},a_i])_{i=1}^{\kappa},
  \end{equation} 
  and hence that
  \begin{equation}
    \mathrm{E}^{1,2}_{\ren}(u_{\ast}) + \mathrm{H}([u_{\ast},a_i])_{i=1}^{\kappa} \leq \mathrm{E}^{1,2}_{\ren}(\bar u) + \mathrm{H}([\bar u,\bar a_i])_{i=1}^{\bar\kappa}.
  \end{equation} 
  While we selected a particular competitor, the same argument goes through for any
  $\bar u \in \WW^{1,2}_{\ren}(\Omega;\mathcal N)$ satisfying $\Tr_{\partial\Omega}\bar u=g$, thereby establishing \eqref{item:mineren}.
  Moreover, if we take $\bar u = u_{\ast}$ in \eqref{eq:almostmin_ineq}, combining with \eqref{eq:lsc_renormalisable} we deduce \eqref{item:approxoferenmin}.
  For \eqref{item:convofmass} we combine \eqref{item:approxoferenmin} and Proposition \ref{prop:upperBound} applied with $u_{\ast}$ which gives
  \begin{equation}
    \begin{split}
   \lim_{n \to \infty} \left[\int_{\Omega} f(\lvert \D u_{\ast}\rvert) - f_n(\lvert \D u_n\rvert) \d x\right] 
    &=\lim_{n \to \infty} \left[\int_{\Omega} f(\lvert \D u_{\ast}\rvert) \d x - \mathcal V(f_n)\Esgnu(g)\right]\\
    &\quad- \lim_{n \to \infty} \left[\int_{\Omega} f(\lvert \D u_n\rvert) \d x - \mathcal V(f_n)\Esgnu(g)\right] 
  \end{split} 
  \end{equation}
  which is zero, as required.

  Now, combining \eqref{item:approxoferenmin} and Theorem \ref{thm:convsub}\eqref{item:isboundaway}, for all $0<\rho<\rho_{\Omega}(a_1,\dotsc,a_{\kappa})$ we have
  \begin{equation}
  \lim_{n \to \infty} \int_{\Omega \setminus \bigcup_{i=1}^{\kappa} \B(a_i,\rho)} f_n(\lvert \D u_n\rvert) \d x = \int_{\Omega\setminus\bigcup_{i=1}^{\kappa} \B(a_i,\rho)} \frac12 \lvert \D u\rvert^2 \d x.
  \end{equation} 
  Hence applying Proposition \ref{prop:strong_conv} in $\Omega \setminus \bigcup_{i=1}^{\kappa} \B(a_i,\rho)$ we infer both \eqref{item:strongconvawayfromsing} and that $\D u_n \to \D u$ strongly in $\LL^1(\Omega \setminus \bigcup_{i=1}^{\kappa} \B(a_i,\rho))$.
  In particular, by applying this to $\rho_k \searrow 0$ and passing to a diagonal subsequence, we deduce that $\D u_n \to \D u$ almost everywhere.
  Since $(\D u_n)_n$ is equi-integrable in $\Omega$ by Theorem \ref{thm:convsub}\eqref{item:weakconvergence}, strong $\LL^1$ convergence in $\Omega$ follows, thereby establishing \eqref{item:L1convergence}.

  Finally combining \eqref{item:approxoferenmin} and Theorem \ref{thm:convsub}\eqref{item:isboundloc}
  we have for $0 < \rho < \rho_{\Omega}(a_1,\dotsc,a_{\kappa})$ and  $i = 1,\dotsc,\kappa$,
  \begin{equation}
    \begin{split}
   &\lim_{n \to \infty} \left[\int_{\B(a_i,\rho)} f_n(\lvert \D u_n\rvert) \d x - \frac{\lambda([u_{\ast},a_i])^2}{4\pi} \mathcal V(f_n) \right] \\
   &\qquad= \Eren^{1,2}(u_{\ast};\B(a_i,\rho)) + \frac{\lambda([u_{\ast},a_i])^2}{4\pi} \log \frac{2\pi}{\lambda([u_{\ast},a_i])},
  \end{split}
  \end{equation} 
  which together with Lemma \ref{lemma:entropy_error} gives \eqref{item:localised_conv}.
\end{Proof}

\begin{Rmk}\label{rmk:firstorder_proof}
  We will sketch how the above proof can be adapted to prove Theorem \ref{thm:firstorder} in the absence of the \ref{hyp:cdec2} condition.
  Let $(\bar f_n)_n$ be the modified family from Example \ref{eq:weak2strong_dec2}, and let $\bar{\mathrm I}_n$ as in \eqref{eq:min_problem} with $\bar f_n$ in place of $f_n$.
  Then for any $\bar u \in \WW^{1,2}_{\ren}(\Omega;\mathcal N)$ such that $\Tr_{\partial\Omega}\bar u = g$ whose $\bar\kappa \in \mathbb N$ singularities $\bar a_1,\cdots,\bar a_{\bar\kappa}$ forms a minimal topological resolution of $g$ in $\Omega$, by Proposition \ref{prop:upperBound} we have
  \begin{equation}
    \lim_{n \to \infty} \left[ \int_{\Omega} \bar f_n(\lvert \D \bar u\rvert) \,\d x - \mathcal V(\bar f_n) \Esgnu(g) \right] = \Eren^{1,2}(\bar u) +\mathrm{H}([\bar u,\bar a_i])_{i=1}^{\bar\kappa}.
  \end{equation} 
  Therefore since $\mathrm{I}_n \leq \bar{\mathrm I}_n$ for all $n$, using \eqref{eq:vint_diff} we have
  \begin{equation}\label{eq:thm1_upper}
    \begin{split}
    \varlimsup_{n \to \infty} \left[ \mathrm{I}_n - \mathcal V(f_n) \Esgnu(g)\right] 
    &\leq \varlimsup_{n \to \infty} \left[\bar{\mathrm{I}}_n  - \mathcal V(\bar f_n)\Esgnu(g)\right] + \frac12 \Esgnu(g) \lim_{n \to \infty} \bar f_n(1)  \\
    &\leq \Eren^{1,2}(u_{\ast}) + \mathrm{H}([u_{\ast},a_i])_{i=1}^{\kappa} + \frac14 \Esgnu(g),
    \end{split}
  \end{equation} 
  which in particular is finite.
  Now taking a sequence $(u_n)_n$ satisfying \eqref{eq:uscren}, we can apply Theorem \ref{thm:convsub} with the family $(f_n)_n$, where it is sufficient to assume \ref{hyp:ndec2}.
  In particular using Theorem \ref{thm:convsub}\eqref{item:limit_renormalisable} we obtain a limit map $u_{\ast} \in \WW^{1,2}_{\ren}(\Omega;\mathcal N)$ such that $\Tr_{\partial\Omega}u_{\ast} = g$ with $\kappa \in \mathbb N$ singularities $a_1,\cdots,a_{\kappa} \in \Omega$ such that
  \begin{equation}\label{eq:thm1_lower}
    \varliminf_{n \to \infty} \left[ \mathrm{I}_n - \mathcal V(f_n) \Esgnu(g)\right] \geq \Eren^{1,2}(u_{\ast}) + \mathrm{H}([u_{\ast},a_i])_{i=1}^{\kappa}.
  \end{equation} 
  Hence combining \eqref{eq:thm1_upper} and \eqref{eq:thm1_lower} we deduce that
  \begin{equation}
    \lim_{n \to \infty} \frac1{\mathcal V(f_n)}\, \mathrm{I}_n = \Esgnu(g),
  \end{equation} 
  thereby proving Theorem \ref{thm:firstorder}.

  We point out that this argument does not allow us to relax the \ref{hyp:cdec2} assumption in Theorem \ref{thm:2ndorder}; for this \eqref{eq:vint_diff} is insufficient, we would instead need an approximating sequence $(\tilde f_n)_n$ satisfying $f_n \leq \tilde f_n$ for all $n$ and such that $\displaystyle\lim_{n \to \infty} \big[ \mathcal V(\tilde f_n) - \mathcal V(f_n)\big] = 0$.
\end{Rmk}

\section{Universality of the vortex map}\label{sec:universality-of-vortex-map}

In this section we will show, under strong conditions on the boundary datum and domain, that the full sequence of minimisers $(u_n)_n$ converges to a unique minimiser of vortex type. 
Since the obtained limit map is independent of the choice of approximating integrands $(f_n)_n$, this shows mappings of type $g(\frac{x}{\lvert x\rvert})$ are in some sense \emph{universal}.
An analogous result for the Ginzburg-Landau approximation was obtained in \cite[Thm.~8.1]{MonteilEtAl2021}, and  we will employ a similar strategy of proof.

We use the synharmonic distance between two maps $g,g' \in \WW^{\sfrac{1}{2},2}(\mathbb S^1; \mathcal N)$, which is defined as
\begin{equation}\label{eq:synh_dist}
\dist_{\mathrm{synh}}(g,g') \doteq \inf \, \biggl\{ \int_{\mathbb S^1} \int_0^T \biggl(\frac{|{\D u}|^2}{2} - \frac{\lambda([g])^2}{8 \pi^2} \biggr)\d t\d\mathcal H^1 \biggr\},
\end{equation} 
where the infimum is taken over all $T>0$ and $u \in \WW^{1,2}(\mathbb S^1 \times [0,T];\mathcal N)$ satisfying $\tr_{\mathbb S^1 \times \{0\}} = g,$ $\tr_{\mathbb S^1 \times \{T\}} = g'$ on $\mathbb S^1$.
This was introduced in \cite[Def.~3.2]{MonteilEtAl2022}, where we refer the reader for further properties and examples. Note this infimum is empty if $g$ and $g'$ are not homotopic, in which case we have $\dist_{\mathrm{synh}}(g,g')=\infty$.

\begin{theorem}\label{thm:universalityofthevortexmap}
  Let $g \in \CC^1(\mathbb S^1; \mathcal N)$ be a minimising geodesic satisfying the following:
	\begin{enumerate}[(a)]
		\item \label{item:hypmorehtanatom} If $(\gamma_1,\dots,\gamma_k)$ is  a minimal topological resolution of $g$, then $k = 1$ and $\gamma_1$ is homotopic to $g$.
		\item \label{item:hypsynharmony} If $g' \in \WW^{\sfrac12,2}(\mathbb S^1; \mathcal N)$ is a minimising geodesic in the same homotopy class as $g$, then \(\dist_{\mathrm{synh}}(g,g')  = 0\).
	\end{enumerate}
	Let $(f_n)_n$ be a sequence of approximating integrands,
	and $(u_n)_n \subset \WW^{1,1}(\B; \mathcal N)$ be an asymptotically minimising sequence in that $\tr_{\mathbb S^1}u_n = g$ for all $n$ and
	\begin{equation}\label{eq:uscrenII}
		\int_{\B_1} f_n(\lvert \D u_n\rvert) \d x - \inf\left\{ \int_{\B_1} f_n(\lvert\D v\rvert)\d x : \begin{matrix} v \in \WW^{1,1}(\B_1;\mathcal N)\\ \Tr_{\partial\Omega}v=g\end{matrix}\right\} \xrightarrow{n \to \infty} 0.
	\end{equation} 
  Then we have $u_n \to u_{\ast}$ in $\WW^{1,1}(\B;\mathcal N)$ as $n \to \infty$, where $u_*(x) = g(\frac{x}{|x|})$.
\end{theorem}
Recall that $\B$ denotes the closed unit ball in $\mathbb R^2$. Assumption \eqref{item:hypmorehtanatom} implies the atomicity (defined above in Definition \ref{def:singen}) of the geodesic \(g\), however it is in general a stronger condition, as \cite[\S 9.3.5]{MonteilEtAl2022} illustrates. 
 
In the case $\mathcal N = \mathbb S^1$, up to reparametrisation there are two minimising geodesics that are atomic; namely $g(z) = z$ and $\bar{g}(z) = \bar z$, where we identify $z \in \mathbb S^1 \subset \mathbb C$.
Since $g$ and $\bar{g}$ lie in different homotopy classes, and since any reparametrisation of $g$ or $\bar{g}$ only differ by a rotation, by \cite[Prop.~3.8]{MonteilEtAl2022} both mappings satisfy \eqref{item:hypmorehtanatom} and \eqref{item:hypsynharmony}.
Therefore Theorem \ref{thm:universalityofthevortexmap} applies, and in particular we infer the universality of the vortex map $u_{\mathrm{V}}(x) = \frac{x}{\lvert x\rvert}$, proving Theorem \ref{thm:univrotex}.
In the specific case of \(f_n(z) = \sfrac{|z|^{p_n}}{p_n}\), where \(p_n \nearrow 2\), this behaviour is expected because the vortex map is the unique \(p\)-minimiser for every \(p \in (1, 2)\); see \cite{hardt1998penergy} and \cite[Thm.~13.6]{brezis2021sobolev}.
 
Theorem \ref{thm:universalityofthevortexmap} will largely follow from properties of the geometric energy introduced in \cite[(2.11)]{MonteilEtAl2022}, which we recall is defined as follows: given a bounded Lipschitz domain \(\Omega \subset \mathbb{R}^2\), an integer \(k \in \mathbb{N}\setminus\{0\}\), distinct points \(a_1, \dotsc, a_k \in \Omega\), a radius \(\rho \in (0, \Bar{\rho}(a_1, \dotsc, a_k))\), a map \(g \in \WW^{1/2,2}(\partial \Omega, \mathcal{N})\) and a topological resolution \(\gamma_1,  \dotsc, \gamma_k\in \WW^{1/2,2}(\mathbb{S}^1, \mathcal{N})\) of \(g\), the \emph{geometrical renormalised energy} is given by:
 \begin{multline}
 	\label{eq:def_renorm_geom}
\mathcal{E}^{\mathrm{geom}}_{ g, \gamma_1, \dotsc, \gamma_k} (a_1, \dotsc, a_k)
  \doteq  \lim_{\rho \searrow 0}\inf \biggl\{ \int_{\Omega \setminus \bigcup_{i = 1}^k \B_\rho (a_i)} \frac{|\D u|^2}{2} \d x - \sum_{i=1}^{k}\frac{\lambda(\gamma_i)^2}{4 \pi^2} \log \frac1{\rho} \\
    : \begin{matrix} u \in \WW^{1, 2} \big(\Omega \setminus \bigcup_{i = 1}^k \B_\rho (a_i); \mathcal{N}\big),\tr_{\partial \Omega} u = g\\
    \tr_{\mathbb{S}^1} u (a_i + \rho\,\cdot) = \gamma_i \ \ \text{for all } i = 1,\dotsc,k \end{matrix}
 	\biggr\}.
 \end{multline}
 We will also need two auxiliary results; the first will be a consequence of Theorem \ref{thm:strongconv}\eqref{item:mineren}.

 \begin{lemma}\label{coro:geoandtopinitemi}
   Under the assumptions of Theorem \ref{thm:strongconv}, the limit map $u_{\ast}$ satisfies
 	\begin{equation}
 		\Eren^{1,2}(u_*) + \mathrm{H}([u_*,a_i])_{i = 1}^{\kappa}=
    \inf\{\mathcal{E}^{\mathrm{geom}}_{g,\bar\gamma_1,\dots,\bar\gamma_{\bar\kappa}}(\bar a_1,\dots,\bar a_{\bar\kappa}) + \mathrm{H}([\bar\gamma_i])_{i = 1}^{\bar\kappa} \} 
 	\end{equation}
  where the infimum is taken over all $\bar\kappa\in\mathbb N$, $\bar a_1,\dotsc,\bar a_{\bar\kappa} \in \Omega$ distinct and minimal topological resolutions $\bar\gamma_1,\dotsc,\bar\gamma_{\bar\kappa} \in \CC^1(\mathbb S^1;\mathcal N)$ of $g$.
  Moreover, the infimum is attained in that there exists $\gamma_1,\dotsc,\gamma_{\kappa} \in \CC^1(\mathbb S^1;\mathcal N)$ such that
  \begin{equation}\label{eq:eren_egeom}
    \Eren^{1,2}(u_*) = \mathcal{E}^{\mathrm{geom}}_{g,\gamma_1,\dots,\gamma_{\kappa}}(a_1,\dots,a_{\kappa}).
  \end{equation} 
\end{lemma}

 \begin{Proof}{Lemma}{coro:geoandtopinitemi}
   Given any $\bar u \in \WW^{1,2}_{\mathrm{ren}}(\Omega;\mathcal N)$ with $\bar\kappa$ singularities at $\bar a_1,\dotsc,\bar a_{\bar\kappa}$ with singularities $\bar\gamma_i$ at each $\bar a_i$,
   by \cite[Prop.~7.2(v)]{MonteilEtAl2022} we have $\Eren^{1,2}(\bar u) \geq \mathcal{E}^{\mathrm{geom}}_{g,\bar\gamma_1,\dotsc,\bar\gamma_{\bar\kappa}}(\bar a_1,\dots,\bar a_{\bar\kappa})$.
   Thus combining with Theorem \ref{thm:strongconv}\eqref{item:mineren} we obtain the lower bound.

   Conversely let $u_{\ast} \in \WW^{1,2}_{\ren}(\Omega;\mathcal N)$ be a limiting map from Theorem \ref{thm:strongconv}, and suppose $u_{\ast}$ has $\kappa \in \mathbb N$ singularities at $a_1,\dotsc,a_{\kappa}$ with limiting geodesics $\gamma_1,\dotsc,\gamma_{\kappa}$, in that they can be obtained as the limit of suitable subsequences of $\Tr_{\partial \B(a_i,\rho_{\ell})}u$ as $\rho_\ell \searrow 0$ (see \emph{e.g.}\,\cite[Prop.~7.2(iii)]{MonteilEtAl2022}).
   Then for each $0 < \rho < \rho_{\Omega}(a_1,\dotsc,a_{\kappa})$, there exists $u_{\rho} \in \WW^{1,2}(\Omega \setminus \bigcup_{i=1}^{\kappa}\B(a_i,\rho);\mathcal N)$ minimising $\int_{\Omega \setminus \bigcup_{i=1}^{\kappa}\B(a_i,\rho)} \frac12 \lvert \D u\rvert^2 \d x$ subject to $\Tr_{\partial\Omega}u_{\rho} = g$ and $\Tr_{\mathbb S^1}u_{\rho}(a_i+\rho\,\cdot) = \gamma_i$ for each $i=1,\dotsc,\kappa$. We then extend $u_{\rho}$ to $\Omega$ by setting $u_{\rho} = \gamma_i(\frac{\cdot-a_i}{\rho})$ in each $\B(a_i,\rho)$, which defines a map $u_{\rho} \in \WW^{1,2}_{\ren}(\Omega;\mathcal N)$ such that
   \begin{equation}\label{eq:urho_eren}
     \Eren^{1,2}(u_{\rho}) = \int_{\Omega \setminus \bigcup_{i=1}^{\kappa} \B(a_i,\rho)} \frac12 \lvert \D u_{\rho}\rvert^2 - \Esgnu(g) \log\frac1{\rho},
   \end{equation} 
   noting the contributions on each $\B(a_i,\rho)$ vanish, since each $\gamma_i$ is a minimising geodesic. 
   Hence using Theorem \ref{thm:strongconv}\eqref{item:mineren} to compare $u_{\ast}$ with $u_{\rho}$ sending $\rho \searrow 0$, using \eqref{eq:def_renorm_geom}, \eqref{eq:urho_eren}, we have
   \begin{equation}
     \begin{split}
     \Eren^{1,2}(u_{\ast})+ \mathrm{H}([u_{\ast},a_i])_{i=1}^{\kappa}
     &\leq \lim_{\rho \to 0} \Eren^{1,2}(u_{\rho}) +\mathrm{H}([u_{\ast},a_i])_{i=1}^{\kappa}\\
     &= \mathcal{E}^{\mathrm{geom}}_{g,\gamma_1,\dotsc,\gamma_{\kappa}}(a_1,\dotsc,a_{\kappa})+\mathrm{H}([u_{\ast},a_i])_{i=1}^{\kappa}.
   \end{split}
   \end{equation} 
   This establishes the reverse inequality, and that equality is attained with $(\gamma_1,\dotsc,\gamma_{\kappa})$.
 \end{Proof}

\begin{lemma}\label{lemma:blackboxlemma} Given the unit ball $\B \subset \mathbb R^2$, let \(g \in \CC^{1} (\partial \B, \mathcal{N})\) be a minimising geodesic.	
  \begin{enumerate}[\rm(i)] 
		\item \label{item:synhmaronicgeo} If \(\gamma\) is a minimising geodesic homotopic to $g$ such that 
		\(
    \dist_{\mathrm{synh}}(g,\gamma) = 0,
		\) then 
		\(
  \mathcal{E}^{\mathrm{geom}}_{g,g}  (a)
    = \mathcal{E}^{\mathrm{geom}}_{g,\gamma}(a)
		\) for all $a \in \B^{\circ}$.

		\item \label{item:inuqesolblackbox}
      We have $\mathcal{E}^{\mathrm{geom}}_{g,g}(0) \leq \mathcal{E}^{\mathrm{geom}}_{g,g}(a)$  for all $a \in \B^{\circ}$, with equality if and only if $a=0$.
		
    \item\label{item:findvortexmap} If $v \in \WW^{1,2}_{\ren}(\B;\mathcal N)$ satisfies $\Tr_{\partial \B}v=g$ and $\Erennu(v, \B) = 0$, then $v(x) =g(\frac{x}{|x|})$.
	\end{enumerate}
\end{lemma}
\begin{Proof}{Lemma}{lemma:blackboxlemma}
We have \eqref{item:synhmaronicgeo} follows from \cite[Prop.~3.10]{MonteilEtAl2022}, and \eqref{item:inuqesolblackbox} follows from \cite[Thm.~10.1]{MonteilEtAl2022} (combined with \cite[\S VIII.4]{bethuel2017ginzburg}).
  Item \eqref{item:findvortexmap} follows from the proof of \cite[Thm.~8.1]{MonteilEtAl2021}; here one observes that $\Eren^{1,2}(v) = 0$ implies that $\partial_ru \equiv 0$ on $\B$, from which the claimed form follows.
\end{Proof}

\begin{Proof}{Theorem}{thm:universalityofthevortexmap}
  By Theorem \ref{thm:strongconv} we know $(u_n)_n$ admits a subsequence converging to some $u_{\ast} \in \WW^{1,2}_\ren(\B,\mathcal N)$, so it suffices to show this limit map is uniquely determined as $g(\frac{x}{\lvert x\rvert})$.
  By assumption \eqref{item:hypmorehtanatom}, we have $u_{\ast}$ has only singularity at some $a \in \B^{\circ}$, and let $\gamma \in \CC^1(\mathbb S^1;\mathcal N)$ be homotopic to $g$ such that $\Eren^{1,2}(u_{\ast}) = \mathcal{E}^{\mathrm{geom}}_{g,\gamma}(a)$, which exists by Lemma \ref{coro:geoandtopinitemi}.

  Now by assumption \eqref{item:hypsynharmony} we have $\dist_{\mathrm{synh}}(g,\gamma)=0$ so by Lemma \ref{lemma:blackboxlemma}\eqref{item:synhmaronicgeo} we have $\mathcal{E}^{\mathrm{geom}}_{g,g}(x) = \mathcal{E}^{\mathrm{geom}}_{g,\gamma}(x)$ for all $x \in \B^{\circ}$, and since Lemma \ref{coro:geoandtopinitemi} asserts this is minimised at $x=a$, by Lemma \ref{lemma:blackboxlemma}\eqref{item:inuqesolblackbox} we have $a=0$. 
  Using $g(\frac{x}{\lvert x\rvert})$ as a competitor, we see that $\mathcal{E}^{\mathrm{geom}}_{g,g}(0)= 0$, so by Lemma \ref{lemma:blackboxlemma}\eqref{item:findvortexmap} we deduce that $u_{\ast}(x) = g(\frac{x}{\lvert x\rvert})$.
\end{Proof}

\hypersetup{bookmarksdepth=-1}
\subsection*{Acknowledgements}
The authors would like to thank the anonymous referees for their careful reading and suggestions, which have improved this manuscript.

\subsection*{Funding}
BVV is supported by a FRIA fellowship; cette publication bénéficie du soutien de la Communauté française de Belgique dans le cadre du financement d’une bourse FRIA.

\hypersetup{bookmarksdepth}
\bibliographystyle{abbrv} 
\bibliography{ref_pub.bib}

\begin{thebibliography}{10}

\bibitem{AmbrosioDalMaso1992}
L.~Ambrosio and G.~Dal~Maso.
\newblock On the relaxation in {$\mathrm{BV}(\Omega;\mathrm{R}^m)$} of quasi-convex integrals.
\newblock {\em J. Func. Anal.}, 109(1):76--97, Oct. 1992.

\bibitem{attouch2014variational}
H.~Attouch, G.~Buttazzo, and G.~Michaille.
\newblock {\em Variational analysis in {S}obolev and {BV} spaces}, volume~17 of {\em MOS-SIAM Series on Optimization}.
\newblock Society for Industrial and Applied Mathematics (SIAM), Philadelphia, PA; Mathematical Optimization Society, Philadelphia, PA, second edition, 2014.
\newblock Applications to PDEs and optimization.

\bibitem{beaufort2017computing}
P.-A. Beaufort, J.~Lambrechts, F.~Henrotte, C.~Geuzaine, and J.-F. Remacle.
\newblock Computing cross fields : A {PDE} approach based on the {G}inzburg-{L}andau theory.
\newblock {\em Procedia Engineering}, 203:219--231, 2017.
\newblock 26th International Meshing Roundtable, IMR26, 18-21 September 2017, Barcelona, Spain.

\bibitem{bethuel1991approximation}
F.~Bethuel.
\newblock The approximation problem for {S}obolev maps between two manifolds.
\newblock {\em Acta Math.}, 167:153--206, 1991.

\bibitem{bethuel2017ginzburg}
F.~Bethuel, H.~Brezis, and F.~H\'elein.
\newblock {\em Ginzburg-{L}andau vortices}.
\newblock Modern Birkh\"auser Classics. Birkh\"auser/Springer, Cham, 1994.

\bibitem{bethuel1995extensions}
F.~Bethuel and F.~Demengel.
\newblock Extensions for {S}obolev mappings between manifolds.
\newblock {\em Calc. Var. Partial Differential Equations}, 3(4):475--491, 1995.

\bibitem{Bildhauer2003}
M.~Bildhauer.
\newblock {\em Convex {{Variational Problems}}}, volume 1818 of {\em Lecture {{Notes}} in {{Mathematics}}}.
\newblock Springer Berlin Heidelberg, Berlin, Heidelberg, 2003.

\bibitem{bourgain2001anther}
J.~Bourgain, H.~Brezis, and P.~Mironescu.
\newblock Another look at {S}obolev spaces.
\newblock In {\em Optimal control and partial differential equations}, pages 439--455. IOS, Amsterdam, 2001.

\bibitem{brezis1993convergence}
H.~Brezis.
\newblock Convergence in {${\mathscr D}'$} and in {$L^1$} under strict convexity.
\newblock In {\em Boundary value problems for partial differential equations and applications}, volume~29 of {\em RMA Res. Notes Appl. Math.}, pages 43--52. Masson, Paris, 1993.

\bibitem{Brezis2011}
H.~Brezis.
\newblock {\em Functional Analysis, {{Sobolev}} Spaces and Partial Differential Equations}.
\newblock Springer New York, New York, NY, 2011.

\bibitem{brezis2021sobolev}
H.~Brezis and P.~Mironescu.
\newblock {\em Sobolev maps to the circle---from the perspective of analysis, geometry, and topology}, volume~96 of {\em Progress in Nonlinear Differential Equations and their Applications}.
\newblock Birkh\"{a}user/Springer, New York, 2021.

\bibitem{BrezisNirenberg1995}
H.~Brezis and L.~Nirenberg.
\newblock Degree theory and {{BMO}}; part {{I}}: {{Compact}} manifolds without boundaries.
\newblock {\em Sel. Math.}, 1(2):197--263, Sept. 1995.

\bibitem{bulanyi2023singular}
B.~Bulanyi and J.~{Van Schaftingen}.
\newblock Singular extension of critical {S}obolev mappings under an exponential weak-type estimate.
\newblock {\em J. Func. Anal.}, 288(1), Jan. 2025.

\bibitem{canevari2028defects}
G.~Canevari and A.~Segatti.
\newblock Defects in nematic shells: a {$\Gamma$}-convergence discrete-to-continuum approach.
\newblock {\em Arch. Ration. Mech. Anal.}, 229(1):125--186, 2018.

\bibitem{EvansGariepy2015}
L.~C. Evans and R.~F. Gariepy.
\newblock {\em Measure Theory and Fine Properties of Functions}.
\newblock Chapman and Hall/CRC, 2015.

\bibitem{FonsecaLeoni2007}
I.~Fonseca and G.~Leoni.
\newblock {\em Modern Methods in the Calculus of Variations: {${L}^p$} Spaces}.
\newblock Springer {M}onographs in {M}athematics. Springer, New York, 2007.

\bibitem{GoffmanSerrin1964}
C.~Goffman and J.~Serrin.
\newblock Sublinear functions of measures and variational integrals.
\newblock {\em Duke Math. J.}, 31(1):159--178, Mar. 1964.

\bibitem{Grisvard2011}
P.~Grisvard.
\newblock {\em Elliptic Problems in Nonsmooth Domains}.
\newblock {Society for Industrial and Applied Mathematics}, Jan. 2011.

\bibitem{hardt1998penergy}
R.~Hardt, F.~Lin, and C.~Wang.
\newblock The {$p$}-energy minimality of {$x/|x|$}.
\newblock {\em Comm. Anal. Geom.}, 6(1):141--152, 1998.

\bibitem{HardtLin1987}
R.~Hardt and F.-H. Lin.
\newblock Mappings minimizing the {{$L^p$}} norm of the gradient.
\newblock {\em Comm. Pure Appl. Math.}, 40(5):555--588, 1987.

\bibitem{HarjulehtoHasto2019}
P.~Harjulehto and P.~H{\"a}st{\"o}.
\newblock {\em Orlicz {{Spaces}} and {{Generalized Orlicz Spaces}}}.
\newblock Springer International Publishing, Cham, 2019.

\bibitem{Jerrard1999}
R.~L. Jerrard.
\newblock Lower {{Bounds}} for {{Generalized Ginzburg--Landau Functionals}}.
\newblock {\em SIAM J. Math. Anal.}, 30(4):721--746, Jan. 1999.

\bibitem{lee2018riem}
J.~M. Lee.
\newblock {\em Introduction to {R}iemannian manifolds}, volume 176 of {\em Graduate Texts in Mathematics}.
\newblock Springer, Cham, 2018.

\bibitem{Leoni2017}
G.~Leoni.
\newblock {\em A First Course in {{Sobolev}} Spaces}.
\newblock Number 181 in Graduate Studies in Mathematics. AMS, Providence, 2nd ed edition, 2017.

\bibitem{lorentz1950somenew}
G.~G. Lorentz.
\newblock Some new functional spaces.
\newblock {\em Ann. of Math. (2)}, 51:37--55, 1950.

\bibitem{Maligranda1989}
L.~Maligranda.
\newblock {\em Orlicz Spaces and Interpolation}.
\newblock Number~5 in Seminars in {{Mathematics}}. Departamento de Matem{\'a}tica, Universidade Estadual de Campinas, 1989.

\bibitem{marcinkiewicz1939interpolation}
J.~Marcinkiewicz.
\newblock Sur l'interpolation d'op\'erations.
\newblock {\em CR Acad. Sci. Paris.}, 208:1272--1273, 1939.

\bibitem{MonteilEtAl2021}
A.~Monteil, R.~Rodiac, and J.~Van~Schaftingen.
\newblock Ginzburg\textendash{{Landau Relaxation}} for {{Harmonic Maps}} on {{Planar Domains}} into a {{General Compact Vacuum Manifold}}.
\newblock {\em Arch Rational Mech Anal}, 242(2):875--935, Nov. 2021.

\bibitem{MonteilEtAl2022}
A.~Monteil, R.~Rodiac, and J.~Van~Schaftingen.
\newblock Renormalised energies and renormalisable singular harmonic maps into a compact manifold on planar domains.
\newblock {\em Math. Ann.}, 383(3):1061--1125, Aug. 2022.

\bibitem{nash1956imbedding}
J.~Nash.
\newblock The imbedding problem for {R}iemannian manifolds, 1956.

\bibitem{palmer2024lifting}
D.~Palmer, A.~Chern, and J.~Solomon.
\newblock Lifting directional fields to minimal sections.
\newblock {\em ACM Transactions on Graphics (TOG)}, 43(4):1--20, 2024.

\bibitem{petrache2017controlled}
M.~Petrache and J.~Van~Schaftingen.
\newblock Controlled singular extension of critical trace {S}obolev maps from spheres to compact manifolds.
\newblock {\em Int. Math. Res. Not. IMRN}, (12):3647--3683, 2017.

\bibitem{polya1951isoperimetric}
G.~P{\'o}lya and G.~Szeg{\"o}.
\newblock {\em Isoperimetric {{Inequalities}} in {{Mathematical Physics}}}.
\newblock Number~27 in Annals of {{Mathematics Studies}}. Princeton University Press, Princeton, NJ, 1951.

\bibitem{rindler2018calculus}
F.~Rindler.
\newblock {\em Calculus of variations}.
\newblock Universitext. Springer Cham, 2018.

\bibitem{Rockafellar1997}
R.~T. Rockafellar.
\newblock {\em Convex analysis}, volume No. 28 of {\em Princeton Mathematical Series}.
\newblock Princeton University Press, Princeton, NJ, 1970.

\bibitem{Sandier1998}
E.~Sandier.
\newblock Lower {{Bounds}} for the {{Energy}} of {{Unit Vector Fields}} and {{Applications}}.
\newblock {\em J. Func. Anal.}, 152(2):379--403, Feb. 1998.

\bibitem{SchoenUhlenbeck1983}
R.~Schoen and K.~Uhlenbeck.
\newblock Boundary regularity and the {{Dirichlet}} problem for harmonic maps.
\newblock {\em J. Diff. Geom.}, 18(2):253--268, Jan. 1983.

\bibitem{solci2024nonlocal}
M.~Solci.
\newblock Nonlocal-interaction vortices.
\newblock {\em SIAM J. Math. Anal.}, 56(3):3430--3451, 2024.

\bibitem{van2024extension}
J.~Van~Schaftingen.
\newblock The extension of traces for {S}obolev mappings between manifolds.
\newblock {\em arXiv preprint arXiv:2403.18738}, 2024.

\bibitem{vanschaftingen2023asymptotic}
J.~Van~Schaftingen and B.~Van~Vaerenbergh.
\newblock Asymptotic behavior of minimizing {$p$}-harmonic maps when {$p \nearrow 2$} in dimension 2.
\newblock {\em Calc. Var. Partial Differential Equations}, 62(8):Paper No. 229, 45, 2023.

\bibitem{willem2013functional}
M.~Willem.
\newblock {\em Functional analysis}.
\newblock Cornerstones. Birkh\"{a}user/Springer, New York, 2013.
\newblock Fundamentals and applications.

\end{thebibliography}

\end{document}